\documentclass[12pt, a4paper]{amsart}
\usepackage{amsthm}
\theoremstyle{}
{\theoremstyle{definition}
\newtheorem{dfn}{Definition}[section]}
\newtheorem{prop}[dfn]{Proposition}

\newtheorem{thm}[dfn]{Theorem}
\newtheorem{rem}[dfn]{Remark}
\newtheorem{lem}[dfn]{Lemma}
\newtheorem{cor}[dfn]{Corollary}
\newtheorem{conj}[dfn]{Conjecture}
\newtheorem{exa}[dfn]{Example}
\usepackage[top=30truemm,bottom=30truemm,left=17truemm,right=17truemm]{geometry}
\usepackage{amsmath,amssymb,amscd,bm,graphicx}
\usepackage[all]{xy}

\renewcommand{\lim}{\varprojlim}
\begin{document}

\title[Equivalences of   derived  factorization categories of gauged LG models]{Equivalences of   derived  factorization categories of gauged Landau-Ginzburg models}
\author[Y.~Hirano]{Yuki Hirano}
\date{}

\maketitle{}

\begin{abstract}
For a given Fourier-Mukai equivalence of bounded derived categories of coherent sheaves on smooth quasi-projective varieties, we construct Fourier-Mukai equivalences of  derived factorization categories of gauged Landau-Ginzburg (LG) models.

 As an application, we obtain some equivalences of  derived factorization categories of K-equivalent gauged LG models. This result is an equivariant version of the result of Baranovsky and Pecharich, and  it also gives a partial answer to Segal's conjecture. As another application, we prove that if the  kernel of the Fourier-Mukai equivalence  is linearizable with respect to  a reductive affine algebraic group action, then the derived categories of equivariant coherent sheaves on  the varieties are equivalent. This result is shown by Ploog for finite groups case.
\end{abstract} 
\tableofcontents

\section{Introduction}
\subsection{Background and motivation}
A \emph{gauged Landau-Ginzburg (LG) model}  is  data $(X,W)^G$ consisting  of a smooth variety (or an algebraic stack) $X$  with a group $G$ action  and a semi invariant regular function $W$ on $X$. For a gauged LG model $(X,W)^G$, we consider a triangulated category $${\rm Dcoh}_G(X,W)$$ which is  called  \emph{the derived factorization category of the gauged LG model} $(X,W)^G$. If the group $G$ is trivial, omitting the subscript $G$ from the notation, we  denote by $(X,W)$ (resp. ${\rm Dcoh}(X,W)$) the gauged LG model (resp. its derived factorization category), and call it  \emph{the LG model}.  

Derived factorization categories of gauged LG models play  an important role in  Homological Mirror Symmetry for non-Calabi-Yau varieties \cite{orlov},  and are useful to study derived categories of coherent sheaves on algebraic stacks. For example, studying  windows in derived factorization categories gives a new technique to obtain some equivalences or semi-orthogonal decompositions of derived categories of algebraic stacks \cite{segal}, \cite{vgit}. 

The triangulated category ${\rm Dcoh}_G(X,W)$ is a generalization of the bounded derived category of coherent sheaves on $X$. In fact, for a gauged LG model $(X,0)^{\mathbb{G}_m}$ with trivial $\mathbb{G}_m$-action,  its derived factorization category is equivalent to the bounded derived category of coherent sheaves on $X$, namely $${\rm Dcoh}_{\mathbb{G}_m}(X,0)\cong{\rm D^b}({\rm coh}X).$$
Hence, it is natural to expect similarities between derived categories and derived factorization categories;  such similarities are observed in  \cite{velez}, \cite{bp},  \cite{ls}, for example.   In the present paper, we obtain  equivalences between derived factorization categories of certain gauged LG models from  equivalences between derived  categories of smooth quasi-projective varieties. 

\subsection{Main results}
Let $X_1$ and $X_2$ be smooth quasi-projective varieties over an algebraically closed field $k$ of characteristic zero, and  $G$ be a reductive  affine algebraic group acting on each $X_i$. Let $W_i:X_i\rightarrow \mathbb{A}^1$ be a $\chi$-semi invariant regular function on $X_i$ for some character $\chi:G\rightarrow\mathbb{G}_m$, and  $\pi_i:X_1\times X_2\rightarrow X_i$ be the projections. Consider the fibre product 
\[\xymatrix{
&&X_1\times_{\mathbb{A}^1} X_2\ar[lld]_{}\ar[rrd]^{}&&\\
X_1\ar[rrd]^{W_1}&&&&X_2\ar[lld]_{W_2}\\
&&\mathbb{A}^1&&
}\]
and let $j:X_1\times_{\mathbb{A}^1} X_2\hookrightarrow X_1\times X_2$ be the embedding.

 An object $P\in{\rm D^b}({\rm coh}X_1\times_{\mathbb{A}^1} X_2)$ whose support is proper over $X_2$ defines the integral functor 
$$\Phi_{j_*(P)}:{\rm D^b}({\rm coh}X_1)\rightarrow{\rm D^b}({\rm coh}X_2)\qquad (-)\longmapsto{\bf R}\pi_{2*}(\pi_1^*(-)\otimes^{\bf L}j_*(P)).$$
 On the other hand, the object $P$ induces an object $\widetilde{P}\in{\rm Dcoh}(X_1\times X_2,\pi_2^*W_2-\pi_1^*W_1^*)$ and it defines the integral functor 
$$\Phi_{\widetilde{P}}:{\rm Dcoh}(X_1,W_1)\rightarrow{\rm Dcoh}(X_2,W_2)\qquad(-)\longmapsto{\bf R}\pi_{2*}(\pi_1^*(-)\otimes^{\bf L}\widetilde{P}).$$
 Furthermore, if the object  $P$ is $G$-linearizable, i.e. it is in the essential image of the forgetful functor $$\Pi:{\rm D^b}({\rm coh}_GX_1\times_{\mathbb{A}^1} X_2)\rightarrow{\rm D^b}({\rm coh}X_1\times_{\mathbb{A}^1} X_2),$$ then the object $P$ induces an object $\widetilde{P}_G\in {\rm Dcoh}_G(X_1\times X_2,\pi_2^*W_2-\pi_1^*W_1^*)$ and it defines the integral functor 
$$\Phi_{\widetilde{P}_G}:{\rm Dcoh}_G(X_1,W_1)\rightarrow{\rm Dcoh}_G(X_2,W_2)\qquad(-)\longmapsto{\bf R}\pi_{2*}(\pi_1^*(-)\otimes^{\bf L}\widetilde{P}_G).$$

The main result of the present paper is the following:

\begin{thm}[Theorem \ref{main theorem 2}]\label{1.1}Let $P\in{\rm D^b}({\rm coh}X_1\times_{\mathbb{A}^1} X_2)$ be a $G$-linearizable object whose support is proper over $X_1$ and $X_2$.  If the integral functor  $\Phi_{j_*(P)}:{\rm D^b}({\rm coh}X_1)\rightarrow{\rm D^b}({\rm coh}X_2)$ is an equivalence (resp. fully faithful), then so is $\Phi_{\widetilde{P}_G}:{\rm Dcoh}_G(X_1,W_1)\rightarrow{\rm Dcoh}_G(X_2,W_2)$.
\end{thm}

This theorem is proved when the group $G$ is trivial, the functions $W_i$ are flat, and $X_i$ are smooth Deligne-Mumford stacks, in \cite{bp}.
Combining Theorem \ref{1.1} with  the result in \cite{bridgeland},  we obtain the following:

\begin{cor}\label{cor1.2}
Let $X$ and $X^{+}$ be smooth quasi-projective threefolds, and let the diagram $$X\xrightarrow{f} Y\xleftarrow{f^+} X^{+}$$ be a flop.
Let $G$ be a reductive affine algebraic group acting on $X$, $X^{+}$ and $Y$ with the morphisms $f$ and $f^+$ equivariant. Take a semi invariant regular function $W_Y:Y\rightarrow\mathbb{A}^1$, and set $W:=f^*W_Y$ and $W^{{+}}:=f^{+*}W_Y$. Then we have an equivalence
$${\rm Dcoh}_G(X,W)\cong{\rm Dcoh}_G(X^{+},W^{+}).$$ 
\end{cor}
The gauged LG models $(X,W)^G$ and $(X^{+},W^{+})^G$ in Corollary \ref{cor1.2} are \emph{K-equivalent}. Here, {\it K}-equivalence means that there  exists a common equivariant resolution of the varieties such that the pull-backs of the functions of LG models, and the classes of  canonical divisors, coincide. We expect the following conjecture, which is a generalization of \cite[Conjecture 2.15]{segal}: 

\begin{conj}\label{conj1.3}
If two gauged LG models are K-equivalent,  then their derived factorization categories are equivalent.
\end{conj}
Conjecture \ref{conj1.3}  for gauged LG models with trivial $\mathbb{G}_m$-actions and trivial functions is proposed in \cite{kawamata}.

As another corollary of Theorem \ref{1.1}, we obtain the following result.
\begin{cor}\label{cor1.4}
Let $P\in{\rm D^b}({\rm coh}X_1\times X_2)$ be an  $G$-linearizable object whose support is proper over $X_1$ and $X_2$. Let $P_G\in{\rm D^b}({\rm coh}_GX_1\times X_2)$ be an object with $\Pi(P_G)\cong P$, where $\Pi$ is the forgetful functor.  If the integral functor $\Phi_P:{\rm D^b}({\rm coh}X_1)\rightarrow{\rm D^b}({\rm coh}X_2)$ is an equivalence (resp.~fully faithful), then so is $\Phi_{P_G}:{\rm D^b(coh}_GX_1)\rightarrow{\rm D^b(coh}_GX_2)$.
 \end{cor}
 
Corollary \ref{cor1.4} is obtained in the case of smooth projective varieties  with finite group actions by \cite[Lemma 5]{ploog}; see also \cite{ks}. We  can also prove Corollary \ref{cor1.4} for a finite group $G$ by using \cite[Theorem 5.2]{elagin2}.

\subsection{Sketch of the proof of Theorem 1.1}
We divide the proof into two steps. The first step is to prove the equivalence of derived factorization categories of LG models without $G$-actions. In the second step, we equivariantize the equivalence obtained in the first step by considering categories of comodules over comonads. The idea of using comodules over comonads  comes  from Elagin's works \cite{elagin1}, \cite{elagin2}.

\textbf{1st step:} In the first step, using the technique in the proof of the  main result of \cite{bdfik}, we show that if $$\Phi_P:{\rm D^b}({\rm coh}X_1)\rightarrow{\rm D^b}({\rm coh}X_2)$$ is fully faithful, then the induced integral functor $$\Phi_{\widetilde{P}}:{\rm Dcoh}(X_1,W_1)\rightarrow{\rm Dcoh}(X_2,W_2)$$ is also fully faithful. 

We have an object $P_L\in{\rm D^b}({\rm coh}X_1\times_{\mathbb{A}^1} X_2)$ such that its support is proper over $X_1$ and the integral functor 
$$\Phi_{j_*(P_L)}:{\rm D^b}({\rm coh}X_2)\rightarrow{\rm D^b}({\rm coh}X_1)$$
is  left adjoint  to $\Phi_{j_*(P)}$.
  Furthermore, the object $P_L$  induces an object $\widetilde{P_L}\in{\rm Dcoh}(X_1\times X_2, \pi_1^*W_1-\pi_2^*W_2)$ with the integral functor $$\Phi_{\widetilde{P_L}}:{\rm Dcoh}(X_2,W_2)\rightarrow{\rm Dcoh}(X_1,W_1)$$  left adjoint to $\Phi_{\widetilde{P}}$. 
Note that the functor $\Phi_P$ is fully faithful if and only if the adjunction $\Phi_{P_L}\circ\Phi_P\rightarrow{\rm id}$ is an isomorphism. By the argument in the proof of the main result of \cite{bdfik}, we see that  the isomorphism $\Phi_{P_L}\circ\Phi_P\xrightarrow{\sim}{\rm id}$  implies an isomorphism of functors  $\Phi_{\widetilde{P_L}}\circ\Phi_{\widetilde{P}}\cong{\rm id}$. Since $\Phi_{\widetilde{P_L}}\dashv\Phi_{\widetilde{P}}$, we see that $\Phi_{\widetilde{P}}$ is fully faithful. 
If  $\Phi_{j_*(P)}$ is an equivalence, its left adjoint $\Phi_{j_*(P_L)}$ is fully faithful. By the above argument, we see that $\Phi_{\widetilde{P_L}}$ is also fully faithful, which implies the equivalence of $\Phi_{\widetilde{P}}$.

\textbf{2nd step:} We prove that if  $$\Phi_{\widetilde{P}}:{\rm Dcoh}(X_1,W_1)\rightarrow{\rm Dcoh}(X_2,W_2)$$ is an equivalence (resp. fully faithful), then so is $$\Phi_{\widetilde{P}_G}:{\rm Dcoh}_G(X_1,W_1)\rightarrow{\rm Dcoh}_G(X_2,W_2).$$ 

The key tools are categories of  comodules over comonads (see section 2 for the definitions).
A $G$-action on $X_i$ induces a comonad $\mathbb{T}_i$ on the category ${\rm DQcoh}(X_i,W_i)$.  Consider the category, ${\rm DQcoh}(X_i,W_i)_{\mathbb{T}_i}$, of comodules over the comonad $\mathbb{T}_i$. Then,  there are natural functors $$\Gamma_i:{\rm DQcoh}_G(X_i,W_i)\rightarrow{\rm DQcoh}(X_i,W_i)_{\mathbb{T}_i},$$  which are fully faithful, since $G$ is reductive, and we have an induced functor $$(\Phi_{\widetilde{P}})_{\mathbb{T}}:{\rm DQcoh}(X_1,W_1)_{\mathbb{T}_1}\rightarrow{\rm DQcoh}(X_2,W_2)_{\mathbb{T}_2}$$ between categories of comodules  such that  the following diagram is commutative:
\[\xymatrix{
{\rm DQcoh}(X_1,W_1)_{\mathbb{T}_1}^{\rm coh}\ar[rr]^{(\Phi_{\widetilde{P}})_{\mathbb{T}}}&&{\rm DQcoh}(X_2,W_2)_{\mathbb{T}_2}^{\rm coh}\\
{\rm Dcoh}_G(X_1,W_1)\ar[rr]^{\Phi_{\widetilde{P}_G}}\ar@{^{(}-{>}}[u]^{\Gamma_1}&&{\rm Dcoh}_G(X_2,W_2)\ar@{^{(}-{>}}[u]_{\Gamma_2}
}\]
where ${\rm DQcoh}(X_i,W_i)_{\mathbb{T}_i}^{\rm coh}$ is a certain full subcategory of ${\rm DQcoh}(X_i,W_i)_{\mathbb{T}_i}$, defined in section 4.3.
By the full-faithfulness of $\Phi_{\widetilde{P}}$ and the existence of the left adjoint functor $\Phi_{\widetilde{P_L}}$, we see that the functor $(\Phi_{\widetilde{P}})_{\mathbb{T}}:{\rm DQcoh}(X_1,W_1)_{\mathbb{T}_1}^{\rm coh}\rightarrow{\rm DQcoh}(X_2,W_2)_{\mathbb{T}_2}^{\rm coh}$ is fully faithful, whence $\Phi_{\widetilde{P}_G}$ is also fully faithful. If $\Phi_{\widetilde{P}}$ is an equivalence, we see that the composition $(\Phi_{\widetilde{P}})_{\mathbb{T}}\circ(\Phi_{\widetilde{P_L}})_{\mathbb{T}}$  is isomorphic to the identity functor on  ${\rm DQcoh}(X_2,W_2)_{\mathbb{T}_2}^{\rm coh}$. This means that the composition $\Phi_{\widetilde{P}_G}\circ\Phi_{(\widetilde{P_L})_G}$ is isomorphic to the identity functor on ${\rm Dcoh}_G(X_2,W_2)$. Hence,  $\Phi_{\widetilde{P}_G}$ is essentially surjective, which implies its equivalence.

\subsection{Plan of the paper}In section 2, we recall  the foundations of  the theory of comodules over comonads. In section 3 and 4, we introduce definitions and basic properties of  derived factorization categories and functors between them. In section 5, we prove the main theorem  and give applications of it.

\subsection{Notation and convention}

For a functor $F:\mathcal{A}\rightarrow\mathcal{B}$ from a category $\mathcal{A}$ to a category $\mathcal{B}$, and a full subcategory $\mathcal{C}\subset\mathcal{A}$, we denote by $F(\mathcal{C})$ the full subcategory of $\mathcal{B}$ whose objects are of the form $F(C)$ for some $C\in\mathcal{C}$, and denote by Im($F|_{\mathcal{C}}$) the essential image of $F|_{\mathcal{C}}$, i.e. the full subcategory of $\mathcal{B}$ whose objects are isomorphic to objects of $F(\mathcal{C})$.

Let $F, G:\mathcal{A}\rightarrow\mathcal{B}$ be functors and let $\alpha:F\rightarrow G$ be a functor morphism. For a functor $P:\mathcal{C}\rightarrow\mathcal{A}$, we denote by $\alpha P:F\circ P\rightarrow G\circ P$ the functor morphism  given by $\alpha P(C):=\alpha(P(C))$ for any object $C\in\mathcal{C}$. For a functor $Q:\mathcal{B}\rightarrow\mathcal{D}$, we denote by $Q\alpha:Q\circ F\rightarrow Q\circ G$ the functor morphism defined by $Q\alpha(A):=Q(\alpha(A))$ for any object $A\in\mathcal{A}$.

For any exact category $\mathcal{E}$, ${\rm Ch^b}(\mathcal{E})$ denotes the category of bounded complexes in $\mathcal{E}$, ${\rm K^b}(\mathcal{E})$ denotes the bounded homotopy category of $\mathcal{E}$, and ${\rm D^b}(\mathcal{E})$ denotes the bounded derived category of $\mathcal{E}$.

\subsection{Acknowledgments}
The author would like to express his deep gratitude to his supervisor Hokuto Uehara for his valuable advice 
and many suggestions to improve this paper.
The author is a Research Fellow of Japan Society for the Promotion of Science. He is partially supported by Grand-in-Aid for JSPS Fellows $\sharp$26-6240.


\section{Comodules over comonads}

Categories of comodules over comonads are main tools to obtain the main result.  In this section, we recall the definitions of comonads and comodules over comonads, and provide basic properties of them, following \cite{elagin2}.

\subsection{Comodules over comonads}

Let $\mathcal{C}$ be a category. We start by recalling the definitions of comonads  on $\mathcal{C}$ and comodules over a comonad.

\vspace*{2mm}
\begin{dfn}
A \textbf{comonad} $\mathbb{T}=(T,\varepsilon,\delta)$ on the category $\mathcal{C}$ consists of a functor $T:\mathcal{C}\rightarrow\mathcal{C}$ and functor morphisms $\varepsilon:T\rightarrow {\rm id}_{\mathcal{C}}$ and $\delta:T\rightarrow T^2$ such that the following diagrams are commutative$:$

\[\xymatrix{
&T \ar[r]^{\delta} \ar[d]_{\delta} \ar@{=}[dr]^{{\rm id}_T} & T^2 \ar[d]^{T\varepsilon}&&&T \ar[r]^{\delta} \ar[d]_{\delta} & T^2 \ar[d]^{T\delta}\\
&T^2 \ar[r]^{\varepsilon T} & T && &T^2 \ar[r]^{\delta T} & T^3
}\]

\end{dfn}

\vspace*{2mm}
\begin{exa}\label{adj}
Let $P=(P^*\dashv P_*)$ be an adjoint pair of functors $P^*:\mathcal{C}\rightarrow\mathcal{D}$ and $P_*:\mathcal{D}\rightarrow\mathcal{C}$, and let $\eta_{P} :{\rm id}_{\mathcal{C}}\rightarrow P_*P^*$ and $\varepsilon_{P}:P^*P_*\rightarrow {\rm id}_{\mathcal{D}}$ be the adjunction morphisms. Set $T_P:=P^*P_*$ and $\delta_P:=P^*\eta_P P_*$. Then   $\mathbb{T}(P):=(T_P,\varepsilon_{P},\delta_P)$ is a comonad on $\mathcal{D}$.
\end{exa}

\vspace*{2mm}
\begin{dfn} \label{comodule}
Let $\mathbb{T}=(T,\varepsilon,\delta)$ be a comonad on $\mathcal{C}$. A \textbf{comodule} over $\mathbb{T}$ is a pair $(C,\theta_C)$ of an object $C\in\mathcal{C}$ and a morphism $\theta_{C}:C\rightarrow T(C)$ such that 

$(1)$ $\varepsilon(C)\circ\theta_C={\rm id}_C$

$(2)$  the following diagram is commutative:

$$\begin{CD}
C@>{\theta_C}>>T(C)\\
@V{\theta_C}VV @VVT({\theta_C})V\\
T(C)@>{\delta(C)}>>T^2(C).
\end{CD}$$

\end{dfn}

\vspace*{3mm}
Given a comonad $\mathbb{T}$ on $\mathcal{C}$, we define  the category $\mathcal{C}_{\mathbb{T}}$ of comodules over the comonad $\mathbb{T}$ as follows:


\begin{dfn}
Let $\mathbb{T}=(T,\varepsilon,\delta)$ be a comonad on $\mathcal{C}$.  The category $\mathcal{C}_{\mathbb{T}}$ of comodules over $\mathbb{T}$ on $\mathcal{C}$ is the category whose  objects are comodules over $\mathbb{T}$ and  whose sets of morphisms  are defined as follows;
$${\rm Hom}((C_1,\theta_{C_1}),(C_2,\theta_{C_2})):=\{f:{C_1}\rightarrow {C_2} \mid  T(f)\circ\theta_{C_1}=\theta_{C_2}\circ f \}.$$
For a full subcategory $\mathcal{B}\subset\mathcal{C}$, we define the full subcategory  $\mathcal{C}_{\mathbb{T}}^{\mathcal{B}}\subset\mathcal{C}_{\mathbb{T}}$ as
$$\textrm{Ob}(\mathcal{C}_{\mathbb{T}}^{\mathcal{B}}):=\{(C,\theta_C)\in\textrm{Ob}(\mathcal{C}_{\mathbb{T}})\mid C\cong B  \textrm{ for some } B\in\mathcal{B}\}.$$

\end{dfn}

\vspace*{2mm}
\begin{rem}\label{rem2.5}
Let $(C,\theta_C)\in\mathcal{C}_{\mathbb{T}}^{\mathcal{B}}$. By definition, there exist an object $B\in\mathcal{B}$ and an isomorphism $\varphi : C\xrightarrow{\sim} B$. If we set $\theta_B:=T(\varphi)\theta_C\varphi^{-1}$, then the pair $(B,\theta_B)$ is an object of $\mathcal{C}_{\mathbb{T}}^{\mathcal{B}}$ and $\varphi$ gives an isomorphism from $(C,\theta_C)$ to $(B,\theta_B)$  in $\mathcal{C}_{\mathbb{T}}^{\mathcal{B}}$.

\end{rem}

\vspace{3mm} For a comonad which is given by an adjoint pair $(P^*\dashv P_*)$, we have a canonical functor, called comparison functor, from  the domain of $P^*$ to the category of comodules over the comonad.
\vspace*{2mm}
\begin{dfn}\label{comparison functor}
The notation is the same as in  Example \ref{adj}. For an adjoint pair $P=(P^*\dashv P_*)$, we define a  functor $$\Gamma_{P}:\mathcal{C}\rightarrow\mathcal{D}_{\mathbb{T}(P)}$$ as follows: For any $C\in\mathcal{C}$ and for any  morphism $f$ in $\mathcal{C}$
$$\Gamma_{P}(C):=(P^*(C),P^*(\eta_{P}(C))) \hspace{3mm}\textrm{and}\hspace{3mm}\Gamma_{P}(f):=P^*(f).$$
This functor is called \textbf{the comparison functor} of $P$. Restricting $\Gamma_{P}$ to a full subcategory $\mathcal{B}\subset\mathcal{C}$, we have a restricted functor $$\Gamma_{P}|_{\mathcal{B}}:\mathcal{B}\rightarrow\mathcal{D}_{\mathbb{T}(P)}^{P^*(\mathcal{B})}.$$\end{dfn}

The following proposition gives  sufficient conditions for a comparison functor to be fully faithful or an equivalence.

\vspace*{2mm}
\begin{prop}[\cite{elagin1} Theorem 3.9, Corollary 3.11]\label{comparison theorem}
The notation is the same as in  {\rm Example} {\rm\ref{adj}}.
\begin{itemize}

\item[$(1)$] If for any $C\in\mathcal{C}$, the morphism $\eta_P(C):C\rightarrow P_*P^*(C)$ is a split mono, i.e. there is a  morphism $\zeta_C:P_*P^*(C)\rightarrow C$ such that $\zeta\circ\eta_P(C)={\rm id}_{{C}}$, then the comparison functor $\Gamma_{P}:\mathcal{C}\rightarrow\mathcal{D}_{\mathbb{T}(P)}$ is fully faithful. 
\item[$(2)$] If  $\mathcal{C}$ is idempotent complete and the functor morphism $\eta_P :{\rm id}_{\mathcal{C}}\rightarrow P_*P^*$ is split mono, i.e. there exists a functor morphism $\zeta:P_*P^*\rightarrow {\rm id}_{\mathcal{C}}$ such that $\zeta\circ\eta={\rm id}$, then $\Gamma_{P}:\mathcal{C}\rightarrow\mathcal{D}_{\mathbb{T}(P)}$ is an equivalence.

\end{itemize}
\end{prop}

\vspace*{3mm}
\subsection{Functors between categories of comodules} 
We introduce the notion of linearizable functors which induce natural functors between categories of comodules.
Let $\mathcal{A}$ $($resp. $\mathcal{B}$ and $\mathcal{C})$ be a category and let $\mathbb{T}_{\mathcal{A}}=(T_{\mathcal{A}},\varepsilon_{\mathcal{A}},\delta_{\mathcal{A}})$ $($resp. $\mathbb{T}_\mathcal{B}=(T_{\mathcal{B}},\varepsilon_{\mathcal{B}},\delta_{\mathcal{B}})$ and $\mathbb{T}_\mathcal{C}=(T_{\mathcal{C}},\varepsilon_{\mathcal{C}},\delta_{\mathcal{C}}))$ be a comonad on $\mathcal{A}$ $($resp.   $\mathcal{B}$ and $\mathcal{C})$.

\vspace*{2mm}
\begin{dfn}\label{equivariant}A functor $F:\mathcal{A}\rightarrow\mathcal{B}$ is called \textbf{linearizable} with respect to $\mathbb{T}_{\mathcal{A}}$ and $\mathbb{T}_{\mathcal{B}}$, or just \textbf{linearizable}, if there exists an isomorphism of functors $$\Omega:F  T_{\mathcal{A}}\xrightarrow{\sim}T_{\mathcal{B}}  F$$ such that the following two diagrams of functor morphisms are commutative :
\[\xymatrix{
(1)& FT_{\mathcal{A}} \ar[rr]^{\Omega} \ar[rd]_{F\varepsilon_{\mathcal{A}}} & &  T_{\mathcal{B}}F \ar[ld]^{\varepsilon_{\mathcal{B}}F}&&(2)&FT_{\mathcal{A}} \ar[rrr]^{\Omega} \ar[d]_{F\delta_{\mathcal{A}}} &&&   T_{\mathcal{B}}F \ar[d]^{\delta_{\mathcal{B}}F} \\ && F &&&&FT^2_{\mathcal{A}} \ar[rrr]^{T_{\mathcal{B}}\Omega\circ\Omega T_{\mathcal{A}}} & && T^2_{\mathcal{B}}F  \\
}\]
\vspace*{2mm}We call  the pair $(F,\Omega)$  a \textbf{linearized functor} with respect to $\mathbb{T}_{\mathcal{A}}$ and $\mathbb{T}_{\mathcal{B}}$, and the isomorphism of functors $\Omega$ is called a \textbf{linearization} of $F$ with respect to $\mathbb{T}_{\mathcal{A}}$ and $\mathbb{T}_{\mathcal{B}}$.
\end{dfn}

\vspace{8mm}
If $F:\mathcal{A}\rightarrow\mathcal{B}$ is a linearizable functor with a linearization $\Omega:F T_{\mathcal{A}}\xrightarrow{\sim}T_{\mathcal{B}} F$, we have an induced functor 
$${F_{\Omega}}:\mathcal{A}_{\mathbb{T}_{\mathcal{A}}}\rightarrow\mathcal{B}_{\mathbb{T}_{\mathcal{B}}}$$
defined by $$F_{\Omega}(A,\theta_A):=(F(A),\Omega(A)\circ F(\theta_A))\hspace{3mm}\textrm{and}\hspace{3mm}F_{\Omega}(f):=F(f).$$

\vspace{2mm}
\begin{lem}\label{composition lemma}
Let $F:\mathcal{A}\rightarrow\mathcal{B}$ and $G:\mathcal{B}\rightarrow\mathcal{C}$ be linearizable functors with linearizations $\Phi:FT_{\mathcal{A}}\xrightarrow{\sim}T_{\mathcal{B}}F$ and  
$\Psi:GT_{\mathcal{B}}\xrightarrow{\sim}T_{\mathcal{C}}G$ respectively. Then the composition $GF$ is a linearizable functor with linearization $\Omega:=\Psi F \circ G \Phi$ and $(GF)_{\Omega}=G_{\Psi}F_{\Phi}$.

\begin{proof}
By definition it is sufficient to prove that $$ GF\varepsilon_{\mathcal{A}}=\varepsilon_{\mathcal{C}}GF\circ\Omega \hspace{3mm}\textrm{and}\hspace{3mm} T_{\mathcal{C}}\Omega \circ \Omega T_{\mathcal{A}}\circ GF \delta_{\mathcal{A}}=\delta_{\mathcal{C}}GF \circ \Omega.$$ 
The former one of the above equations follows from easy diagram chasing as follows.
\begin{align}
GF\varepsilon_{\mathcal{A}}&=G(\varepsilon_{\mathcal{B}}F\circ\Phi)=G\varepsilon_{\mathcal{B}}F\circ G\Phi=(\varepsilon_{C}G\circ\Psi)F\circ G\Phi=\varepsilon_{\mathcal{C}}GF\circ\Omega,\notag
\end{align}
where the first and third equations  follow from the commutativity of the diagrams corresponding to $(1)$ in Definition \ref{equivariant}.
The latter one is verified as follows;
\begin{align}
&\hspace{5mm}T_{\mathcal{C}}\Omega \circ \Omega T_{\mathcal{A}}\circ GF \delta_{\mathcal{A}}\notag\\
&= T_{\mathcal{C}}\Psi F\circ T_{\mathcal{C}}G\Phi \circ \Psi FT_{\mathcal{A}}\circ G\Phi T_{\mathcal{A}}\circ GF \delta_{\mathcal{A}}\notag\\
&=T_{\mathcal{C}}\Psi F\circ\Psi T_{\mathcal{B}}F\circ GT_{\mathcal{B}}\Phi\circ G\Phi T_{\mathcal{A}}\circ GF\delta_{\mathcal{A}}\notag\\
&=T_{\mathcal{C}}\Psi F\circ\Psi T_{\mathcal{B}}F\circ G(T_{\mathcal{B}}\Phi\circ\Phi T_{\mathcal{A}}\circ F\delta_{\mathcal{A}})\notag\\
&=T_{\mathcal{C}}\Psi F\circ\Psi T_{\mathcal{B}}F\circ G(\delta_{\mathcal{B}}F\circ\Phi)\notag\\
&=(T_{\mathcal{C}}\Psi\circ\Psi T_{\mathcal{B}}\circ G\delta_{\mathcal{B}})F\circ G\Phi\notag\\
&=(\delta_{\mathcal{C}}G\circ\Psi) F\circ G\Phi\notag\\
&=\delta_{\mathcal{C}}GF\circ\Omega,\notag
\end{align}
where the second equation follows from the functoriality of $\Psi$, and the fourth and the sixth equations  follow from the commutativity of the diagrams corresponding to $(2)$ in Definition \ref{equivariant}.
\end{proof}
\end{lem}

\vspace{2mm}
The next proposition gives a sufficient condition for a restriction of  the  functor $F_{\Omega}$ associated with a linearized functor $(F,\Omega)$ to be fully faithful or an equivalence.

\vspace*{2mm}
\begin{prop}\label{main prop}

Let $F:\mathcal{A}\rightarrow\mathcal{B}$ be a linearizable functor with a linearization $\Omega:F\circ T_{\mathcal{A}}\xrightarrow{\sim}T_{\mathcal{B}}\circ F$. Let $\mathcal{C}\subset\mathcal{A}$ be a full subcategory of $\mathcal{A}$ and let $\mathcal{D}\subset\mathcal{B}$ be a full subcategory of $\mathcal{B}$ containing $F(\mathcal{C})$.
Assume the following condition:

\begin{itemize}
\item[$(\ast)$:] 
${\rm Hom}(C,T^n_{\mathcal{A}}(C'))\xrightarrow{F}{\rm Hom}(F(C),F(T^n_{\mathcal{A}}(C')))$ is an isomorphism for any $C, C'\in\mathcal{C}$ and $n=1,2$.
\end{itemize}
 If $F|_{\mathcal{C}}:\mathcal{C}\rightarrow\mathcal{D}$ is fully faithful (resp. an equivalence), then the functor
$$F_{\Omega}|_{\mathcal{A}^{\mathcal{C}}_{\mathbb{T}_{\mathcal{A}}}}:\mathcal{A}^{\mathcal{C}}_{\mathbb{T}_{\mathcal{A}}}\rightarrow\mathcal{B}^{\mathcal{D}}_{\mathbb{T}_{\mathcal{B}}}$$
is also fully faithful (resp. an equivalence).

\begin{proof}Assume $F|_{\mathcal{C}}$ is fully faithful. 
At first we show that $F_{\Omega}$ is  fully faithful on $\mathcal{A}^{\mathcal{C}}_{\mathbb{T}_{\mathcal{A}}}$. 

Let $\widetilde{C}:=(C,\theta_C)$ and $\widetilde{C'}:=(C',\theta_{C'})$ be objects of $\mathcal{A}^{\mathcal{C}}_{\mathbb{T}_{\mathcal{A}}}$. By Remark \ref{rem2.5}, we may assume that $C$ and $C'$ are objects of $\mathcal{C}$. For $f,g\in{\rm Hom}(\widetilde{C},\widetilde{C'})\subset{\rm Hom}_{\mathcal{A}}({C},{C'})$, if $F_{\Omega}(f)=F_{\Omega}(g)$, then $F(f)=F(g)$ as morphisms in $\mathcal{B}$. Since $F$ is fully faithful on $\mathcal{C}$, this implies that $f=g$ as morphisms in $\mathcal{A}$, whence $f=g$ in $\mathcal{A}_{\mathbb{T}_{\mathcal{A}}}$. Hence $F_{\Omega}$ is faithful. 

Take any morphism $h\in{\rm Hom}(F_{\Omega}(\widetilde{C}),F_{\Omega}(\widetilde{C'}))$. Since $F$ is full on $\mathcal{C}$, there exists a morphism $f\in{\rm Hom}(C,C')$ such that $F(f)=h$, and we have the following commutative diagram:
$$\begin{CD}
F(C) @>F(\theta_C)>> F(T_{\mathcal{A}}(C)) @>\Omega(C)>>T_{\mathcal{B}}(F(C))\\
@VF(f)VV && @VVT_{\mathcal{B}}(F(f))V\\
F(C') @>F(\theta_{C'})>> F(T_{\mathcal{A}}(C')) @>\Omega(C')>>T_{\mathcal{B}}(F(C'))
\end{CD}$$
 By the functoriality of $\Omega$, the following diagram is commutative:
 $$\begin{CD}
 F(T_{\mathcal{A}}(C)) @>\Omega(C)>>T_{\mathcal{B}}(F(C))\\
 @VF(T_{\mathcal{A}}(f))VV  @VVT_{\mathcal{B}}(F(f))V\\
 F(T_{\mathcal{A}}(C')) @>\Omega(C')>>T_{\mathcal{B}}(F(C'))
 \end{CD}$$
Combining commutativity of the above diagrams, we have $$F(T_{\mathcal{A}}(f)\circ\theta_C)=F(\theta_C'\circ f)$$  since $\Omega(C')$ is an isomorphism.
 By the condition $(\ast)$ in the assumption, we see that $T_{\mathcal{A}}(f)\circ\theta_C=\theta_C'\circ f$, which implies that $f$ is a morphism in ${\mathcal{A}^{\mathcal{C}}_{\mathbb{T}_{\mathcal{A}}}}$. Hence $F_{\Omega}$ is full.
 
  Assume $F|_{\mathcal{C}}$ is an equivalence. We verify that the functor $F_{\Omega}|_{\mathcal{A}^{\mathcal{C}}_{\mathbb{T}_{\mathcal{A}}}}:\mathcal{A}^{\mathcal{C}}_{\mathbb{T}_{\mathcal{A}}}\rightarrow\mathcal{B}^{\mathcal{D}}_{\mathbb{T}_{\mathcal{B}}}$ is essentially surjective. 
Since $F|_{\mathcal{C}}$ is an equivalence, it is sufficient to prove that for any object $(B,\theta_B)\in\mathcal{B}^{\mathcal{D}}_{\mathbb{T}_{\mathcal{B}}}$ with $B=F(C)$ for some $C\in\mathcal{C}$, there exists an object $(C,\theta_{C})\in{\mathcal{A}^{\mathcal{C}}_{\mathbb{T}_{\mathcal{A}}}}$ such that $F_{\Omega}(C,\theta_C)=(B,\theta_B)$. By the condition $(\ast)$, we know that there exists a morphism $\theta_C:C\rightarrow T_{\mathcal{A}}(C)$ such that $F(\theta_C)=\Omega(C)^{-1}\circ\theta_{F(C)}:F(C)\rightarrow F(T_{\mathcal{A}}(C))$. To show that the pair $(C,\theta_C)$ is an object of $\mathcal{A}_{\mathbb{T}_{\mathcal{A}}}$, we check  two conditions in Definition \ref{comodule}. Considering the following commutative diagram;
 \[\xymatrix{
 F(C) \ar[rrrrrr]^{F(\theta_C)} \ar[rrrd]^{\theta_{F(C)}} \ar@{=}[rrrddd] & & & & & & F(T_{\mathcal{A}}(C))\ar[llld]_{\Omega(C)} \ar[lllddd]^{F(\varepsilon_{\mathcal{A}}(C))}\\
 & & & T_{\mathcal{B}}(F(C)) \ar[dd]^{\varepsilon_{\mathcal{B}}(F(C))} & & & \\
 & & \\
 & & & F(C) & & &
 }\]
 we see that $F(\varepsilon_{\mathcal{A}}(C)\circ\theta_C)={\rm id}_{F(C)}$. Since $F|_{\mathcal{C}}$ is fully faithful, we obtain
 $$\varepsilon_{\mathcal{A}}(C)\circ\theta_C={\rm id}_C,$$which is the first condition in Definition \ref{comodule}. By the following commutative diagram;
  
  \[\xymatrix{
 & & & & & F(T_{\mathcal{A}}(C)) \ar[ld]^{\Omega(C)} \ar[ddd]_{F(T_{\mathcal{A}}(\theta_C))}&\\
  & & F(C) \ar[rr]_{\theta_{F(C)}} \ar[d]^{\theta_{F(C)}} \ar[urrr]^{F(\theta_C)} \ar[ldd]_{F(\theta_{C})}& & T_{\mathcal{B}}(F(C)) \ar[d]_{T_{\mathcal{B}}(\theta_{F(C)})} & & &\\
  & & T_{\mathcal{B}}(F(C)) \ar[rr]_{\delta_{\mathcal{B}}(F(C))} & & T^2_{\mathcal{B}}(F(C)) & & &\\
 &  F(T_{\mathcal{A}}(C)) \ar[ur]_{\Omega(C)} \ar[rrrr]_{F(\delta_{\mathcal{A}}(C))}  & & & & F( T^2_{\mathcal{A}}(C)) \ar[ul]^{T_{\mathcal{B}}(\Omega(C))\circ\Omega(T_{\mathcal{A}}(C))},&
 }\]
 we see that $F(\delta_{\mathcal{A}}(C)\circ\theta_C)=F(T_{\mathcal{A}}(\theta_C)\circ\theta_C).$
 By the condition $(\ast)$, we obtain $$\delta_{\mathcal{A}}(C)\circ\theta_C=T_{\mathcal{A}}(\theta_C)\circ\theta_C,$$ which is the second condition in Definition \ref{comodule}. Hence, the pair $(C,\theta_C)$ is a comodule over $\mathbb{T}_{\mathcal{A}}$, and we see that $F_{\Omega}(C,\theta_C)=(F(C),\theta_{F(C)})$ by the construction of $(C,\theta_C)$. 
\end{proof}

\end{prop}

\vspace{2mm}
The following lemma gives a useful criteria for a functor to be linearizable with respect to comonads which are constructed from $^{\rotatebox[origin=C]{180}{"}}$compatible'' adjoint pairs.


\vspace*{2mm}
\begin{lem}\label{compatible lemma}

Assume that we have the following diagram of functors;
 \[\xymatrix{
\mathcal{A}\ar@/^/[dd]^{P_*} \ar[rr]^{F} &&  \mathcal{B} \ar@/^/[dd]^{Q_*} \\\\
\mathcal{A}' \ar@/^/[uu]^{P^*} \ar[rr]^{F'}&& \mathcal{B'} \ar@/^/[uu]^{Q^*}
 }\]
 where $P:=(P^*\dashv P_*)$ and $Q:=(Q^*\dashv Q_*)$ are adjoint pairs.  Assume that we have two isomorphisms of functors $\Omega^*:F P^*\xrightarrow{\sim}Q^*F'$ and $\Omega_*:F' P_*\xrightarrow{\sim}Q_* F$. Let $\Omega:FT_P\xrightarrow{\sim}T_{Q}F$ be  the composition  of functor morphisms $Q^*\Omega_*\circ\Omega^*P_*$. Consider the following two diagrams of functor morphisms:
 \[\xymatrix{
({\rm i}) &Q^*F'P_* \ar[r]^{Q^*\Omega_*} & Q^*Q_*F \ar[d]^{\varepsilon_{Q}F}&&({\rm ii}) &F' \ar[r]^{\eta_{Q}F'} \ar[d]_{F'\eta_P} & Q_*Q^*F'  \\
&FP^*P_* \ar[u]^{\Omega^*P_*} \ar[r]^{F\varepsilon_{P}} & F &&&F'P_*P^* \ar[r]^{\Omega_*P^*} & Q_*FP^* \ar[u]_{Q_*\Omega^*} }\]
If the above two diagrams are commutative, then $(F,\Omega)$ is a linearized functor with respect to $\mathbb{T}(P)$ and $\mathbb{T}(Q)$, and  there exists an isomorphism of functors $\Sigma:F_{\Omega}\Gamma_{P}\xrightarrow{\sim}\Gamma_{Q}F'$.  
 
 \begin{proof}
 
 We verify that the diagrams corresponding to ones in Definition $\ref{equivariant}$ are commutative. The commutativity of (i) immediately implies the commutativity of the diagram corresponding  to (1) in Definition $\ref{equivariant}$. We show that if the diagram of (ii) is commutative,  the diagram corresponding to (2) in Definition $\ref{equivariant}$ is commutative. By the functoriality of $\Omega^*$ and $\eta_{Q}$,  the following diagrams of functor morphisms are commutative;
 
 \[\xymatrix{
 FP^*P_*\ar[rr]^{\Omega^*P_*} \ar[d]_{FP^*\eta_{P}P_*} && Q^*F'P_* \ar[d]^{Q^*F'\eta_{P}P_*} \\
 FP^*P_*P^*P_* \ar[rr]^{\Omega^*P_*P^*P_*} && Q^*F'P_*P^*P_*
 }\]
 and 
 \[\xymatrix{
 F'P_* \ar[rr]^{\eta_{P'}F'P_*} \ar[d]_{\Omega_*} &&Q_*Q^*F'P_* \ar[d]^{Q_*Q^*\Omega_*} \\
 Q_*F(B) \ar[rr]^{\eta_{Q}Q_*F} && Q_*Q^*Q_*F.
 }\]
 Hence, we have equations of  functor morphisms $$(a):\hspace{5mm}\Omega^*P_*P^*P_*\circ FP^*\eta_{P}P_*=Q^*F'\eta_{P}P_*\circ\Omega^*P_*$$ and 
 $$(b):\hspace{5mm}\eta_{Q}Q_*F\circ\Omega_*=Q_*Q^*\Omega_*\circ\eta_{Q}F'P_*.$$
 We see that the diagram corresponding to  (2) in Definition $\ref{equivariant}$  is commutative as follows;
 
\begin{align}
&\hspace{5mm}T_{Q}\Omega\circ\Omega T_P\circ F\delta_P\notag\\
&=T_{Q}(Q^*\Omega_*\circ\Omega^*P_*)\circ(Q^*\Omega_*\circ\Omega^*P_*)T_P\circ F\delta_P\notag\\
&=Q^*Q_*(Q^*\Omega_*\circ\Omega^*P_*)\circ(Q^*\Omega_*\circ\Omega^*P_*)P^*P_*\circ FP^*\eta_P P_*\notag\\
&=Q^*Q_*Q^*\Omega_*\circ Q^*Q_*\Omega^*P_*\circ Q^*\Omega_*P^*P_*\circ(\Omega^*P_*P^*P_*\circ FP^*\eta_P P_*)\notag\\
(a)\rightarrow&=Q^*Q_*Q^*\Omega_*\circ Q^*Q_*\Omega^*P_*\circ Q^*\Omega_*P^*P_*\circ(Q^*F'\eta_{P}P_*\circ\Omega^*P_*)\notag\\
&=Q^*Q_*Q^*\Omega_*\circ (Q^*Q_*\Omega^*P_*\circ Q^*\Omega_*P^*P_*\circ Q^*F'\eta_{P}P_*)\circ\Omega^*P_*\notag\\
&=Q^*Q_*Q^*\Omega_*\circ Q^*(Q_*\Omega^*\circ \Omega_*P^*\circ F'\eta_{P})P_*\circ\Omega^*P_*\notag\\
({\rm ii})\rightarrow&=Q^*Q_*Q^*\Omega_*\circ Q^*\eta_{Q}F'P_*\circ\Omega^*P_*\notag\\
&=Q^*(Q_*Q^*\Omega_*\circ \eta_{Q}F'P_*)\circ\Omega^*P_*\notag\\
(b)\rightarrow&=Q^*(\eta_{Q}Q_*F\circ \Omega_*)\circ\Omega^*P_*\notag\\
&=Q^*\eta_{Q}Q_*F\circ Q^*\Omega_*\circ\Omega^*P_*\notag\\
&=\delta_{Q}F\circ\Omega,\notag
\end{align}
where the fourth, seventh and ninth equations follow from the above equation $(a)$, the commutativity of (ii) and the above equation $(b)$ respectively.  
Hence $(F,\Omega)$ is a linearized functor.

For any $A\in\mathcal{A}'$, let $\Sigma(A):=\Omega^*(A)$. By constructions, we have $F_{\Omega}\Gamma_{P}(A)=(FP^*(A),\Omega P^*(A)\circ FP^*\eta_P(A))$ and $\Gamma_{Q}F'(A)=(Q^*F'(A),Q^*\eta_{Q}F'(A))$. We show that $\Sigma(-)$ defines a functor morphism $\Sigma:F_{\Omega}\Gamma_{P}\rightarrow\Gamma_{Q}F'$. So we have to verify that $\Omega^*(A)$ is a morphism in  ${\mathcal{B}}_{\mathbb{T}(Q)}$ for each $A\in\mathcal{A}'$, i.e., verify the following diagram is commutative:
 \[\xymatrix{
 FP^*(A)\ar[rrr]^{\Omega P^*(A)\circ FP^*\eta_P(A)} \ar[d]_{\Omega^*(A)}&&&T_{Q}(FP^*(A)) \ar[d]^{T_{Q}(\Omega^*(A))}\\
 Q^*F'(A) \ar[rrr]^{Q^*\eta_{Q}F'(A)}&&&T_{Q}(Q^*F'(A))
 }\]
By the functoriality of $\Omega^*$ and the commutativity of (ii), we see that the above diagram is commutative as follows:
\begin{align}
&\hspace{5mm}T_{Q}(\Omega^*(A))\circ\Omega P^*(A)\circ FP^*\eta_P(A)\notag\\
&=Q^*Q_*\Omega^*(A)\circ Q^*\Omega_*P^*(A)\circ\Omega^*P_*P^*(A)\circ FP^*\eta_P(A)\notag\\
&=Q^*Q_*\Omega^*(A)\circ Q^*\Omega_*P^*(A)\circ\{\Omega^*(P_*P^*(A))\circ FP^*(\eta_P(A))\}\notag\\
{\rm functoriality}\hspace{2mm} {\rm of}\hspace{2mm} \Omega^* \rightarrow&=Q^*Q_*\Omega^*(A)\circ Q^*\Omega_*P^*(A)\circ\{Q^*F'(\eta_P(A))\circ\Omega^*(A)\}\notag\\
&=Q^*\{Q_*\Omega^*(A)\circ\Omega_*P^*(A)\circ F'(\eta_P(A))\}\circ\Omega^*(A)\notag\\
({\rm ii})\rightarrow&=Q^*\eta_{Q}F'(A)\circ\Omega^*(A).\notag
\end{align}
Hence $\Sigma(-)$ defines a functor morphism, and it is an isomorphism.
\end{proof}
\end{lem}

\vspace{5mm}
In the following, we give an important lemma to prove the main theorem.
Notation is  same as the above lemma. Let $G:\mathcal{B}\rightarrow\mathcal{A}$ and $G':\mathcal{B}'\rightarrow\mathcal{A}'$ be functors. Let $\mathcal{C}\subset\mathcal{A}$, $\mathcal{D}\subset\mathcal{B}$, $\mathcal{C}'\subset\mathcal{A}'$ and $\mathcal{D}'\subset\mathcal{B}'$ be full subcategories with $F(\mathcal{A})\subset\mathcal{D}$, $G(\mathcal{D})\subset\mathcal{C}$, $P^*(\mathcal{C}')\subset\mathcal{C}$ and $Q^*(\mathcal{D}')\subset\mathcal{D}$. Now we have the following diagram of functors;
 \[\xymatrix{
 \mathcal{C}\ar@/^/[rrrr]^{F|_{\mathcal{C}}} \ar@{^{(}-{>}}[rd]&&&&\mathcal{D}\ar@/^/[llll]^{G|_{\mathcal{D}}} \ar@{^{(}-{>}}[ld]\\
 &\mathcal{A}\ar@/^/[rr]^{F}\ar@/^/[d]^{P_*}  &&  \mathcal{B} \ar@/^/[d]^{Q_*} \ar@/^/[ll]^{G}& \\
&\mathcal{A}' \ar@/^/[u]^{P^*} \ar@/^/[rr]^{F'}&& \mathcal{B'} \ar@/^/[u]^{Q^*} \ar@/^/[ll]^{G'}&\\
\mathcal{C}'\ar@{^{(}-{>}}[ru] \ar@/^/[rrrr]^{F'|_{\mathcal{C}'}}  \ar[uuu]^{P^*|_{\mathcal{C}'}}&&&&\mathcal{D}' \ar@/^/[llll]^{G'|_{\mathcal{D}'}} \ar@{^{(}-{>}}[lu] \ar[uuu]_{Q^*|_{\mathcal{D}'}}
 }\]
 Let $\Omega_{F}^*:F P^*\xrightarrow{\sim}Q^*F'$,  $\Omega_{F*}:F' P_*\xrightarrow{\sim}Q_* F$, $\Omega_{G}^*:G Q^*\xrightarrow{\sim}P^*G'$ and  $\Omega_{G*}:G' Q_*\xrightarrow{\sim}P_* G$ be isomorphisms of functors such  that the diagrams corresponding to  $({\rm i})$ and $({\rm ii})$ in Lemma \ref{compatible lemma}, namely the following diagrams,  are commutative. 
  \[\xymatrix{
Q^*F'P_*\ar[rr]^{Q^*\Omega_*}&&Q^*Q_*F\ar[d]^{\varepsilon_QF}& &F'\ar[rr]^{\eta_QF'}\ar[d]_{F'\eta_P}&&Q^*Q_*F'\\
FP_*P^*\ar[rr]^{F\varepsilon_P}\ar[u]^{\Omega^*P_*}&&F& &F'P_*P^*\ar[rr]^{\Omega_*P^*}&&Q_*FP^*\ar[u]_{Q_*\Omega^*},
 }\]
  \[\xymatrix{
P^*G'Q_*\ar[rr]^{P^*\Omega_*}&&P^*P_*G\ar[d]^{\varepsilon_PG}& &G'\ar[rr]^{\eta_PG'}\ar[d]_{G'\eta_Q}&&P^*P_*G'\\
GQ_*Q^*\ar[rr]^{G\varepsilon_Q}\ar[u]^{\Omega^*Q_*}&&G& &G'Q_*Q^*\ar[rr]^{\Omega_*Q^*}&&P_*GQ^*\ar[u]_{P_*\Omega^*},
 }\]
 Set $\Omega_F:=Q^*\Omega_{F*}\circ\Omega_F^*P_*$ and $\Omega_G:=P^*\Omega_{G*}\circ\Omega_G^*Q_*$. 

\vspace{2mm} 
\begin{lem}\label{main lemma}
Notation is same as above. Assume that the adjunction morphisms $\eta_P:{\rm id}\rightarrow P_*P^*$ and $\eta_Q:{\rm id}\rightarrow Q_*Q^*$  are split mono, and for any $D\in\mathcal{D}$ and $A\in\mathcal{A}$ there is a natural isomorphism
 $$\Sigma(D,A):{\rm Hom}_{\mathcal{B}}(D,F(A))\cong{\rm Hom}_{\mathcal{A}}(G(D),A)$$ which is functorial in  $D$ and $A$. Then, if $F|_{\mathcal{C}}:\mathcal{C}\rightarrow\mathcal{D}$ is fully faithful, so is $F'|_{\mathcal{C}'}:\mathcal{C}'\rightarrow\mathcal{D}'$. Moreover, if  $F|_{\mathcal{C}}$ and $G'|_{\mathcal{D}'}$ are fully faithful and the following diagram $(\star)$ of functor morphisms is commutative, $F'|_{\mathcal{C}'}$ is an equivalence. Define a diagram by
\[\xymatrix{
 (\star): &GFP^*|_{\mathcal{C}'} \ar[rr]^{GFP^*\eta_P}\ar[d]_{\omega P^*}&&GFP^*P_*P^*|_{\mathcal{C}'} \ar[r]^{G\Omega_F P^*}&GQ^*Q_*FP^*|_{\mathcal{C}'}  \ar[r]^{\Omega_GFP^*} &P^*P_*GFP^*|_{\mathcal{C}'} \ar[d]^{P^*P_*\omega P^*}\\
 & P^*|_{\mathcal{C}'} \ar[rrrr]^{P^*\eta_P} &&&& P^*P_*P^*|_{\mathcal{C}'}, 
 }\]
 where $\omega:GF|_{\mathcal{C}}\rightarrow{\rm id}_{\mathcal{C}}$ is the adjunction morphism of the adjoint pair $(G|_{\mathcal{D}}\dashv F|_{\mathcal{C}})$ given by $\Sigma(-,*)$.
\begin{proof}
By the assumption and Lemma \ref{compatible lemma}, $(F,\Omega_F)$ and $(G,\Omega_G)$ are linearized functor, and  
 we have the following commutative diagram of functors 
\[\xymatrix{
\mathcal{A}_{\mathbb{T}(P)}^{\mathcal{C}}\ar[rr]^{F_{\Omega_{F}}|_{\mathcal{A}_{\mathbb{T}(P)}^{\mathcal{C}}}}&&\mathcal{B}_{\mathbb{T}(Q)}^{\mathcal{D}}\\
\mathcal{C}'\ar[rr]^{F'|_{\mathcal{C}'}} \ar[u]^{\Gamma_{P}|_{\mathcal{C}'}}&&\mathcal{D}' \ar[u]_{\Gamma_{Q}|_{\mathcal{D}'}}.
 }\]
Since the adjunction morphisms $\eta_P$ and $\eta_Q$ are split mono, the comparison functors $\Gamma_P:\mathcal{A}'\rightarrow\mathcal{A}_{\mathbb{T}(P)}$ and $\Gamma_Q:\mathcal{B}'\rightarrow\mathcal{B}_{\mathbb{T}(Q)}$ are fully faithful functors by Proposition \ref{comparison theorem}. 

 We show that if $F|_{\mathcal{C}}$ is fully faithful, then the condition $(\ast)$ in Proposition \ref{main prop} is satisfied, i.e. the map $F:{\rm Hom}(C_1,T^n_P(C_2))\rightarrow{\rm Hom}(F(C_1),F(T^n_P(C_2)))$ is bijective for any $C_i\in\mathcal{C}$ and $n=1,2$. Consider the following commutative diagram of maps
 \[\xymatrix{
{\rm Hom}(C_1,T^n_P(C_2)) \ar[rr]^{(-)\circ\omega(C_1)} \ar[rd]_{F}&&{\rm Hom}(G(F(C_1)),T^n_P(C_2))\\
&{\rm Hom}(F(C_1),F(T^n_P(C_2))) \ar[ru]_{\Sigma(F(C_1),T^n_P(C_2))}&
 }\]
Since $F|_{\mathcal{C}}$ is fully faithful, $\omega(C_1)$ is an isomorphism, whence  the maps in the above diagram except for $F$ are bijective. Hence, the condition $(\ast)$ in Proposition \ref{main prop} is satisfied, and  we see that if $F|_{\mathcal{C}}$ is fully faithful, then  $F'|_{\mathcal{C}'}$ is also fully faithful by Proposition \ref{main prop}.

Assume that $F|_{\mathcal{C}}$ and $G'|_{\mathcal{D}'}$ are fully faithful and that the diagram $(\star)$ is commutative. Since the diagram $(\star)$ is commutative, the functor morphism $\omega P^*|_{\mathcal{C}'}:GF|_{\mathcal{C}}P^*|_{\mathcal{C}'}\rightarrow P^*|_{\mathcal{C}'}$ induces a functor morphism  $\omega':G_{\Omega_G}\circ F_{\Omega_F}\circ\Gamma_P|_{\mathcal{C}'}\xrightarrow{\sim}\Gamma_P|_{\mathcal{C}'}$.  Since $F|_{\mathcal{C}}$ is fully faithful, $\omega'$ is an isomorphism of functors. Since we have $G_{\Omega_G}\circ F_{\Omega_F}\circ\Gamma_P|_{\mathcal{C}'}\cong \Gamma_P\circ G'\circ F'|_{\mathcal{C}'}$ and $\Gamma_P$ is fully faithful, the functor isomorphism $\omega'$ implies an isomorphism of functors $G'F'|_{\mathcal{C}'}\xrightarrow{\sim}{\rm id}_{\mathcal{C}'}$. Hence  $G'|_{\mathcal{D}'}:\mathcal{D}'\rightarrow\mathcal{C}'$ is an equivalence, and therefore, $F'|_{\mathcal{C}'}$ is also an equivalence.
\end{proof}
\end{lem}

\vspace{7mm}
\section{Derived factorization categories}
\vspace{2mm}
In this section, we give definitions and foundations of categories with potentials, and construct  derived factorization categories of them. We also construct functors between factorization categories from cwp-functors.

\vspace*{2mm}
\subsection{Factorization categories}Let $\mathcal{A}$ be an exact category in the sense of  Quillen (see \cite{quillen}). First of all, we define potentials on $\mathcal{A}$.

\vspace*{2mm}
\begin{dfn}
A \textbf{potential} of $\mathcal{A}$ is a pair $(\Phi,W)$ of an exact equivalence $\Phi:\mathcal{A}\xrightarrow{\sim}\mathcal{A}$ and a functor morphism $W:{{\rm id}_{\mathcal{A}}}\rightarrow\Phi$ such that $\Phi W=W\Phi$.
The triple $(\mathcal{A},\Phi,W)$ is called a $\textbf{category with a potential}$.
\end{dfn}

\vspace*{2mm}
Let  $(\Phi,W)$ be a potential of $\mathcal{A}$. A \textbf{factorization} of $(\Phi,W)$ is a sequence in $\mathcal{A}$
$$A=\Bigl(A_1\xrightarrow{\varphi^A_1}A_0\xrightarrow{\varphi^A_0}\Phi(A_1)\Bigr)$$
such that $\varphi^A_0\circ\varphi^A_1=W(A_1)$ and $\Phi(\varphi^A_1)\circ\varphi^A_0=W(A_0)$. Objects $A_1$ and $A_0$ in the above sequence are called  \textbf{components} of the factorization $A$.

\vspace*{2mm}
\begin{dfn}
For a category with a potential $(\mathcal{A},\Phi,W)$, we define a dg-category $\mathfrak{F}(\mathcal{A},\Phi,W)$, whose objects are factorizations of $(\Phi,W)$, as follows. 
For two factorizations $A,B\in\mathfrak{F}(\mathcal{A},\Phi,W)$,  the set of morphisms  ${\rm Hom}(A,B)$ is a complex
$${\rm Hom}(A,B):=\bigoplus_{n\in\mathbb{Z}}{\rm Hom}(A,B)^n$$
with a differential $d$ on ${\rm Hom}(A,B)$ given by
$$d(f):=\varphi^B\circ f-(-1)^{{\rm deg}(f)}f\circ \varphi^A \hspace{3mm}{\rm if}\hspace{1mm} f\in{\rm Hom}(A,B)^{{\rm deg}(f)},$$
where $${\rm Hom}(A,B)^{2n}:={\rm Hom}(A_1,\Phi^n(B_1))\oplus{\rm Hom}(A_0,\Phi^n(B_0))$$
$${\rm Hom}(A,B)^{2n+1}:={\rm Hom}(A_1,\Phi^{n}(B_0))\oplus{\rm Hom}(A_0,\Phi^{n+1}(B_1)).$$
We call $\mathfrak{F}(\mathcal{A},\Phi,W)$ the \textbf{factorization category} of $(\mathcal{A},\Phi,W)$.
\end{dfn}

\vspace*{3mm}For any dg-category $\mathcal{D}$, we define two categories $Z^0(\mathcal{D})$ and $H^0(\mathcal{D})$ whose objects are same as $\mathcal{D}$ and whose morphisms are defined as follows;  $${\rm Hom}_{Z^0(\mathcal{D})}(A,B):=Z^0({\rm Hom}_{\mathcal{D}}(A,B))$$
$${\rm Hom}_{H^0(\mathcal{D})}(A,B):=H^0({\rm Hom}_{\mathcal{D}}(A,B)),$$
where ${\rm Hom}_{\mathcal{D}}(A,B)$ in the right hand sides are considered as complexes.

\vspace{2mm}
\begin{rem} The categories $Z^0(\mathfrak{F}(\mathcal{A},\Phi,W))$ and $H^0(\mathfrak{F}(\mathcal{A},\Phi,W))$ are generalizations of categories of classical matrix factorizations introduced by Eisenbud \cite{eisenbud}.

  Let $A, B$ be objects in $Z^0(\mathfrak{F}(\mathcal{A},\Phi,W))$. Then the set of morphisms from $A$ to $B$ can be  described as follows:
\begin{center}
{\rm Hom}$_{Z^0(\mathfrak{F}(\mathcal{A},\Phi,W))}(A,B)\cong\{(f_1,f_0)\mid f_i:A_i\rightarrow B_i$ and the diagram $(\star)$ is commutative.$\}$
\end{center}
\[\xymatrix{
(\star):& A_1\ar[rr]^{\varphi_1^A}\ar[d]_{f_1}&&A_0\ar[rr]^{\varphi_0^A}\ar[d]^{f_0}&&\Phi(A_1)\ar[d]^{\Phi(f_1)}\\
&B_1\ar[rr]^{\varphi_1^B}&&B_0\ar[rr]^{\varphi_0^B}&&\Phi(B_1)
 }\]

The set of morphisms in the category $H^0(\mathfrak{F}(\mathcal{A},\Phi,W))$ can be described as the set of homotopy equivalence classes of {\rm Hom}$_{Z^0(\mathfrak{F}(\mathcal{A},\Phi,W))}(A,B)$$;$
$${\rm Hom}_{H^0(\mathfrak{F}(\mathcal{A},\Phi,W))}(A,B)\cong{\rm Hom}_{Z^0(\mathfrak{F}(\mathcal{A},\Phi,W))}(A,B)/\sim.$$
Two morphisms $f=(f_1,f_0)$ and $g=(g_1,g_0)$ in ${\rm Hom}_{Z^0(\mathfrak{F}(\mathcal{A},\Phi,W))}(A,B)$ are \textbf{homotopy equivalence} if there exist morphisms 
$$h_0:A_0\rightarrow B_1\hspace{3mm}{\rm and}\hspace{3mm}h_1:\Phi(A_1)\rightarrow B_0$$
such that $f_0=\varphi_1^Bh_0+h_1\varphi_0^A$ and $\Phi(f_1)=\varphi_0^Bh_1+\Phi(h_0)\Phi(\varphi_1^A)$.
\end{rem}

\begin{dfn}
For each $i=0,1$, we have a natural exact functor
$$(-)_i:Z^0(\mathfrak{F}(\mathcal{A},\Phi,W))\rightarrow\mathcal{A}$$
defined by $(A_1\xrightarrow{\varphi^A_1} A_0\xrightarrow{\varphi^A_0} \Phi(A_1))_i:=A_i$. This functor extends to an exact functor of their derived categories,
$$(-)_i:{\rm D^b}(Z^0(\mathfrak{F}(\mathcal{A},\Phi,W)))\rightarrow{\rm D^b}(\mathcal{A}).$$
\end{dfn}

\vspace{2mm}

\vspace{2mm}

\begin{prop}\label{Z^0exact}
The category $Z^0(\mathfrak{F}(\mathcal{A},\Phi,W))$ is an exact category. 
Furthermore, if $\mathcal{A}$ is abelian category, then $Z^0(\mathfrak{F}(\mathcal{A},\Phi,W))$ is an abelian category.  

\begin{proof}
Assume that $\mathcal{A}$ is abelian category. At first, we show that  $Z^0(\mathfrak{F}(\mathcal{A},\Phi,W))$ is an abelian category.  
 For any morphism $f=(f_1,f_0):A\rightarrow B$ in $Z^0(\mathfrak{F}(\mathcal{A},\Phi,W))$, let $$k_i:K_i\hookrightarrow A_i$$ be the kernel of $f_i:A_i\rightarrow B_i$ for each $i=0,1$. By the universal property of kernels, there exist morphisms $\varphi_1^K:K_1\rightarrow K_0$ and $\varphi_0^K:K_0\rightarrow \Phi(K_1)$ such that the following diagram is commutative:
\[\xymatrix{
 K_1\ar[rr]^{\varphi_1^K}\ar[d]_{k_1}&&K_0\ar[rr]^{\varphi_0^K}\ar[d]^{k_0}&&\Phi(K_1)\ar[d]^{\Phi(k_1)}\\
A_1\ar[rr]^{\varphi_1^A}&&A_0\ar[rr]^{\varphi_0^A}&&\Phi(A_1)
 }\]
Since we have an equality $\Phi(k_1)\circ(\varphi_0^K\circ\varphi_1^K)=\Phi(k_1)\circ W(K_1)$, and $\Phi(k_1)$ is injective, we have $\varphi_0^K\circ\varphi_1^K=W(K_1)$.  Similarly, we see that $\Phi(\varphi_1^K)\circ\varphi_0^K=W(K_0)$. Hence, $$K:=(K_1\xrightarrow{\varphi_1^K}K_0\xrightarrow{\varphi_0^K}\Phi(K_1))$$ is an object of $Z^0(\mathfrak{F}(\mathcal{A},\Phi,W))$. Since, $K_i$ is the kernel of $f_i$, $K$ is the kernel of $f$. Similarly, we see that $Z^0(\mathfrak{F}(\mathcal{A},\Phi,W))$ admits cokernel of any morphism, and we obtain a natural isomorphism ${\rm Im}(f)\cong {\rm Coim}(f)$. Hence, $Z^0(\mathfrak{F}(\mathcal{A},\Phi,W))$ is an abelian category.

Next, we show that $Z^0(\mathfrak{F}(\mathcal{A},\Phi,W))$ is an exact category. Let $\overline{\mathcal{A}}$ be the category of left exact functors from $\mathcal{A}^{\rm op}$ to the category of abelian groups (in a fixed universe containing $\mathcal{A}$). By \cite{quillen}, the category $\overline{\mathcal{A}}$ is an abelian category, and  we have a fully faithful functor $$h:\mathcal{A}\rightarrow\overline{\mathcal{A}},$$
such that $h$ embeds $\mathcal{A}$ as a full subcategory of $\overline{\mathcal{A}}$ closed under extensions, and a sequence 
$$A'\rightarrow A\rightarrow A''$$
in $\mathcal{A}$ is exact if and only if $h$ carries it into an exact sequence in $\overline{\mathcal{A}}$ (the category $\overline{\mathcal{A}}$ is called {\it the abelian envelope} of $\mathcal{A}$).
We define an exact autoequivalence $\overline{\Phi}:\overline{\mathcal{A}}\rightarrow\overline{\mathcal{A}}$ and a functor morphism $\overline{W}:{\rm id}_{\overline{\mathcal{A}}}\rightarrow\overline{\Phi}$ as follows:
For an object $F\in\overline{\mathcal{A}}$, we define $\overline{\Phi}(F):=F\circ(\Phi^{\rm op})^{-1}\in\overline{\mathcal{A}}$ and $\overline{W}(F):=FW^{\rm op}(\Phi^{\rm op})^{-1}:F\rightarrow F\circ(\Phi^{\rm op})^{-1}$ where $W^{\rm op}(\Phi^{\rm op})^{-1}:{\rm id}_{\mathcal{A}^{\rm op}}\rightarrow(\Phi^{\rm op})^{-1}$ is the composition $${\rm id}_{\mathcal{A}^{\rm op}}\xrightarrow{\sim}\Phi^{\rm op}\circ(\Phi^{\rm op})^{-1}\xrightarrow{W^{\rm op}(\Phi^{\rm op})^{-1}}(\Phi^{\rm op})^{-1}$$
Since the functor $h$ is compatible with potentials,  it induces a fully faithful functor $$Z^0(\mathfrak{F}(\mathcal{A},\Phi,W))\rightarrow Z^0(\mathfrak{F}(\overline{\mathcal{A}},\overline{\Phi},\overline{W})).$$
By this embedding, we obtain a natural structure of exact category on $Z^0(\mathfrak{F}(\mathcal{A},\Phi,W))$.
\end{proof}

\end{prop}

\vspace{4mm}
For an object $A\in Z^0(\mathfrak{F}(\mathcal{A},\Phi,W))$, we can construct a twisted-periodic infinite sequence Com$(A)=({\rm Com}(A)^{\text{\tiny{\textbullet}}}, d_A^{\text{\tiny{\textbullet}}})$  in $\mathcal{A}$ with $d_A^{i+1}\circ d_A^i=W({\rm Com}(A)^i)$ as follows;

$${\rm Com}(A)^{2i}:=\Phi^i(A_0) , \hspace{3mm}\hspace{3mm}{\rm Com}(A)^{2i-1}:=\Phi^i(A_1),$$
$$d_A^{2i}:=\Phi^i(\phi_0^A),\hspace{3mm}\hspace{3mm}d_A^{2i-1}:=\Phi^i(\phi_1^A).$$
For a morphism $f=(f_1,f_0)\in{{\rm Hom}_{Z^0(\mathfrak{F}(\mathcal{A},\Phi,W))}(A,B)}\subset{{\rm Hom}(A_1,B_1)}\oplus{{\rm Hom}(A_0,B_0)}$, we define  
 a morphism Com$(f)=({\rm Com}(f)^{\text{\tiny{\textbullet}}})$ from Com$(A)$ to Com$(B)$ as follows:

$${\rm Com}(f)^{2i}:=\Phi^i(f_0)\hspace{6mm}{\rm Com}(f)^{2i-1}:=\Phi^i(f_1)$$

\vspace*{2mm}
\begin{dfn}
Let $C^{\text{\tiny{\textbullet}}}=(\cdot\cdot\cdot\rightarrow C^i\xrightarrow{\delta_{C^{\text{\tiny{\textbullet}}}}^i}C^{i+1}\rightarrow\cdot\cdot\cdot)$ be a bounded complex of $Z^0(\mathfrak{F}(\mathcal{A},\Phi,W))$. We define the {\bf totalization} of $C^{\text{\tiny{\textbullet}}}$ as an object Tot$(C^{\text{\tiny{\textbullet}}})\in Z^0(\mathfrak{F}(\mathcal{A},\Phi,W))$ in a similar way to construct the total complex of a double complex, i.e.,
$${\rm Tot}(C^{\text{\tiny{\textbullet}}}):=(T_1\xrightarrow{t_1}T_0\xrightarrow{t_0}\Phi(T_1)),$$
where
$$T_l:=\bigoplus_{i+j=-l}{\rm Com}(C^i)^j,$$
$$t_l|_{{\rm Com}(C^i)^j}:={\rm Com}(\delta_{C^{\text{\tiny{\textbullet}}}}^i)^j+(-1)^id_{C^i}^j.$$

Let $\varphi^{\text{\tiny{\textbullet}}}:C^{\text{\tiny{\textbullet}}}\rightarrow D^{\text{\tiny{\textbullet}}}$ be a morphism of complexes of $Z^0(\mathfrak{F}(\mathcal{A},\Phi,W))$. We define a morphism Tot$(\varphi^{\text{\tiny{\textbullet}}}): {\rm Tot}(C^{\text{\tiny{\textbullet}}})\rightarrow{\rm Tot}(D^{\text{\tiny{\textbullet}}})$ in $Z^0(\mathfrak{F}(\mathcal{A},\Phi,W))$ as 
$${\rm Tot}(\varphi^{\text{\tiny{\textbullet}}}):=(\tau_1,\tau_0),$$
where $$\tau_l|_{{\rm Com}(C^i)^j}:={\rm Com}(\varphi^i)^j.$$
Taking totalizations gives an exact functor $${\rm Tot}:{\rm Ch^b}(Z^0(\mathfrak{F}(\mathcal{A},\Phi,W)))\rightarrow Z^0(\mathfrak{F}(\mathcal{A},\Phi,W)).$$
\end{dfn}

\vspace{4mm}
In what follows, we will see that the  category $H^0(\mathfrak{F}(\mathcal{A},\Phi,W))$ has a structure of a triangulated category.

\vspace*{2mm}
\begin{dfn}
We define an automorphism $T$ on $H^0(\mathfrak{F}(\mathcal{A},\Phi,W))$, which is called \textbf{the shift functor}, as follows.
For an object $A\in H^0(\mathfrak{F}(\mathcal{A},\Phi,W))$, we define an object $T(A)$ as
$$T(A):=(A_0\xrightarrow{-\varphi^A_0}\Phi(A_1)\xrightarrow{-\Phi(\varphi^A_1)}\Phi(A_0))$$
and for a morphism $f\in {\rm Hom}(A,B)$, a morphism $T(f)\in{\rm Hom}(T(A),T(B))$ is suitably defined. For any integer $n\in\mathbb{Z}$, denote by $(-)[n]$ the functor $T^n(-)$.
\end{dfn}

\vspace*{2mm}
\begin{dfn}
Let $f : A\rightarrow B$ be a morphism in $Z^0(\mathfrak{F}(\mathcal{A},\Phi,W))$. We define its {\bf mapping cone} Cone$(f)$ to be the totalization of the complex $$(\cdot\cdot\cdot\rightarrow0\rightarrow A\xrightarrow{f} B\rightarrow0\rightarrow\cdot\cdot\cdot)$$ with $B$ in degree zero.

A sequence in $H^0(\mathfrak{F}(\mathcal{A},\Phi,W))$ of the form

$$A\xrightarrow{f}B\xrightarrow{i}{\rm Cone}(f)\xrightarrow{p}A[1],$$
where $i$ is the natural injection and $p$ is the natural projection, is called a $\textbf{standard triangle}$ and  a sequence which is isomorphic to a standard triangle is called $\textbf{distinguished triangle}$.
\end{dfn}

\vspace*{2mm}
\begin{prop}
$H^0(\mathfrak{F}(\mathcal{A},\Phi,W))$ is a triangulated category with respect to its shift functor  and its distinguished triangles. 
\begin{proof}This follows from an argument similar to a proof showing that homotopy categories of exact categories are triangulated categories.
\end{proof}
\end{prop}

\vspace*{2mm}
Following Positselski (cf. \cite{posi} or \cite{efi-posi}), we define derived  factorization categories.


\begin{dfn}
Denote by  ${\rm Acycl}^{\rm abs}(\mathcal{A},\Phi,W))$ the smallest thick subcategory of $H^0(\mathfrak{F}(\mathcal{A},\Phi,W))$ containing all totalizations of short exact sequences in $Z^0(\mathfrak{F}(\mathcal{A},\Phi,W)$.   $E\in H^0(\mathfrak{F}(\mathcal{A},\Phi,W))$ is called $\textbf{absolutely acyclic}$  if it lies in ${\rm Acycl}^{\rm abs}(\mathcal{A},\Phi,W))$.
The $\textbf{absolute derived factorization}$ \textbf{category} of $(\mathcal{A},\Phi,W)$ is the Verdier quotient 
$${\rm D^{abs}}(\mathcal{A},\Phi,W):=H^0(\mathfrak{F}(\mathcal{A},\Phi,W))/{\rm Acycl}^{\rm abs}(\mathcal{A},\Phi,W)$$
\end{dfn}

\vspace*{2mm}
\begin{dfn}
Assume $\mathcal{A}$ admits small coproducts. Denote ${\rm Acycl}^{\rm co}(\mathcal{A},\Phi,W))$ the smallest thick subcategory of $H^0(\mathfrak{F}(\mathcal{A},\Phi,W))$ containing all totalizations of short exact sequences in $Z^0(\mathfrak{F}(\mathcal{A},\Phi,W)$ and closed under taking small coproducts.
$E\in H^0(\mathfrak{F}(\mathcal{A},\Phi,W))$ is called $\textbf{co-acyclic}$  if it lies in ${\rm Acycl}^{\rm co}(\mathcal{A},\Phi,W)$.
The $\textbf{co-derived factorization category}$ of $(\mathcal{A},\Phi,W)$ is the Verdier quotient 
$${\rm D^{co}}(\mathcal{A},\Phi,W):=H^0(\mathfrak{F}(\mathcal{A},\Phi,W))/{\rm Acycl}^{\rm co}(\mathcal{A},\Phi,W)$$
\end{dfn}

\begin{rem}
$(1)$Let $\mathcal{E}$ be an exact category, and take a  complex $E^{\text{\tiny{\textbullet}}}$ in $\mathcal{E}$; 
$$E^{\text{\tiny{\textbullet}}}=\cdot\cdot\cdot \rightarrow E^{n-1}\xrightarrow{d^{n-1}} E^n\xrightarrow{d^n}E^{n+1}\rightarrow\cdot\cdot\cdot.$$
We say that the complex $E^{\text{\tiny{\textbullet}}}$ is \textbf{exact} if all kernels and images of differentials exist, and for any $n\in\mathbb{Z}$, we have natural isomorphisms $${\rm Im}(d^{n-1})\cong {\rm Ker}(d^{n}).$$
Let $\mathcal{B}$ be an abelian category, and let $\mathcal{C}$ be a strictly full additive subcategory of $\mathcal{B}$ which is closed under extensions. The category $\mathcal{C}$ has a natural structure of an exact category. 
If $\mathcal{C}$ admits either all kernels or all cokernels, then a bounded complex in $\mathcal{C}$ is exact in the above sense if and only if the complex is exact  in $\mathcal{B}$. \\
$(2)$Note that in the definitions of ${\rm Acycl}^{\rm abs}(\mathcal{A},\Phi,W)$ and  ${\rm Acycl}^{\rm co}(\mathcal{A},\Phi,W)$, we can replace the words $^{\rotatebox[origin=C]{180}{"}}$totalizations of short exact sequences" with $^{\rotatebox[origin=C]{180}{"}}$totalizations of bounded exact sequences". 

\end{rem}

\vspace{2mm}






By the next lemma, we see that the totalization functor $${\rm Tot}:{\rm Ch^b}(Z^0(\mathfrak{F}(\mathcal{A},\Phi,W)))\rightarrow Z^0(\mathfrak{F}(\mathcal{A},\Phi,W))$$ induces a functor $${\rm Tot}:{\rm K^b}(Z^0(\mathfrak{F}(\mathcal{A},\Phi,W)))\rightarrow H^0(\mathfrak{F}(\mathcal{A},\Phi,W))$$
which is an exact functor of triangulated categories. This functor naturally induces an exact functor 
$${\rm Tot}: {\rm D^b}(Z^0(\mathfrak{F}(\mathcal{A},\Phi,W)))\rightarrow{\rm D^{abs}}(\mathcal{A},\Phi,W).$$

\vspace*{2mm}
\begin{lem}
Let $\varphi^{\text{\tiny{\textbullet}}}:C^{\text{\tiny{\textbullet}}}\rightarrow D^{\text{\tiny{\textbullet}}}$ be a morphism in  ${\rm Ch^b}(Z^0(\mathfrak{F}(\mathcal{A},\Phi,W)))$. If $\varphi^{\text{\tiny{\textbullet}}}$ is homotopic to zero, i.e. $\varphi^{\text{\tiny{\textbullet}}}=0$ in ${\rm K^b}(Z^0(\mathfrak{F}(\mathcal{A},\Phi,W)))$, then {\rm Tot}$(\varphi^{\text{\tiny{\textbullet}}})=0$ in $H^0(\mathfrak{F}(\mathcal{A},\Phi,W))$.
\begin{proof}

Let $\delta_{C^{\text{\tiny{\textbullet}}}}^i:C^i\rightarrow C^{i+1}$ and $\delta_{D^{\text{\tiny{\textbullet}}}}^i:D^i\rightarrow D^{i+1}$ be differentials of complexes ${C^{\text{\tiny{\textbullet}}}}$ and ${D^{\text{\tiny{\textbullet}}}}$, and set 
$$S=(S_1\xrightarrow{s_1}S_0\xrightarrow{s_0}\Phi(S_1)):={\rm Tot}(C^{\text{\tiny{\textbullet}}}),$$ 
$$T=(T_1\xrightarrow{t_1}T_0\xrightarrow{t_0}\Phi(T_1)):={\rm Tot}(D^{\text{\tiny{\textbullet}}})$$ and 
$$\tau=(\tau_1,\tau_0):={\rm Tot}(\varphi^{\text{\tiny{\textbullet}}}),$$ where $\tau_l:S_l\rightarrow T_l$.
If  $\varphi^{\text{\tiny{\textbullet}}}=0$ in ${\rm K^b}(Z^0(\mathfrak{F}(\mathcal{A},\Phi,W)))$, then there exist  morphisms $h^i : C^i\rightarrow D^{i-1}$ such that $\varphi^i=\delta_{D^{\text{\tiny{\textbullet}}}}^{i-1} h^i+h^{i+1}\delta_{C^{\text{\tiny{\textbullet}}}}^i$.   We define two morphisms $\sigma_0: S_0\rightarrow T_1$ and $\sigma_1:\Phi(S_1)\rightarrow T_0$ in  $\mathcal{A}$ as 
$$\sigma_l|_{{\rm Com}(C^i)^j}:={\rm Com}(h^i)^j$$
for each $l=0,1$. Then we have 
\begin{align}
&\hspace{5mm}\bigl(s_1\sigma_0+\Phi(\sigma_1)t_0\bigr)|_{{\rm Com}(C^i)^j}\notag\\
&=\bigl({\rm Com}(\delta_{D^{\text{\tiny{\textbullet}}}}^{i-1})^j+(-1)^{i-1}d_{D^{i-1}}^j\bigr){\rm Com}(h^i)^j+\Phi(\sigma_1)\bigl({\rm Com}(\delta_{C^{\text{\tiny{\textbullet}}}}^i)^j+(-1)^id_{C^i}^j\bigr) \notag\\
&={\rm Com}(\delta_{D^{\text{\tiny{\textbullet}}}}^{i-1})^j{\rm Com}(h^i)^j+(-1)^{i-1}d_{D^{i-1}}^j{\rm Com}(h^i)^j+{\rm Com}(h^{i+1})^j{\rm Com}(\delta_{C^{\text{\tiny{\textbullet}}}}^i)^j+(-1)^i{\rm Com}(h^i)^{j+1}d_{C^i}^j\notag\\
&={\rm Com}(\delta_{D^{\text{\tiny{\textbullet}}}}^{i-1}h^i+h^{i+1}\delta_{C^{\text{\tiny{\textbullet}}}}^i)^j\notag\\
&=\tau_0|_{{\rm Com}(C^i)^j},\notag
\end{align}
where $d_{D^{i-1}}^j$ and $d_{C^i}^j$ are morphisms in the infinite sequences ${\rm Com}(D^{i-1})$ and ${\rm Com}(C^i)$ respectively. Hence, we have $\tau_0=s_1\sigma_0+\Phi(\sigma_1)t_0$. Similarly, we obtain $\Phi(\tau_1)=s_0\sigma_1+\Phi(\sigma_0)\Phi(t_1)$. Hence, ${\rm Tot}(\varphi^{\text{\tiny{\textbullet}}})=0$ in $H^0(\mathfrak{F}(\mathcal{A},\Phi,W))$.
\end{proof}
\end{lem}

\vspace{3mm}

Consider an exact functor of exact categories
$$\tau:\mathcal{A}\rightarrow Z^0(\mathfrak{F}(\mathcal{A},\Phi,0)),$$
which is defined by 
$$\tau(A):=(0\longrightarrow A  \longrightarrow0).$$
Then this functor induces an exact functor of triangulated categories
$$\tau:{\rm D^b}(\mathcal{A})\rightarrow {\rm D^b}(Z^0(\mathfrak{F}(\mathcal{A},\Phi,0))).$$

\begin{dfn}\label{upsilon}
We define an exact functor $$\Upsilon:{\rm D^b}(\mathcal{A})\rightarrow {\rm D^{abs}}(\mathcal{A},\Phi,0)$$
as the composition
$${\rm D^b}(\mathcal{A})\xrightarrow{\tau} {\rm D^b}(Z^0(\mathfrak{F}(\mathcal{A},\Phi,0)))\xrightarrow{\rm Tot} {\rm D^{abs}}(\mathcal{A},\Phi,0).$$
\end{dfn}

\vspace*{3mm}
\subsection{cwp-functors}

Let $(\mathcal{A},\Phi_{\mathcal{A}},W_{\mathcal{A}})$, $(\mathcal{B},\Phi_{\mathcal{B}},W_{\mathcal{B}})$ and $(\mathcal{C},\Phi_{\mathcal{C}},W_{\mathcal{C}})$  be categories with potentials.

\vspace*{2mm}
\begin{dfn}\label{cwp-functor}
Let $F:\mathcal{A}\rightarrow\mathcal{B}$ be an additive functor. $F$ is \textbf{compatible with potentials}  with respect to $(\Phi_{\mathcal{A}},W_{\mathcal{A}})$ and $(\Phi_{\mathcal{B}},W_{\mathcal{B}})$ if there exists a functor isomorphism $\sigma:F\Phi_{\mathcal{A}}\xrightarrow{\sim}\Phi_{\mathcal{B}}F$ such that $W_{\mathcal{B}}F=\sigma\circ FW_{\mathcal{A}}$. We call the pair $(F,\sigma)$  a \textbf{cwp-functor} and write $$(F,\sigma):(\mathcal{A},\Phi_{\mathcal{A}},W_{\mathcal{A}})\rightarrow(\mathcal{B},\Phi_{\mathcal{B}},W_{\mathcal{B}}).$$
We just say $F$ is a cwp-functor and write $F:(\mathcal{A},\Phi_{\mathcal{A}},W_{\mathcal{A}})\rightarrow(\mathcal{B},\Phi_{\mathcal{B}},W_{\mathcal{B}})$ when there is no confusion about what $\sigma$ is.
\end{dfn}

\vspace*{2mm}
 A cwp-functor $(F,\sigma):(\mathcal{A},\Phi_{\mathcal{A}},W_{\mathcal{A}})\rightarrow(\mathcal{B},\Phi_{\mathcal{B}},W_{\mathcal{B}})$ induces a natural dg-functor $$\mathfrak{F}(F,{\sigma}):\mathfrak{F}(\mathcal{A},\Phi_{\mathcal{A}},W_{\mathcal{A}})\rightarrow\mathfrak{F}(\mathcal{B},\Phi_{\mathcal{B}},W_{\mathcal{B}})$$ as follows. For objects $A,B\in\mathfrak{F}(\mathcal{A},\Phi_{\mathcal{A}},W_{\mathcal{A}})$ and for a morphism $f\in{\rm Hom}(A_i,\Phi_{\mathcal{A}}^n(B_j))$, we define
 $$\mathfrak{F}(F,{\sigma})(A):=\Bigl(F(A_1)\xrightarrow{F(\varphi^A_1)}F(A_0)\xrightarrow{\sigma(A_1)\circ F(\varphi^A_0)}\Phi_{\mathcal{B}}(F(A_1))\Bigr)$$
 and
 $$\mathfrak{F}(F,{\sigma})(f):=\sigma^n(B_j)\circ F(f)\in{\rm Hom}(F(A_i),\Phi_{\mathcal{B}}^n(F(B_j))),$$
 where $\sigma^n:F\Phi_{\mathcal{A}}^n\xrightarrow{\sim}\Phi_{\mathcal{B}}^nF$ is the functor isomorphism induced by $\sigma$.
 By the construction, we see that  the morphism $\mathfrak{F}(F,{\sigma}):{\rm Hom}(A,B)\rightarrow{\rm Hom}(\mathfrak{F}(F,{\sigma})(A),\mathfrak{F}(F,{\sigma})(B))$ preserves degrees of complexes and is compatible with differentials.
 
 \vspace{3mm}
 In the following lemma, we give  fundamental properties of  dg-functors give as $\mathfrak{F}(-)$. Since the proof is straightforward, we skip the proof.
 
  \vspace*{2mm}
 \begin{lem}\label{cwp-functor lemma} Let $(F,\sigma):(\mathcal{A},\Phi_{\mathcal{A}},W_{\mathcal{A}})\rightarrow(\mathcal{B},\Phi_{\mathcal{B}},W_{\mathcal{B}})$ and $(G,\tau):(\mathcal{B},\Phi_{\mathcal{B}},W_{\mathcal{B}})\rightarrow(\mathcal{C},\Phi_{\mathcal{C}},W_{\mathcal{C}})$ be cwp-functors. Then we have

$(1)$ $(G\circ F,\tau F\circ G\sigma)$ is a cwp-functor, and we have $$\mathfrak{F}(G\circ F,{\tau F\circ G\sigma})=\mathfrak{F}(G,{\tau})\circ \mathfrak{F}(F,{\sigma}).$$

$(2)$ If $F$ is fully faithful, so is $\mathfrak{F}(F,{\sigma})$. 

$(3)$ If $F$ is an equivalence, so is $\mathfrak{F}(F,{\sigma})$.
 \end{lem}

 \vspace*{2mm}
\begin{dfn}Let $(F,\sigma), (F',\sigma'):(\mathcal{A},\Phi_{\mathcal{A}},W_{\mathcal{A}})\rightarrow(\mathcal{B},\Phi_{\mathcal{B}},W_{\mathcal{B}})$  be cwp-functors and let $\alpha:F\rightarrow F'$ be functor morphism. We say that $\alpha$ is a \textbf{cwp-functor morphism} if the following diagram of functor morphisms is commutative.
 \[\xymatrix{
F\Phi_{\mathcal{A}}\ar[r]^{\sigma} \ar[d]_{\alpha\Phi_{\mathcal{A}}}&\Phi_{\mathcal{B}}F\ar[d]^{\Phi_{\mathcal{B}}\alpha}\\
F'\Phi_{\mathcal{A}}\ar[r]^{\sigma'}&\Phi_{\mathcal{B}}{F'}
 }\]
If $\alpha$ is a cwp-functor morphism, we write $\alpha:(F,\sigma)\rightarrow(F',\sigma').$
\end{dfn}

\vspace*{2mm}
A cwp-functor morphism $\alpha:(F,\sigma)\rightarrow(F',\sigma')$ induces a functor morphism
$$\mathfrak{F}(\alpha):\mathfrak{F}(F,{\sigma})\rightarrow \mathfrak{F}(F',{\sigma'})$$ 
defined by $$\mathfrak{F}(\alpha)(A):=(\alpha(A_1),\alpha(A_0))\in{\rm Hom}(F(A_1),F'(A_1))\oplus{\rm Hom}(F(A_0),F'(A_0)).$$
Since $\alpha$ is a cwp-functor morphism, the following diagram is commutative,
 \[\xymatrix{
F(A_1)\ar[rr]^{F(\varphi^A_1)}\ar[d]_{\alpha(A_1)}&&F(A_0)\ar[rr]^{\sigma(A_1)\circ F(\varphi^A_0)}\ar[d]_{\alpha(A_0)}&&\Phi_{\mathcal{B}}(F(A_1))\ar[d]^{\Phi_{\mathcal{B}}(\alpha(A_1))}\\
F'(A_1)\ar[rr]^{F'(\varphi^A_1)}&&F'(A_0)\ar[rr]^{\sigma'(A_1)\circ F'(\varphi^A_0)}&&\Phi_{\mathcal{B}}(F'(A_1)),
}\]
which means that $\mathfrak{F}(\alpha)(A)\in Z^0({\rm Hom}(F(A),F'(A)))$ for any $A\in\mathfrak{F}(\mathcal{A},\Phi_{\mathcal{A}},W_{\mathcal{A}})$.

\vspace*{2mm}
\begin{dfn}
Let $(F,\sigma), (F',\sigma'):(\mathcal{A},\Phi_{\mathcal{A}},W_{\mathcal{A}})\rightarrow(\mathcal{B},\Phi_{\mathcal{B}},W_{\mathcal{B}})$  be cwp-functors and let $\alpha:(F,\sigma)\rightarrow (F',\sigma')$ be a cwp-functor morphism. For a cwp-functor $(G,\mu):(\mathcal{C},\Phi_{\mathcal{C}},W_{\mathcal{C}})\rightarrow(\mathcal{A},\Phi_{\mathcal{A}},W_{\mathcal{A}})$, we define a cwp-functor morphism $$\alpha(G,\mu):(F\circ G,\mu F\circ G\sigma)\rightarrow(F'\circ G,\mu F'\circ G\sigma')$$ as $\alpha(G,\mu)(C):=\alpha(G(C))$ for any $C\in\mathcal{C}$.
Similarly, for a cwp-functor $(H,\nu):(\mathcal{B},\Phi_{\mathcal{B}},W_{\mathcal{B}})\rightarrow(\mathcal{C},\Phi_{\mathcal{C}},W_{\mathcal{C}})$, we define a cwp-functor morphism $$(H,\nu)\alpha:(H\circ F,\nu F\circ H\sigma)\rightarrow(H\circ F',\nu F'\circ H\sigma')$$ as $(H,\nu)\alpha(A):=H(\alpha(A))$ for any $A\in\mathcal{A}$.

\end{dfn}

\vspace{3mm}
The next lemma gives fundamental properties of functor morphisms given as $\mathfrak{F}(-)$. The proof is  left to the reader.

\vspace*{2mm}
\begin{lem}\label{cwp-functor morphism}
Let $(F,\sigma)$, $(F',\sigma')$ and $(F'',\sigma'')$ be cwp-functors from $(\mathcal{A},\Phi_{\mathcal{A}},W_{\mathcal{A}})$ to $(\mathcal{B},\Phi_{\mathcal{B}},W_{\mathcal{B}})$ and let  $\alpha:(F,\sigma)\rightarrow(F',\sigma')$ and $\beta:(F',\sigma')\rightarrow(F'',\sigma'')$ be cwp-functor morphisms. Then

\begin{itemize}
\item[$(1)$]$\beta\circ\alpha:(F,\sigma)\rightarrow(F'',\sigma'')$ is a cwp-functor morphism, and we have
$$\mathfrak{F}(\beta\circ\alpha)=\mathfrak{F}(\beta)\circ\mathfrak{F}(\alpha)$$

\item[$(2)$]If $\alpha$ is an isomorphism of functors, so is $\mathfrak{F}(\alpha)$.

\item[$(3)$] For a cwp-functor $(G,\mu):(\mathcal{C},\Phi_{\mathcal{C}},W_{\mathcal{C}})\rightarrow(\mathcal{A},\Phi_{\mathcal{A}},W_{\mathcal{A}})$, we have $$\mathfrak{F}(\alpha(G,\mu))=\mathfrak{F}(\alpha)\mathfrak{F}(G,\mu).$$

Similarly, for a cwp-functor $(H,\nu):(\mathcal{B},\Phi_{\mathcal{B}},W_{\mathcal{B}})\rightarrow(\mathcal{C},\Phi_{\mathcal{C}},W_{\mathcal{C}})$, we have 
$$\mathfrak{F}((H,\nu)\alpha)=\mathfrak{F}(H,\nu)\mathfrak{F}(\alpha).$$

\end{itemize}
\end{lem}

\vspace{3mm}
Next, we introduce the notion of  cwp-adjunction of cwp-functors.

\vspace*{2mm}
\begin{dfn}
Let $(F,\sigma):(\mathcal{A},\Phi_{\mathcal{A}},W_{\mathcal{A}})\rightarrow(\mathcal{B},\Phi_{\mathcal{B}},W_{\mathcal{B}})$ and $(G,\tau):(\mathcal{B},\Phi_{\mathcal{B}},W_{\mathcal{B}})\rightarrow(\mathcal{A},\Phi_{\mathcal{A}},W_{\mathcal{A}})$ be cwp-functors. We say that $(F,\sigma)$ is \textbf{left cwp-adjoint} to $(G,\tau)$, denoted by $(F,\sigma)\dashv(G,\tau)$, if $F$ is left adjoint to $G$ and adjunction morphisms are cwp-functor morphisms.
\end{dfn}

\vspace*{2mm}
\begin{lem}
In the same notation as above, assume $(F,\sigma)\dashv(G,\tau)$ and let $\varepsilon:(FG,\sigma G\circ F\tau)\rightarrow{\rm id}_{\mathcal{B}}$ and $\eta:{\rm id}_{\mathcal{A}}\rightarrow (GF,\tau F\circ G\sigma)$ be adjunction morphisms which are cwp-functor morphisms. Then $\mathfrak{F}(F,\sigma)\dashv\mathfrak{F}(G,\tau)$ and 
$$\mathfrak{F}(\varepsilon):\mathfrak{F}(F,\sigma)\circ\mathfrak{F}(G,\tau)\rightarrow{\rm id}_{\mathfrak{F}(\mathcal{B},\Phi_{\mathcal{B}},W_{\mathcal{B}})}$$ 
$$\mathfrak{F}(\eta):{\rm id}_{\mathfrak{F}(\mathcal{A},\Phi_{\mathcal{A}},W_{\mathcal{A}})}\rightarrow\mathfrak{F}(G,\tau)\circ{\mathfrak{F}(F,\sigma)}$$
are the adjunction morphisms of the adjoint pair $\mathfrak{F}(F,\sigma)\dashv\mathfrak{F}(G,\tau)$.
\begin{proof} Since $\varepsilon:FG\rightarrow{\rm id}_{\mathcal{B}}$ and $\eta:{\rm id}_{\mathcal{A}}\rightarrow GF$ are adjunction morphisms of the adjoint pair $F\dashv G$,  the following compositions  are identities of functors;
$$F\xrightarrow{F\eta}FGF\xrightarrow{\varepsilon F}F\hspace{3mm}{\rm and}\hspace{3mm}G\xrightarrow{\eta G}GFG\xrightarrow{G\varepsilon}G.$$
By Lemma \ref{cwp-functor lemma} and Lemma \ref{cwp-functor morphism}, the following compositions are also identities of dg-functors;
$$\mathfrak{F}(F,\sigma)\xrightarrow{\mathfrak{F}(F,\sigma)\mathfrak{F}(\eta)}\mathfrak{F}(F,\sigma)\mathfrak{F}(G,\tau)\mathfrak{F}(F,\sigma)\xrightarrow{\mathfrak{F}(\varepsilon) \mathfrak{F}(F,\sigma)}\mathfrak{F}(F,\sigma)$$
and
$$\mathfrak{F}(G,\tau)\xrightarrow{\mathfrak{F}(\eta) \mathfrak{F}(G,\tau)}\mathfrak{F}(G,\tau)\mathfrak{F}(F,\sigma)\mathfrak{F}(G,\tau)\xrightarrow{\mathfrak{F}(G,\tau)\mathfrak{F}(\varepsilon)}\mathfrak{F}(G,\tau).$$
Hence, we have an adjunction $\mathfrak{F}(F,\sigma)\dashv\mathfrak{F}(G,\tau)$, and $\mathfrak{F}(\varepsilon)$ and $\mathfrak{F}(\eta)$ are adjunction morphisms.
\end{proof}
\end{lem}

\vspace{2mm}
We give definitions of relative adjoint functors and basic properties of it after \cite{ulmer}.

\vspace*{2mm}
\begin{dfn}\label{rel ad}
Let $\mathcal{C}_1$, $\mathcal{C}_2$ and $\mathcal{D}$ be categories and let $F:\mathcal{C}_1\rightarrow\mathcal{D}$, $G:\mathcal{D}\rightarrow\mathcal{C}_2$ and $J:\mathcal{C}_1\rightarrow\mathcal{C}_2$ be functors. $F$ is called \textbf{left} $J$\textbf{-relative adjoint to} $G$ (or $J$\textbf{-left adjoint to} $G$) if for each $C\in\mathcal{C}_1, D\in\mathcal{D}$ there is an isomorphism
$${\rm Hom}_{\mathcal{D}}(F(C),D)\cong{\rm Hom}_{\mathcal{C}_2}(J(C),G(D))$$
which is functorial in $C$ and $D$.

Dually, $F$ is called \textbf{right} $J$\textbf{-relative adjoint to} $G$ (or $J$\textbf{-right adjoint to} $G$) if for each $C\in\mathcal{C}_1, D\in\mathcal{D}$ there is an isomorphism
$${\rm Hom}_{\mathcal{D}}(D, F(C))\cong{\rm Hom}_{\mathcal{C}_2}(G(D), J(C))$$
which is functorial in $C$ and $D$.
\end{dfn}

\vspace*{2mm}
\begin{rem}
\hfill\\
$(1)$ Relative adjointness is not symmetric property, i.e. although $F$ is  $J$-left adjoint to $G$, $G$ is not $J$-right adjoint to $F$ in general.\\
$(2)$ If $F$ is $J$-left adjoint to $G$, there is a functor morphism $$\mu:J\rightarrow GF$$ such that $\mu(C):J(C)\rightarrow G(F(C))$ corresponds to ${\rm id}_{F(C)}$.

Similarly, if $F$ is $J$-right adjoint to $G$, there is a functor morphism $$\nu:GF\rightarrow J$$ such that $\nu(C):G(F(C))\rightarrow J(C)$ is corresponding to ${\rm id}_{F(C)}$.

The above functor morphisms $\mu:J\rightarrow GF$ and $\nu:GF\rightarrow J$ are called \textbf{the front adjunction}.

\end{rem}

\vspace{2mm}
By the next lemma, we see that the existence of a front adjunction implies a relative adjunction.

\vspace*{2mm}
\begin{lem}[\cite{ulmer} Lemma 2.7]\label{rel ad lemma}
The notation is the same as in {\rm Definition \ref{rel ad}}.  The functor $F$ is $J$-left adjoint to $G$ if and only if there exists a functor morphism $\mu:J\rightarrow GF$  such that for each $C\in\mathcal{C}_1$ and $D\in\mathcal{D}$
the composition of maps
$${\rm Hom}(F(C),D)\xrightarrow{G(-)}{\rm Hom}(G(F(C)),G(D))\xrightarrow{{\rm Hom}(\mu(C),G(D))}{\rm Hom}(J(C),G(D))$$
is a bijection.

Similarly, $F$ is $J$-right adjoint to $G$ if and only if there exists a functor morphism $\nu:GF\rightarrow J$ such that for each $C\in\mathcal{C}_1$ and $D\in\mathcal{D}$
the composition of maps
$${\rm Hom}(D,F(C))\xrightarrow{G(-)}{\rm Hom}(G(D),G(F(C)))\xrightarrow{{\rm Hom}(G(D),\nu(C))}{\rm Hom}(G(D),J(C))$$
is bijective.
\end{lem}

\vspace{2mm}
Similarly, the notion of relative cwp-adjunction is  given in the following.

\vspace*{2mm}
\begin{dfn}\label{rel cwp-ad}
In the same notation as in {\rm Definition \ref{rel ad}}, let $(\Phi_i,W_i)$ and $(\Psi,V)$ be potentials of $\mathcal{C}_i$ and $\mathcal{D}$ respectively. Let $(F,\sigma):(\mathcal{C}_1,\Phi_1,W_1)\rightarrow(\mathcal{D},\Psi,V)$, $(G,\tau):(\mathcal{D},\Psi,V)\rightarrow\mathcal(\mathcal{C}_2,\Phi_2,W_2)$ and $(J,\eta):(\mathcal{C}_1,\Phi_1,W_1)\rightarrow(\mathcal{C}_2,\Phi_2,W_2)$ be cwp-functors. $(F,\sigma)$ is called $(J,\eta)$\textbf{-left cwp-adjoint to} $(G,\tau)$ if $F$ is $J$-left adjoint to $G$ and the front adjunction is cwp-functor morphism.

Dually, we say $(F,\sigma)$ is $(J,\eta)$\textbf{-right cwp-adjoint to} $(G,\tau)$ if $F$ is $J$-right adjoint to $G$ and the front adjunction is cwp-functor morphism.

\end{dfn}

\vspace*{2mm}
\begin{lem}
Notation is the same as in {\rm Definition \ref{rel cwp-ad}}. If $(F,\sigma)$ is  $(J,\eta)$-left cwp-adjoint to $(G,\tau)$ and $\mu: J\rightarrow GF$ is the front adjunction, then $\mathfrak{F}(F,\sigma)$ is  $\mathfrak{F}(J,\eta)$-left adjoint to $\mathfrak{F}(G,\tau)$ and the front adjunction is  $\mathfrak{F}(\mu): \mathfrak{F}(J,\eta)\rightarrow \mathfrak{F}(G,\tau)\mathfrak{F}(F,\sigma)$. 

Similarly, if $(F,\sigma)$ is  $(J,\eta)$-right cwp-adjoint to $(G,\tau)$ and $\nu$ is the front adjunction, then $\mathfrak{F}(F,\sigma)$ is  $\mathfrak{F}(J,\eta)$-right adjoint to $\mathfrak{F}(G,\tau)$ and the front adjunction is  $\mathfrak{F}(\nu)$.

\begin{proof}
 If $(F,\sigma)$ is  $(J,\eta)$-left cwp-adjoint to $(G,\tau)$, then the front adjunction  $\mu: J\rightarrow GF$ is cwp-functor morphism, and the composition $${\rm Hom}(F(C),D)\xrightarrow{G(-)}{\rm Hom}(G(F(C)),G(D))\xrightarrow{{\rm Hom}(\mu(C),G(D))}{\rm Hom}(J(C),G(D))$$
is a bijection. Hence, the composition of morphisms
$${\rm Hom}(\mathfrak{F}(F,\sigma)(C),D)\xrightarrow{\mathfrak{F}(G,\tau)(-)}{\rm Hom}(\{\mathfrak{F}(G,\tau)\circ\mathfrak{F}(F,\sigma)\}(C),\mathfrak{F}(G,\tau)(D))$$
and 
$${\rm Hom}(\{\mathfrak{F}(G,\tau)\circ\mathfrak{F}(F,\sigma)\}(C),\mathfrak{F}(G,\tau)(D))\xrightarrow{{\rm Hom}(\mathfrak{F}(\mu)(C),\mathfrak{F}(G,\tau)(D))}{\rm Hom}(\mathfrak{F}(J,\eta)(C),\mathfrak{F}(G,\tau)(D))$$
is also bijective. By Lemma \ref{rel ad lemma}, we see that $\mathfrak{F}(F,\sigma)$ is  $\mathfrak{F}(J,\eta)$-left adjoint to $\mathfrak{F}(G,\tau)$, and the front adjunction is  $\mathfrak{F}(\mu): \mathfrak{F}(J,\eta)\rightarrow \mathfrak{F}(G,\tau)\mathfrak{F}(F,\sigma)$.

The latter statement can be proved in a similar way.
\end{proof}

\end{lem}

\vspace{2mm}

In what follows, we define cwp-bifunctors.

\vspace*{2mm}
\begin{dfn}\label{cwp-cobifunctor} Let $P:\mathcal{A}\times\mathcal{B}\rightarrow\mathcal{C}$ be a bifunctor. We say that $P$ is\textbf{ compatible with potentials} with respect to $(\Phi_{\mathcal{A}},W_{\mathcal{A}})$, $(\Phi_{\mathcal{B}},W_{\mathcal{B}})$ and $(\Phi_{\mathcal{C}},W_{\mathcal{C}})$ if there are bifunctor isomorphisms $\sigma_{\mathcal{A}}:P(\Phi_{\mathcal{A}}\times{\rm id}_{\mathcal{B}})\xrightarrow{\sim}\Phi_{\mathcal{C}}P$ and $\sigma_{\mathcal{B}}:P({\rm id}_{\mathcal{A}}\times\Phi_{\mathcal{B}})\xrightarrow{\sim}\Phi_{\mathcal{C}}P$
such that  
$$\sigma_{\mathcal{A}}(A,B)\circ P(W_{\mathcal{A}}(A),B)+\sigma_{\mathcal{B}}(A,B)\circ P(A,W_{\mathcal{B}}(B))=W_{\mathcal{C}}(P(A,B))$$
and 
$$\Phi_{\mathcal{C}}(\sigma_{\mathcal{B}}(A,B))\circ\sigma_{\mathcal{A}}(A,\Phi_{\mathcal{B}}(B))=\Phi_{\mathcal{C}}(\sigma_{\mathcal{A}}(A,B))\circ\sigma_{\mathcal{B}}(\Phi_{\mathcal{A}}(A),B)$$
  for any $A\in\mathcal{A}$ and $B\in\mathcal{B}$. By the latter equation above, $\sigma_{\mathcal{A}}$ and $\sigma_{\mathcal{B}}$ induce a natural functor isomorphism $\sigma^{m,n}:P(\Phi_{\mathcal{A}}^m\times\Phi_{\mathcal{B}}^n)\xrightarrow{\sim}\Phi^{m+n}P$ for any $m,n\in\mathbb{Z}$.
 The triple $(P,\sigma_{\mathcal{A}},\sigma_{\mathcal{B}})$ is called \textbf{cwp-bifunctor} and we write $$(P,\sigma_{\mathcal{A}},\sigma_{\mathcal{B}}):(\mathcal{A},\Phi_{\mathcal{A}},W_{\mathcal{A}})\times(\mathcal{B},\Phi_{\mathcal{B}},W_{\mathcal{B}})\rightarrow(\mathcal{C},\Phi_{\mathcal{C}},W_{\mathcal{C}}).$$
\end{dfn}

\vspace{7mm}For a cwp-bifunctor $(P,\sigma_{\mathcal{A}},\sigma_{\mathcal{B}}):(\mathcal{A},\Phi_{\mathcal{A}},W_{\mathcal{A}})\times(\mathcal{B},\Phi_{\mathcal{B}},W_{\mathcal{B}})\rightarrow(\mathcal{C},\Phi_{\mathcal{C}},W_{\mathcal{C}})$, we define a dg-bifunctor 
$$\mathfrak{F}(P,\sigma_{\mathcal{A}},\sigma_{\mathcal{B}}):\mathfrak{F}(\mathcal{A},\Phi_{\mathcal{A}},W_{\mathcal{A}})\otimes\mathfrak{F}(\mathcal{B},\Phi_{\mathcal{B}},W_{\mathcal{B}})\rightarrow\mathfrak{F}(\mathcal{C},\Phi_{\mathcal{C}},W_{\mathcal{C}})$$
as follows. For any object $A=(A_1\xrightarrow{\varphi^A_1}A_0\xrightarrow{\varphi^A_0}\Phi_{\mathcal{A}}(A_1))\in\mathfrak{F}(\mathcal{A},\Phi_{\mathcal{A}},W_{\mathcal{A}})$ and $B=(B_1\xrightarrow{\psi^B_1}B_0\xrightarrow{\psi^B_0}\Phi_{\mathcal{B}}(B_1))\in\mathfrak{F}(\mathcal{B},\Phi_{\mathcal{B}},W_{\mathcal{B}})$, we define the object $\mathfrak{F}(P,\sigma_{\mathcal{A}},\sigma_{\mathcal{B}})(A,B)\in\mathfrak{F}(\mathcal{C},\Phi_{\mathcal{C}},W_{\mathcal{C}})$ as 

$$\Bigl(P(A_1,B_0)\oplus P(A_0,B_1)\xrightarrow{\omega_1}P(A_0,B_0)\oplus\Phi_{\mathcal{C}}(P(A_1,B_1))\xrightarrow{\omega_0}\Phi_{\mathcal{C}}(P(A_1,B_0))\oplus\Phi_{\mathcal{C}}(P(A_0,B_1))\Bigr),$$
where 
$${\omega_1}=\begin{pmatrix}
P(\varphi^A_1,{\rm id})&P({\rm id},\psi^B_1)\\ -\sigma_{\mathcal{B}}(A_1,B_1)\circ P({\rm id},\psi^B_0)&\sigma_{\mathcal{A}}(A_1,B_1)\circ P(\varphi^A_0,{\rm id})
\end{pmatrix}
$$
and 

$$\omega_0=\begin{pmatrix}
\sigma_{\mathcal{A}}(A_1,B_0)\circ P(\varphi^A_0,{\rm id}) & -\Phi_{\mathcal{C}}(P({\rm id},\psi^B_1))\\
\sigma_{\mathcal{B}}(A_0,B_1)\circ P({\rm id},\psi^B_0) & \Phi_{\mathcal{C}}(P(\varphi^A_1,{\rm id}))
\end{pmatrix}.$$\\
For a morphism $f:(A,B)\rightarrow (A',B')$ in $\mathfrak{F}(\mathcal{A},\Phi_{\mathcal{A}},W_{\mathcal{A}})\otimes\mathfrak{F}(\mathcal{B},\Phi_{\mathcal{B}},W_{\mathcal{B}})$, we define  the morphism $\mathfrak{F}(P,\sigma_{\mathcal{A}},\sigma_{\mathcal{B}})(f):\mathfrak{F}(P,\sigma_{\mathcal{A}},\sigma_{\mathcal{B}})(A,B)\rightarrow\mathfrak{F}(P,\sigma_{\mathcal{A}},\sigma_{\mathcal{B}})(A',B')$ by the following rule, 
$$\mathfrak{F}(P,\sigma_{\mathcal{A}},\sigma_{\mathcal{B}})(g_{i,j}^m\otimes h_{k,l}^n):=\begin{cases}(-1)^{{\rm deg}(h_{1,l}^n)}\Phi_{\mathcal{C}}(\sigma^{m,n}(A'_j,B'_l)\circ P(g_{1,j}^m,h_{1,l}^n)) & {\rm if}\hspace{2mm} i=k=1\\
(-1)^{i{\rm deg}(h_{k,l}^n)}\sigma^{m,n}(A'_j,B'_l)\circ P(g_{i,j}^m,h_{k,l}^n) & {\rm otherwise}
\end{cases},$$
where $g_{i,j}^m\in{\rm Hom}_{\mathcal{A}}(A_i,\Phi_{\mathcal{A}}^{m}(A'_j))$ and $h_{k,l}^n\in{\rm Hom}_{\mathcal{B}}(B_k, \Phi_{\mathcal{B}}^{n}(B'_l))$.


\vspace*{2mm}
\begin{dfn}\label{cwp-contrabifunctor} Let $Q:\mathcal{A}^{\rm op}\times\mathcal{B}\rightarrow\mathcal{C}$ be a bifunctor. We say that $Q$ is\textbf{ compatible with potentials} with respect to $(\Phi_{\mathcal{A}},W_{\mathcal{A}})$, $(\Phi_{\mathcal{B}},W_{\mathcal{B}})$ and $(\Phi_{\mathcal{C}},W_{\mathcal{C}})$ if there are bifunctor isomorphisms $\tau_{\mathcal{A}}:Q((\Phi_{\mathcal{A}}^{\rm op})^{-1}\times{\rm id}_{\mathcal{B}})\xrightarrow{\sim}\Phi_{\mathcal{C}}Q$ and $\tau_{\mathcal{B}}:Q({\rm id}_{\mathcal{A}}\times\Phi_{\mathcal{B}})\xrightarrow{\sim}\Phi_{\mathcal{C}}Q$
such that  
$$-\tau_{\mathcal{A}}(A,B)\circ Q((\Phi_{\mathcal{A}}^{\rm op})^{-1}(W_{\mathcal{A}}^{\rm op}(A)),B)+\tau_{\mathcal{B}}(A,B)\circ Q(A,W_{\mathcal{B}}(B))=W_{\mathcal{C}}(Q(A,B))$$
and 
$$\Phi_{\mathcal{C}}(\tau_{\mathcal{B}}(A,B))\circ\tau_{\mathcal{A}}(A,\Phi_{\mathcal{B}}(B))=\Phi_{\mathcal{C}}(\tau_{\mathcal{A}}(A,B))\circ\tau_{\mathcal{B}}((\Phi_{\mathcal{A}}^{\rm op})^{-1}(A),B)$$
  for any $A\in\mathcal{A}$ and $B\in\mathcal{B}$, where $\Phi_{\mathcal{A}}^{\rm op}:\mathcal{A}^{\rm op}\xrightarrow{\sim}\mathcal{A}^{\rm op}$ is the opposite equivalence of $\Phi_{\mathcal{A}}$ and $W_{\mathcal{A}}^{\rm op}:\Phi_{\mathcal{A}}^{\rm op}\rightarrow{\rm id}_{\mathcal{A}^{\rm op}}$ is the opposite functor morphism of $W_{\mathcal{A}}$. By the latter equation above, $\tau_{\mathcal{A}}$ and $\tau_{\mathcal{B}}$ induce a natural functor isomorphism $\tau^{m,n}:Q(\Phi_{\mathcal{A}}^m\times\Phi_{\mathcal{B}}^n)\xrightarrow{\sim}\Phi^{-m+n}Q$ for any $m,n\in\mathbb{Z}$.
 The triple $(Q,\tau_{\mathcal{A}},\tau_{\mathcal{B}})$ is called \textbf{cwp-bifunctor} and we write $$(Q,\tau_{\mathcal{A}},\tau_{\mathcal{B}}):(\mathcal{A},\Phi_{\mathcal{A}},W_{\mathcal{A}})^{\rm op}\times(\mathcal{B},\Phi_{\mathcal{B}},W_{\mathcal{B}})\rightarrow(\mathcal{C},\Phi_{\mathcal{C}},W_{\mathcal{C}}).$$
\end{dfn}

\vspace{7mm}For a cwp-bifunctor $(Q,\tau_{\mathcal{A}},\tau_{\mathcal{B}}):(\mathcal{A},\Phi_{\mathcal{A}},W_{\mathcal{A}})^{\rm op}\times(\mathcal{B},\Phi_{\mathcal{B}},W_{\mathcal{B}})\rightarrow(\mathcal{C},\Phi_{\mathcal{C}},W_{\mathcal{C}})$, we define a dg-bifunctor 
$$\mathfrak{F}(Q,\tau_{\mathcal{A}},\tau_{\mathcal{B}}):\mathfrak{F}(\mathcal{A},\Phi_{\mathcal{A}},W_{\mathcal{A}})^{\rm op}\otimes\mathfrak{F}(\mathcal{B},\Phi_{\mathcal{B}},W_{\mathcal{B}})\rightarrow\mathfrak{F}(\mathcal{C},\Phi_{\mathcal{C}},W_{\mathcal{C}})$$
as follows. For any object $A=(A_1\xrightarrow{\varphi^A_1}A_0\xrightarrow{\varphi^A_0}\Phi_{\mathcal{A}}(A_1))\in\mathfrak{F}(\mathcal{A},\Phi_{\mathcal{A}},W_{\mathcal{A}})^{\rm op}$ and $B=(B_1\xrightarrow{\psi^B_1}B_0\xrightarrow{\psi^B_0}\Phi_{\mathcal{B}}(B_1))\in\mathfrak{F}(\mathcal{B},\Phi_{\mathcal{B}},W_{\mathcal{B}})$, we define the object $\mathfrak{F}(Q,\tau_{\mathcal{A}},\tau_{\mathcal{B}})(A,B)\in\mathfrak{F}(\mathcal{C},\Phi_{\mathcal{C}},W_{\mathcal{C}})$ as 

$$\Bigl(\Phi_{\mathcal{C}}^{-1}(Q(A_1,B_0))\oplus Q(A_0,B_1)\xrightarrow{\omega_1}Q(A_0,B_0)\oplus Q(A_1,B_1))\xrightarrow{\omega_0}Q(A_1,B_0)\oplus\Phi_{\mathcal{C}}(Q(A_0,B_1))\Bigr),$$
where 
$${\omega_1}=\begin{pmatrix}
Q(\varphi^A_0,{\rm id})\circ((\tau^{1,0})(A_1,B_0))^{-1}&Q({\rm id},\psi^B_1)\\ \Phi_{\mathcal{C}}^{-1}(\tau_{\mathcal{B}}(A_1,B_1)\circ Q({\rm id},\psi^B_0))&Q(\varphi^A_1,{\rm id})
\end{pmatrix}
$$
and 

$$\omega_0=\begin{pmatrix}
-Q(\varphi^A_1,{\rm id}) & Q({\rm id},\psi^B_1)\\
\tau_{\mathcal{B}}(A_0,B_1)\circ Q({\rm id},\psi^B_0) & -\tau_{\mathcal{A}}(A_0,B_1)\circ(Q(\Phi_{A}^{-1}(\varphi^A_0),{\rm id}))
\end{pmatrix}.$$\\
For a morphism $f:(A,B)\rightarrow (A',B')$ in $\mathfrak{F}(\mathcal{A},\Phi_{\mathcal{A}},W_{\mathcal{A}})^{\rm op}\otimes\mathfrak{F}(\mathcal{B},\Phi_{\mathcal{B}},W_{\mathcal{B}})$, we define  the morphism $\mathfrak{F}(Q,\tau_{\mathcal{A}},\tau_{\mathcal{B}})(f):\mathfrak{F}(Q,\tau_{\mathcal{A}},\tau_{\mathcal{B}})(A,B)\rightarrow\mathfrak{F}(Q,\tau_{\mathcal{A}},\tau_{\mathcal{B}})(A',B')$ by the following rule, 

\vspace{3mm}
$\mathfrak{F}(Q,\tau_{\mathcal{A}},\tau_{\mathcal{B}})(g_{i,j}^m\otimes h_{k,l}^n)\\
\hspace{6mm}:=\begin{cases}\Phi_{\mathcal{C}}^{m-1}(\tau^{0,n}(A'_j,B'_l)\circ Q(g_{1,j}^m,h_{0,l}^n)\circ \tau^{m,0}(A_1,B_0)^{-1}) & {\rm if}\hspace{2mm} i=1\hspace{2mm} k=0\\
(-1)^{i-k+1}\Phi_{\mathcal{C}}^{m}(\tau^{0,n}(A'_j,B'_l)\circ Q(g_{i,j}^m,h_{k,l}^n)\circ \tau^{m,0}(A_i,B_k)^{-1}) & {\rm otherwise}
\end{cases}$

where $g_{i,j}^m\in{\rm Hom}_{\mathcal{A}^{\rm op}}(\Phi_{\mathcal{A}}^{m}(A_i),A_j')$ and $h_{k,l}^n\in{\rm Hom}(B_k,\Phi_{\mathcal{B}}^{n}(B'_l))$.

\vspace{3mm}
\subsection{ind/pro-categories and their factorization categories}

In this section, we recall the notion of ind-categories and pro-categories, and study factorization categories of ind/pro-categories. For the detail of ind/pro-categories, see \cite{cp} or \cite{kas}, for example.

At first, we recall the definition and the foundations of ind/pro-categories.

\begin{dfn}
A small category $\mathcal{I}$ is called \textbf{filtering} if  the following properties hold;
\begin{itemize}
\item[$(1)$] For any objects $i,i'\in \mathcal{I}$, there exist an object $j\in\mathcal{I}$ and morphisms $i\rightarrow j$ and $i'\rightarrow j$.
\item[$(2)$] For two morphisms $u,v:k'\rightarrow k$ in $\mathcal{I}$, there exist an object $l\in\mathcal{I}$ and a morphism $w:k\rightarrow l$ such that $w\circ u=w\circ v$.
\end{itemize}

A small category $\mathcal{J}$ is called \textbf{cofiltering} if its opposite category $\mathcal{J}^{\rm op}$ is filtering.

\end{dfn}

\vspace{2mm}
\begin{dfn}
Let $\mathcal{C}$ be a category. \\
(1) We define the \textbf{ind-category} of $\mathcal{C}$, denoted by ${\rm Ind}(\mathcal{C})$, as follows:

An object of  ${\rm Ind}(\mathcal{C})$ is a functor $D:\mathcal{I}\rightarrow\mathcal{C}$ with $\mathcal{I}$ filtering. For two objects $D:\mathcal{I}\rightarrow\mathcal{C}$ and $E:\mathcal{J}\rightarrow\mathcal{C}$, we define the set of morphisms as 
 \begin{eqnarray*}
 {\rm Hom}_{{\rm Ind}(\mathcal{C})}(D,E):=\varprojlim_{i\in\mathcal{I}}\varinjlim_{j\in\mathcal{J}}{\rm Hom}_{\mathcal{C}}(D(i),E(j)).
 \end{eqnarray*}
 (2) We define the  \textbf{pro-category} of $\mathcal{C}$, denoted by ${\rm Pro}(\mathcal{C})$, by the following:
 
 An object of ${\rm Pro}(\mathcal{C})$ is a functor $P:\mathcal{I}\rightarrow\mathcal{C}$ with $\mathcal{I}$ cofiltering. For two objects $P:\mathcal{I}\rightarrow\mathcal{C}$ and $Q:\mathcal{J}\rightarrow\mathcal{C}$, we define the space of morphisms as 
 \begin{eqnarray*}
 {\rm Hom}_{{\rm Pro}(\mathcal{C})}(P,Q):=\varprojlim_{j\in\mathcal{J}}\varinjlim_{i\in\mathcal{I}}{\rm Hom}_{\mathcal{C}}(P(i),Q(j)).
 \end{eqnarray*}

\end{dfn}

\vspace{2mm}
\begin{rem}
$(1)$ We have a natural equivalence $${\rm Pro}(\mathcal{C})\cong{\rm Ind}(\mathcal{C}^{\rm op})^{\rm op}.$$
$(2)$ Let $D:\mathcal{I}\rightarrow\mathcal{C}$ and $E:\mathcal{J}\rightarrow\mathcal{C}$ be objects of ${\rm Ind}(\mathcal{C})$. The set of morphisms $ {\rm Hom}_{{\rm Ind}(\mathcal{C})}(D,E)$ is interpreted as the set of equivalence classes of maps of systems defined as follows:

A \textbf{map of systems} from $D$ to $E$ is a pair $\varphi=(\{\varphi_i\}_{i\in\mathcal{I}}, \theta_{\varphi})$ where $\theta_{\varphi}:{\rm Ob}(\mathcal{I})\rightarrow{\rm Ob}(\mathcal{J})$ is a map from ${\rm Ob}(\mathcal{I})$ to ${\rm Ob}(\mathcal{J})$, and $\varphi_i\in{\rm Hom}_{\mathcal{C}}(D(i), E(\theta_{\varphi}(i)))$, such that  for any morphism $v:i\rightarrow i'$ in $\mathcal{I}$ there are $j\in\mathcal{J}$, $u: \theta_{\varphi}(i)\rightarrow j$ and $u': \theta_{\varphi}(i')\rightarrow j$ such that the following diagram is commutative:
 \[\xymatrix{
D(i)\ar[r]^{\varphi_{i}}\ar[dd]_{D(v)}&E(\theta_{\varphi}(i))\ar[rd]^{E(u)}&\\
&&E(j)\\
D(i') \ar[r]^{\varphi_{i'}}&E(\theta_{\varphi}(i'))\ar[ru]^{E(u')}
  }\]
Two maps of systems $\varphi=(\{\varphi_i\}_{i\in\mathcal{I}}, \theta_{\varphi})$ and $\psi=(\{\psi_i\}_{i\in\mathcal{I}}, \theta_{\psi})$ are \textbf{equivalent} if for each $i\in\mathcal{I}$, there exist $j\in{\rm Ob}(\mathcal{J})$, $u:\theta_{\varphi}(i)\rightarrow j$ and $v:\theta_{\psi}(i)\rightarrow j$ such that the following diagram commutes:
 \[\xymatrix{
&E(\theta_{\varphi}(i))\ar[rd]^{E(u)}&\\
D(i)\ar[ru]^{\varphi_i} \ar[rd]_{\psi_i}&&E(j)\\
&E(\theta_{\psi}(i)) \ar[ru]^{E(v)}&
  }\]
We denote by $[\varphi]$ the morphism from $D$ to $E$ in ${\rm Ind}(\mathcal{C})$ corresponding to the equivalence class of a map of systems $\varphi$. With this notation one can easily write down the composition of $[\varphi]\in{\rm Hom}(D,E)$ and $[\psi]\in{\rm Hom}(E,H)$, where $H:\mathcal{K}\rightarrow\mathcal{C}$. The composition is given by 
$$[\psi]\circ[\varphi]=[(\{\psi_{\theta_{\varphi}(i)}\circ\varphi_i\}_{i\in\mathcal{I}},\theta_{\psi}\circ\theta_{\varphi})].$$  
\\
$(3)$ Let $P:\mathcal{I}\rightarrow\mathcal{C}$ and $Q:\mathcal{J}\rightarrow\mathcal{C}$ be objects of ${\rm Pro}(\mathcal{C})$. Similarly to $(2)$, the set  ${\rm Hom}_{{\rm Pro}(\mathcal{C})}(P,Q)$ is interpreted as the set of equivalence classes of maps of systems  defined by the following:

A \textbf{map of systems} from $P$ to $Q$ is a pair $\varphi=(\{\varphi_j\}_{j\in\mathcal{J}}, \theta_{\varphi})$ where $\theta_{\varphi}:{\rm Ob}(\mathcal{J})\rightarrow{\rm Ob}(\mathcal{I})$ is a map from ${\rm Ob}(\mathcal{J})$ to ${\rm Ob}(\mathcal{I})$ and $\varphi_j\in{\rm Hom}_{\mathcal{C}}(P( \theta_{\varphi}(j)), Q(j))$, such that  for any morphism $v:j\rightarrow j'$ in $\mathcal{J}$ there are $i\in\mathcal{I}$, $u:i\rightarrow \theta_{\varphi}(j)$ and $u':i\rightarrow \theta_{\varphi}(j')$ such that the following diagram is commutative:
 \[\xymatrix{
&P(\theta_{\varphi}(j'))\ar[r]^{\varphi_{j'}}&Q(j')\\
P(i)\ar[ru]^{P(u')} \ar[rd]_{P(u)}&&\\
&P(\theta_{\varphi}(j)) \ar[r]^{\varphi_j}&Q(j)\ar[uu]_{Q(v)}
  }\]
Two maps of systems $\varphi=(\{\varphi_j\}_{j\in\mathcal{J}}, \theta_{\varphi})$ and $\psi=(\{\psi_j\}_{j\in\mathcal{J}}, \theta_{\psi})$ are \textbf{equivalent} if for each $j\in\mathcal{J}$, there exist $i\in{\rm Ob}(\mathcal{I})$, $u:i\rightarrow\theta_{\varphi}(j)$ and $v:i\rightarrow\theta_{\psi}(j)$ such that the following diagram commutes:
 \[\xymatrix{
&P(\theta_{\varphi}(j))\ar[rd]^{\varphi_j}&\\ 
P(i)\ar[ru]^{P(u)} \ar[rd]_{P(v)}&&Q(j)\\
&P(\theta_{\psi}(j)) \ar[ru]^{\psi_j}&
  }\]
We denote by $[\varphi]$ the morphism from $P$ to $Q$ in ${\rm Pro}(\mathcal{C})$ corresponding to the equivalence class of a map of systems $\varphi$. Let  $[\varphi]\in{\rm Hom}(P,Q)$ and $[\psi]\in{\rm Hom}(Q,R)$ be morphisms, where $R:\mathcal{K}\rightarrow\mathcal{C}$. The composition of $[\varphi]$ with $[\psi]$ is given by 
$$[\psi]\circ[\varphi]=[(\{\psi_k\circ\varphi_{\theta_{\psi}(k)}\}_{k\in\mathcal{K}},\theta_{\varphi}\circ\theta_{\psi})].$$  

\end{rem}


\vspace{2mm}
\begin{dfn}
 For $C\in\mathcal{C}$,  $\iota(C):\widetilde{\{1\}}\rightarrow\mathcal{C}$ is the functor from the category 
 $\widetilde{\{1\}}$ with a unique object, $1$, and a unique morphism, ${\rm id}_1$, defined by $\iota(C)(1):=C$. $\iota(-)$ defines  natural functors
 $$\iota_{\rm Ind}:\mathcal{C}\rightarrow{\rm Ind}(\mathcal{C})$$
 $$\iota_{\rm Pro}:\mathcal{C}\rightarrow{\rm Pro}(\mathcal{C}).$$
By the constructions, the functors $\iota_{\rm Ind}$ and $\iota_{\rm Pro}$ are  fully faithful. 
\end{dfn}

\vspace{2mm}
\begin{rem} 
 ${\rm Ind}(-)$ defines an endofunctor on the category of functors, i.e.\\
$(a)$ A functor $F:\mathcal{C}\rightarrow\mathcal{D}$  induces  a natural functor 
$${\rm Ind}(F):{\rm Ind}(\mathcal{C})\rightarrow{\rm Ind}(\mathcal{D})$$
as follows: 
For an object $D:\mathcal{I}\rightarrow\mathcal{C}\in{\rm Ind}(\mathcal{C})$, the object ${\rm Ind}(F)(D)$ is defined by $F\circ D:\mathcal{I}\rightarrow\mathcal{D}$. For another object $D':\mathcal{I}'\rightarrow\mathcal{C}$ and for a morphism $[\varphi]:D\rightarrow D'$,  ${\rm Ind}(F)([\varphi])$ is defined by $[(\{F(\varphi_{i})\}_{i\in\mathcal{I}},\theta_{\varphi})]$. The following diagram is commutative.
\[\xymatrix{
{\rm Ind}(\mathcal{C})\ar[rr]^{{\rm Ind}(F)}&&{\rm Ind}(\mathcal{D})\\
\mathcal{C}\ar[rr]^{F} \ar[u]^{\iota_{\rm Ind}}&&\mathcal{D}\ar[u]_{\iota_{\rm Ind}}
  }\]
  $(b)$ Let $F,G:\mathcal{C}\rightarrow\mathcal{D}$ be functors. A functor morphism $\alpha:F\rightarrow G$  induces a natural functor morphism $${\rm Ind}(\alpha):{\rm Ind}(F)\rightarrow {\rm Ind}(G)$$
  as follows:
 For an object $D:\mathcal{I}\rightarrow\mathcal{C}\in{\rm Ind}(\mathcal{C})$, the morphism ${\rm Ind}(\alpha)(D)$ is defined by $[(\{\alpha D(i)\}_{i\in\mathcal{I}},{\rm id}_{\mathcal{I}})]$.

\vspace{1mm}
 Similarly, ${\rm Pro}(-)$ defines an endofunctor on the category of functors, i.e.\\
$(a')$ A functor $F:\mathcal{C}\rightarrow\mathcal{D}$ induces a natural functor 
$${\rm Pro}(F):{\rm Pro}(\mathcal{C})\rightarrow{\rm Pro}(\mathcal{D})$$
as follows:
For an object $P:\mathcal{I}\rightarrow\mathcal{C}\in{\rm Pro}(\mathcal{C})$, the object ${\rm Pro}(F)(P)$ is defined by $(F\circ P:\mathcal{I}\rightarrow\mathcal{D})$. For another object $P':\mathcal{I}'\rightarrow\mathcal{C}$ and for a morphism $[\varphi]:P\rightarrow P'$,  ${\rm Pro}(F)([\varphi])$ is defined by $[((F(\varphi_{i'}))_{i'\in\mathcal{I}'},\theta_{\varphi})]$. The following diagram is commutative.
\[\xymatrix{
{\rm Pro}(\mathcal{C})\ar[rr]^{{\rm Pro}(F)}&&{\rm Pro}(\mathcal{D})\\
\mathcal{C}\ar[rr]^{F} \ar[u]^{\iota}&&\mathcal{D}\ar[u]_{\iota}
  }\]
  $(b')$ Let $F,G:\mathcal{C}\rightarrow\mathcal{D}$ be functors. A functor morphism $\alpha:F\rightarrow G$  induces a natural functor morphism $${\rm Pro}(\alpha):{\rm Pro}(F)\rightarrow {\rm Pro}(G)$$
as follows:  
 For an object $P:\mathcal{I}\rightarrow\mathcal{C}\in{\rm Pro}(\mathcal{C})$, the morphism ${\rm Pro}(\alpha)(P)$ is defined by $[(\alpha P(i))_{i\in\mathcal{I}},{\rm id}_{\mathcal{I}})]$.
\end{rem}

\vspace{8mm}

\begin{prop} We have the following:
\begin{itemize}
\item[$(1)$] If $\mathcal{C}$ is an abelian category, then the categories ${\rm Ind}(\mathcal{C})$ and ${\rm Pro}(\mathcal{C})$ are abelian categories.
\item[$(2)$] If $\mathcal{E}$ is an exact category, then the categories ${\rm Ind}(\mathcal{E})$ and ${\rm Pro}(\mathcal{E})$ are exact categories.
\item[$(3)$] If $F:\mathcal{A}\rightarrow\mathcal{B}$ is an exact functor of exact categories, then the functors  ${\rm Ind}(F):{\rm Ind}(\mathcal{A})\rightarrow{\rm Ind}(\mathcal{B})$ and ${\rm Pro}(F):{\rm Pro}(\mathcal{A})\rightarrow{\rm Pro}(\mathcal{B})$ are exact functors.
\end{itemize}
\begin{proof}
(1) This  follows from \cite[Theorem 8.6.5.]{kas}\\
(2) This is \cite[Proposition 4.18.]{previdi}\\
(3) Since we can take abelian envelopes of the exact categories $\mathcal{A}$ and $\mathcal{B}$, and extend the functor $F$ to a functor between the abelian envelopes (see the proof of Proposition \ref{Z^0exact}), we may assume that $\mathcal{A}$ and $\mathcal{B}$ are abelian categories. Then we obtain the result by \cite[Corollary 8.6.8.]{kas} 
\end{proof}
\end{prop}

\vspace{2mm}

Let $(\mathcal{A},\Phi_\mathcal{A},W_\mathcal{A})$ be a category with a potential. Then 
$${\rm Ind}(\mathcal{A},\Phi_\mathcal{A},W_\mathcal{A}):=({\rm Ind}(\mathcal{A}),{\rm Ind}(\Phi_\mathcal{A}),{\rm Ind}(W))$$
$${\rm Pro}(\mathcal{A},\Phi_\mathcal{A},W_\mathcal{A}):=({\rm Pro}(\mathcal{A}),{\rm Pro}(\Phi_\mathcal{A}),{\rm Pro}(W))$$
are  categories with potentials. Since the natural functor $\iota_{\rm Ind}:\mathcal{A}\rightarrow{\rm Ind}(\mathcal{A})$ (resp. $\iota_{\rm Pro}:\mathcal{A}\rightarrow{\rm Pro}(\mathcal{A})$) is compatible with potentials with respect to $(\Phi_\mathcal{A},W_\mathcal{A})$ and $({\rm Ind}(\Phi_\mathcal{A}),{\rm Ind}(W))$ (resp. $({\rm Pro}(\Phi_\mathcal{A}),{\rm Pro}(W))$), it induces a natural fully faithful functor
$$\hspace{5mm} \mathfrak{F}({\iota}_{\rm Ind}):\mathfrak{F}(\mathcal{A},\Phi_\mathcal{A},W_\mathcal{A})\rightarrow\mathfrak{F}{\rm Ind}(\mathcal{A},\Phi_\mathcal{A},W_\mathcal{A})$$
$$({\rm resp.}\hspace{3mm} \mathfrak{F}({\iota}_{\rm Pro}):\mathfrak{F}(\mathcal{A},\Phi_\mathcal{A},W_\mathcal{A})\rightarrow\mathfrak{F}{\rm Pro}(\mathcal{A},\Phi_\mathcal{A},W_\mathcal{A})\hspace{1mm}).$$

Let $(F,\sigma):(\mathcal{A},\Phi_\mathcal{A},W_\mathcal{A})\rightarrow(\mathcal{B},\Phi_\mathcal{B},W_\mathcal{B})$ be a cwp-functor. Then 
$${\rm Ind}(F,\sigma):=({\rm Ind}(F),{\rm Ind}(\sigma)):{\rm Ind}(\mathcal{A},\Phi_\mathcal{A},W_\mathcal{A})\rightarrow{\rm Ind}(\mathcal{B},\Phi_\mathcal{B},W_\mathcal{B})$$
$${\rm Pro}(F,\sigma):=({\rm Pro}(F),{\rm Pro}(\sigma)):{\rm Pro}(\mathcal{A},\Phi_\mathcal{A},W_\mathcal{A})\rightarrow{\rm Pro}(\mathcal{B},\Phi_\mathcal{B},W_\mathcal{B})$$ 
are  cwp-functors, and the following diagrams are commutative:
\[\xymatrix{
\mathfrak{F}{\rm Ind}(\mathcal{A},\Phi_\mathcal{A},W_\mathcal{A})\ar[rr]^{\mathfrak{F}{\rm Ind}(F,\sigma)}&&\mathfrak{F}{\rm Ind}(\mathcal{B},\Phi_\mathcal{B},W_\mathcal{B})&\mathfrak{F}{\rm Pro}(\mathcal{A},\Phi_\mathcal{A},W_\mathcal{A})\ar[rr]^{\mathfrak{F}{\rm Pro}(F,\sigma)}&&\mathfrak{F}{\rm Pro}(\mathcal{B},\Phi_\mathcal{B},W_\mathcal{B})\\
\mathfrak{F}(\mathcal{A},\Phi_\mathcal{A},W_\mathcal{A})\ar[rr]^{\mathfrak{F}(F,\sigma)} \ar[u]^{\mathfrak{F}(\iota_{\rm Ind})}&&\mathfrak{F}(\mathcal{B},\Phi_\mathcal{B},W_\mathcal{B})\ar[u]_{\mathfrak{F}(\iota_{\rm Ind})}&\mathfrak{F}(\mathcal{A},\Phi_\mathcal{A},W_\mathcal{A})\ar[rr]^{\mathfrak{F}(F,\sigma)} \ar[u]^{\mathfrak{F}(\iota_{\rm Pro})}&&\mathfrak{F}(\mathcal{B},\Phi_\mathcal{B},W_\mathcal{B})\ar[u]_{\mathfrak{F}(\iota_{\rm Pro})}
  }\]


\vspace{7mm}
\section{Derived  factorization categories of  gauged LG models}
\vspace*{2mm}
Let $X$ be a   quasi-projective variety and let $G$ be an affine  algebraic group acting on $X$ over an algebraically closed field $k$ of characteristic zero. Let $\sigma:G\times X\rightarrow X$ be the morphism defining the action, $\pi:G\times X\rightarrow X$ be a projection and $\iota:X\rightarrow G\times X$ be an embedding given by $x\mapsto(e, x)$, where $e\in G$ is the identity of group $G$.

\subsection{Equivariant sheaves and factorization categories of gauged LG models}

\vspace*{2mm}

\begin{dfn}
A quasi-coherent (resp. coherent) \textbf{$G$-equivariant} sheaf is a pair $(\mathcal{F}, \theta)$ of a quasi-coherent (resp. coherent) sheaf $\mathcal{F}$ and an isomorphism $\theta:\pi^*\mathcal{F}\xrightarrow{\sim}\sigma^*\mathcal{F}$ such that 
$$\iota^*\theta={\rm id}_{\mathcal{F}} \hspace{4mm}\rm{and} \hspace{5mm}\big{(}(1_G\times\sigma)\circ({\it s}\times1_X)\big{)}^*\theta\circ(1_G\times\pi)^*\theta=({\it m}\times1_X)^*\theta,$$
where $m:G\times G\rightarrow G$ is the multiplication and $s:G\times G\rightarrow G\times G$ is the switch of two factors.
A \textbf{$G$-invariant} morphism $\varphi:(\mathcal{F}_1, \theta_1)\rightarrow(\mathcal{F}_2, \theta_2)$ of equivariant sheaves is a morphism of sheaves $\varphi:\mathcal{F}_1\rightarrow\mathcal{F}_2$ which is commutative with $\theta_i$, i.e. $\sigma^*\varphi\circ\theta_1=\theta_2\circ\pi^*\varphi$.

We denote by ${\rm Qcoh}_G(X)$ (resp. ${\rm coh}_G(X)$) the category of quasi-coherent (resp. coherent) $G$-equivariant sheaves on $X$ whose morphisms are $G$-invariant morphisms. And  we denote by ${\rm Inj}_G(X)$, ${\rm LFr}_G(X)$ and ${\rm lfr}_G(X)$ the full subcategories of ${\rm Qcoh}_G(X)$ consisting of injective quasi-coherent equivariant  sheaves, locally free equivariant sheaves and locally free equivariant sheaves of finite ranks.
\end{dfn}

\vspace*{2mm} Let $L\in{\rm Pic}_G(X)$ be a $G$-equivariant invertible sheaf on $X$ and let $W\in H^0(X,L)^G$ be an invariant section of $L$. 

\begin{dfn}
We call the data $(X,L,W)^G$ a \textbf{gauged Landau-Ginzburg model}  or \textbf{gauged LG model}, for short. We sometimes drop the script $L$ from the notion $(X,L,W)^G$, and write $(X,W)^G$ if there is  no confusion.
\end{dfn}

The pair $(L,W):=((-)\otimes L, (-)\otimes W)$ is a potential of ${\rm Qcoh}_G(X)$, ${\rm coh}_G(X)$, ${\rm Inj}_G(X)$, ${\rm LFr}_G(X)$ and ${\rm lfr}_G(X)$, where $W$ is considered as the morphism $W:\mathcal{O}_X\rightarrow L$ corresponding to the section of $L$.

\begin{dfn}
We define factorization categories of $(X, L, W)^G$ as 
$${\rm Qcoh}_G(X, L, W):=\mathfrak{F}({\rm Qcoh}_G(X), L,  W)$$
$${\rm coh}_G(X, L, W):=\mathfrak{F}({\rm coh}_G(X),  L,  W)$$
$${\rm Inj}_G(X, L, W):=\mathfrak{F}({\rm Inj}_G(X),  L,  W)$$
$${\rm LFr}_G(X, L, W):=\mathfrak{F}({\rm LFr}_G(X),  L,  W)$$
$${\rm lfr}_G(X, L, W):=\mathfrak{F}({\rm lfr}_G(X),  L,  W).$$
We define categories of acyclic factorizations as
$${\rm Acycl}_G(X,L,W):={\rm Acycl}^{\rm abs}({\rm Qcoh}_G(X), L,  W)$$
$${\rm Acycl}^{\rm co}_G(X,L,W):={\rm Acycl}^{\rm co}({\rm Qcoh}_G(X), L,  W)$$
  and  derived  factorization categories are defined as
$${\rm DQcoh}_G(X, L, W):={\rm D^{abs}}({\rm Qcoh}_G(X), L,  W)$$
$${\rm Dcoh}_G(X, L, W):={\rm D^{abs}}({\rm coh}_G(X),  L,  W)$$
$${\rm DLFr}_G(X, L, W):={\rm D^{abs}}({\rm LFr}_G(X), L,  W)$$
$${\rm Dlfr}_G(X, L, W):={\rm D^{abs}}({\rm lfr}_G(X),  L,  W).$$
We call the category ${\rm Dcoh}_G(X, L, W)$ \textbf{the derived factorization category} of a gauged LG model $(X,L,W)^G$.
For $E, F\in {\rm Qcoh}_G(X, L, W)$, we say $E$ and $F$ are \textbf{quasi-isomorphic} if $E$ and $F$ are isomorphic in ${\rm DQcoh}_G(X, L, W)$. We denote by ${\rm D_{coh}Qcoh}_G(X, L, W)$ the full subcategory of ${\rm DQcoh}_G(X, L, W)$ whose objects are quasi-isomorphic to objects in ${\rm Dcoh}_G(X, L, W)$.
If $G$ is trivial, we drop the subscript $G$ in the above notations.
\end{dfn}

\vspace{2mm}
\begin{rem}\label{co=abs}
 By  Lemma \ref{cocomplete}, if $X$ is smooth,  then ${\rm Acycl}_G(X,L,W)={\rm Acycl}^{\rm co}_G(X,L,W)$ and hence ${\rm DQcoh}_G(X, L, W)={\rm D^{co}}({\rm Qcoh}_G(X),  L,  W)$.
\end{rem}

\vspace{2mm}
\begin{dfn}
A gauged LG model  $(X,\mathcal{O}(\chi),0)^{G\times\mathbb{G}_m}$ such that  the potential is zero, the character $\chi:G\times\mathbb{G}_m\rightarrow\mathbb{G}_m$ is  projection, and the action of $\mathbb{G}_m$ is trivial, is called \textbf{of $\sigma$-type}. If $G$ is trivial, the gauged LG model $(X,\mathcal{O}(\chi),0)^{\mathbb{G}_m}$ of $\sigma$-type is called \textbf{of trivial $\sigma$-type}.
\end{dfn}

\vspace{2mm}
The derived factorization category of a gauged LG models of $\sigma$-type is equivalent to bounded derived category of coherent sheaves on some algebraic stack.

\begin{prop}[\cite{vgit}, Corollary 2.3.12]\label{sigma LG}
Let $(X,\mathcal{O}(\chi),0)^{G\times\mathbb{G}_m}$ be a gauged LG model of $\sigma$-type. Then we have an equivalence
$${\rm Dcoh}_{G\times\mathbb{G}_m}(X,\mathcal{O}(\chi),0)\cong{\rm D^b}({\rm coh}[X/G]).$$

\end{prop}

\vspace{2mm}
The following lemma is necessary to replace objects of ${\rm DQcoh}_G(X, L, W)$ or ${\rm Dcoh}_G(X, L, W)$ to ones with injective (or locally free) components. These replacements ensure that we can define derived functors between  derived factorization categories from exact functors between  homotopy categories of  factorization categories. 

\begin{lem}[cf. \cite{ls}, Lemma 2.10.]\label{resolutions}
Assume that $X$ is smooth. Then we have
\begin{itemize}
\item[$(1)$]For any $F\in {\rm Qcoh}_G(X,L,W)$ there exists a bounded exact sequence $0\rightarrow F\rightarrow I^0\rightarrow\cdot\cdot\cdot\rightarrow I^n\rightarrow0$ in $Z^0({\rm Qcoh}_G(X,L,W))$ with all $I^m \in {\rm Inj}_G(X,L,W)$. In particular, there is an isomorphism $F\rightarrow {\rm Tot}(I^{\text{\tiny{\textbullet}}})$ in ${\rm DQcoh}_G(X,L,W)$.
\vspace{2mm}
\item[$(2)$]For any object $F$ of ${\rm Qcoh}_G(X,L,W)$ (resp. ${\rm coh}_G(X,L,W)$) there exists a bounded exact sequence $0\rightarrow P^n\rightarrow\cdot\cdot\cdot\rightarrow P^0\rightarrow F\rightarrow0$ in $Z^0({\rm Qcoh}_G(X,L,W))$ (resp. $Z^0({\rm coh}_G(X,L,W))$) with all $P^m$ in ${\rm LFr}_G(X,L,W)$ (resp. ${\rm lfr}_G(X,L,W)$). In particular, we have an isomorphism ${\rm Tot}(P^{\text{\tiny{\textbullet}}})\rightarrow F$.
\end{itemize}
\begin{proof}
This is an equivariant version of \cite[Lemma 2.10]{ls}. Since ${\rm Qcoh}_G(X,L,W)$ has enough injective objects and for any equivariant sheaf $E\in{\rm Qcoh}_G(X)$ there exist an equivariant locally free sheaf $P$ and surjection $P\rightarrow E$ (see e.g. \cite[Proposition 5.1.26]{cg}), the exact sequences can be constructed in a similar way as in \cite[Lemma 2.10]{ls}. 
\end{proof}
\end{lem}

\vspace{2mm}
\begin{lem}[\cite{bfk} Proposition 3.11]\label{inj}
Assume $X$ is smooth. We have
$${\rm Hom}_{H^0({\rm Qcoh}_G(X,L,W))}(A,I)=0$$
for any $A\in {\rm Acycl}_G(X,L,W)$ and $I\in H^0({\rm Inj}_G(X,L,W))$. Moreover, the following compositions are equivalences;
$$H^0({\rm Inj}_G(X,L,W))\rightarrow H^0({\rm Qcoh}_G(X,L,W))\rightarrow {\rm DQcoh}_G(X,L,W)$$
$$H^0({\rm inj}_G(X,L,W))\rightarrow H^0({\rm Qcoh}_G(X,L,W))\rightarrow {\rm Dcoh}_G(X,L,W),$$
where ${\rm inj}_G(X,L,W)$ is the dg-subcategory of ${\rm Inj}_G(X,L,W)$ consisting of factorizations  which are quasi-isomorphic to  factorizations with coherent components.
\end{lem}

\vspace{2mm}
Since the embedding $H^0({\rm inj}_G(X,L,W))\rightarrow H^0({\rm Inj}_G(X,L,W))$ is fully faithful, so is 
$${\rm Dcoh}_G(X,L,W)\rightarrow{\rm DQcoh}_G(X,L,W)$$
by the above lemma. Hence we have a natural equivalence,
$${\rm Dcoh}_G(X,L,W)\xrightarrow{\sim}{\rm D_{coh}Qcoh}_G(X,L,W)$$

\vspace{2mm}
\begin{lem}[\cite{bfk} Proposition 3.14]\label{lfr}
Assume $X$ is smooth. The following natural functors are equivalences:
$${\rm DLFr}_G(X, L, W)\rightarrow{\rm DQcoh}_G(X, L, W)$$
$${\rm Dlfr}_G(X, L, W)\rightarrow{\rm Dcoh}_G(X, L, W)$$
\end{lem}

\vspace{2mm}
\begin{lem}[cf. \cite{ls}, Corollary 2.23.]\label{cocomplete}
Assume $X$ is smooth. The categories $H^0({\rm Qcoh}_G(X,L,W))$, $H^0({\rm Inj}_G(X,L,W))$, ${\rm Acycl}_G(X,L,W)$ and 
${\rm DQcoh}_G(X,L,W)$ are closed under arbitrary direct sums and therefore idempotent complete.
\begin{proof}
We can prove this in a  similar way as in \cite[Corollary 2.23]{ls}.
\end{proof}

\end{lem}

\vspace{2mm}
We define the supports of factorizations and complexes of factorizations as follows: 

\vspace{2mm}
\begin{dfn}\label{supp}
Let $E\in Z^0({\rm coh}_G(X,L,W))$.  \textbf{The support} ${\rm Supp}(E)$ of $E$ is defined as 
$$
{\rm Supp}(E):={\rm Supp}(E_1)\cup{\rm Supp}(E_0).
$$ 
For an object $E^{\text{\tiny{\textbullet}}}\in {\rm D^{b}}(Z^0({\rm coh}_G(X,L,W)))$, we define \textbf{the support} ${\rm Supp}(E^{\text{\tiny{\textbullet}}})$ of $E^{\text{\tiny{\textbullet}}}$ as
$${\rm Supp}(E^{\text{\tiny{\textbullet}}}):=\bigcup\limits_{i\in\mathbb{Z}}{\rm Supp}(H^i(E^{\text{\tiny{\textbullet}}})).$$
\end{dfn}

\vspace{2mm}

\begin{rem}\label{support}
By definition  the support of $E^{\text{\tiny{\textbullet}}}\in{\rm D^{b}}(Z^0({\rm coh}_G(X,L,W)))$ is the union of supports of objects $E^{\text{\tiny{\textbullet}}}_i\in{\rm D^b}(X)$, i.e. ${\rm Supp}(E^{\text{\tiny{\textbullet}}})=\bigcup_{i=0,1}{\rm Supp}(E^{\text{\tiny{\textbullet}}}_i)$, where the support of a complex in ${\rm D^b}(X)$ is defined by the union of supports of its cohomologies.
\end{rem}

\vspace{2mm}

In the following, we define properness of the $^{\rotatebox[origin=C]{180}{"}}$support" of an object in ${\rm Dcoh}_G(X,L,W)$ by using totalization.

\vspace{2mm}
\begin{dfn}

Let $f:X\rightarrow Y$ be a morphism, where $Y$ is another quasi-projective variety. A closed subset $Z$ of $X$ is called \textbf{$f$-proper } if the composition $Z\hookrightarrow X\xrightarrow{f}Y$ is a proper morphism. We denote by ${\rm coh}_{\sqcap G}^f(X,L,W)$ the full subcategory of ${\rm coh}_G(X,L,W)$ consisting of objects whose supports are $f$-proper.

Let $F$ be  an object in ${\rm Dcoh}_G(X,L,W)$. We say $F$ \textbf{has a $f$-proper support} if there exists an object 
$F^{\text{\tiny{\textbullet}}}\in {\rm D^{b}}(Z^0({\rm coh}_G(X,L,W)))$ such that ${\rm Tot}(F^{\text{\tiny{\textbullet}}})$ is isomorphic to $F$ in ${\rm Dcoh}_G(X,L,W)$ and the closed subset
${\rm Supp}(F^{\text{\tiny{\textbullet}}})$ is $f$-proper.

We denote by ${\rm D}^f_{\sqcap}{\rm coh}_G(X,L,W)$ the full subcategory of ${\rm Dcoh}_G(X,L,W)$ consisting of objects which have $f$-proper supports.

\end{dfn}

\vspace{2mm}
\begin{rem}

\hfill\\
$(1)$ ${\rm D}^f_{\sqcap}{\rm coh}_G(X,L,W)$ is strictly full subcategory, i.e. closed under isomorphisms in ${\rm Dcoh}_G(X,L,W)$.\\
$(2)$ If $f$ is proper morphism then  ${\rm D}^f_{\sqcap}{\rm coh}_G(X,L,W)={\rm Dcoh}_G(X,L,W)$.\\
$(3)$ Let $g:Y\rightarrow Z$ be another morphism of quasi-projective varieties. If $F\in{\rm Dcoh}_G(X,L,W)$ has a $g\circ f$-proper support, then $F$ has a $f$-proper support.\\
$(4)$ An object $E\in{\rm Dcoh}_G(X,L,W)$ which is quasi-isomorphic to $F\in {\rm coh}_{\sqcap G}^f(X,L,W)$ has a $f$-proper support. 

\end{rem}


\vspace{3mm}
\subsection{Functors of factorization categories of gauged LG models}
Throughout this section, we assume $X$ is smooth. In what follows, we define exact functors between derived equivariant factorization categories. 

\subsubsection{Derived functors between triangulated categories} In this section, we recall definitions and generalities on derived functors of exact functors of triangulated categories after  \cite{murfet}. Let $\mathcal{D}$ be a triangulated category, and let $\mathcal{C}$ be a full subcategory of $\mathcal{D}$ with Verdier quotient $Q:\mathcal{D}\rightarrow\mathcal{D}/\mathcal{C}$. Throughout this section, all functor morphisms of exact functors are assumed to be commutative with shift functors, i.e. if $\alpha: F\rightarrow G$
is a functor morphism between exact functors $F,G:\mathcal{T}\rightarrow\mathcal{T}'$ of triangulated categories $\mathcal{T}$ and $\mathcal{T}'$ with shift functors $\Sigma:\mathcal{T}\rightarrow\mathcal{T}$ and 
$\Sigma':\mathcal{T}'\rightarrow\mathcal{T}'$, then $\alpha$ satisfies the commutativity of the following diagram of functor morphisms,
 \[\xymatrix{
 F\Sigma\ar[r]^{\sim} \ar[d]_{F\alpha}& \Sigma'F \ar[d]^{\Sigma'\alpha}\\
 G\Sigma\ar[r]^{\sim}& \Sigma'G. 
  }\]

\vspace{2mm}
\begin{dfn}
Let $F:\mathcal{D}\rightarrow \mathcal{T}$ be an exact functor of triangulated categories.
The \textbf{right derived functor of} $F$ (with respect to $\mathcal{C}$) is a pair $({\bf R}F,\zeta)$ of an exact functor ${\bf R}F:\mathcal{D}/\mathcal{C}\rightarrow\mathcal{T}$ and functor morphism $\zeta:F\rightarrow{\bf R}F\circ Q$ with the following universal property: for any exact functor $G:\mathcal{D}/\mathcal{C}\rightarrow\mathcal{T}$ and functor morphism $\rho:F\rightarrow G\circ Q$ there is a unique functor morphism $\eta:{\bf R}F\rightarrow G$ making the following diagram  commute:
 \[\xymatrix{
&F \ar[ld]_{\zeta} \ar[rd]^{\rho}&\\
{\bf R}F\circ Q\ar[rr]^{\eta Q} && G\circ Q.
 }\]
We will often drop the subcategory $\mathcal{C}$ and $\zeta$ from the notation, and say simply that ${\bf R}F$ is right derived functor of $F$.
\end{dfn}

\vspace{2mm}

\begin{rem}
By the definition, if right derived functor exists, it is unique up to natural equivalence.

\end{rem}

\vspace{2mm}
\begin{dfn}
Let $F:\mathcal{D}\rightarrow \mathcal{T}$ be an exact functor. An object $A\in\mathcal{D}$ is \textbf{right} $F$\textbf{-acyclic} with respect to $\mathcal{C}$ if the following condition holds: if $s:A\rightarrow B$ is a morphism with cone in $\mathcal{C}$, there is a morphism $t:B\rightarrow C$ with cone in $\mathcal{C}$ such that $F(ts)$ is an isomorphism. 
\end{dfn}

\begin{rem}
If $A\in\mathcal{D}$ is a right $F$-acyclic with respect to $\mathcal{C}$ and in $\mathcal{C}$, then $F(A)=0$.
\end{rem}

\vspace{2mm}

The following theorem will be applied several times in the following  sections to construct exact functors between derived  factorization categories.

\begin{thm}[\cite{murfet} Theorem 116]\label{derived functor}
Let $F:\mathcal{D}\rightarrow \mathcal{T}$ be an exact functor. Assume $\mathcal{C}$ is a thick subcategory of $\mathcal{D}$. Suppose that for each object $X\in\mathcal{D}$ there exists a right $F$-acyclic object $A_X$ and a morphism $\eta_X:X\rightarrow A_X$ with cone in $\mathcal{C}$. Then $F$ admits a right derived functor $({\bf R}F,\zeta)$ with the following properties
\begin{itemize}
\item[$(1)$] For any object $X\in\mathcal{D}$ we have ${\bf R}F(X)=F(A_X)$ and $\zeta(X)=F(\eta_X)$.
\item[$(2)$] An object $X\in\mathcal{D}$ is right $F$-acyclic if and only if $\zeta(X)$ is an isomorphism in $\mathcal{T}$.
\end{itemize}
\end{thm}

\vspace{2mm}
There are similar definitions and results for left derived functors. See \cite{murfet} for the detail.

\vspace{3mm}
\subsubsection{Direct and inverse image}

Let  $Y$ be another smooth quasi-projective variety with an action of $G$, defined by $\tau: G\times Y\rightarrow Y$, and let $f:X\rightarrow Y$ be an equivariant morphism, i.e. $f\circ\sigma=\tau\circ(1_{G}\times f)$.

For the morphism $f$, the direct image $f_*:{\rm Qcoh}_G(X)\rightarrow{\rm Qcoh}_G(Y)$   and the inverse image $f^*:{\rm Qcoh}_G(Y)\rightarrow{\rm Qcoh}_G(X)$ are defined by 
$$f_*(\mathcal{F}, \theta):=(f_*(\mathcal{F}), (1\times f)_*\theta)\hspace{3mm}{\rm and}\hspace{3mm}f^*(\mathcal{F}, \theta):=(f^*\mathcal{F}, (1\times f)^*\theta).$$

Let $L\in {\rm Pic}_G(Y)$ be an equivariant invertible sheaf on $Y$ and let $W\in H^0(Y, L)^G$ be an invariant section of $L$. Then we have  potentials $(f^*L,f^*W)$ and $(L ,W)$ of ${\rm Qcoh}_G(X)$ and  ${\rm Qcoh}_G(Y)$ respectively.  By the natural isomorphisms of functors $f_*((-)\otimes f^*L)\cong f_*(-)\otimes L$ and $f^*((-)\otimes L)\cong f^*(-)\otimes f^*L$,  we see that the direct image $f_*$ and inverse image $f^*$ are compatible with potentials with respect to  $(f^*L,f^*W)$ and $(L ,W)$ (see Definition \ref{cwp-functor}). So we have direct image $f_*$ and inverse image $f^*$, denoted by the same notation as usual ones, between  factorization categories 
$$f_*: {\rm Qcoh}_G(X, f^*L, f^*W)\rightarrow {\rm Qcoh}_G(Y, L, W)$$
$$f^*:{\rm Qcoh}_G(Y, L, W)\rightarrow{\rm Qcoh}_G(X, f^*L, f^*W).$$
Taking  $H^0(-)$ of these dg-functors, we have exact functors 
$$f_*: H^0({\rm Qcoh}_G(X, f^*L, f^*W))\rightarrow H^0({\rm Qcoh}_G(Y, L, W))$$
$$f^*:H^0({\rm Qcoh}_G(Y, L, W))\rightarrow H^0({\rm Qcoh}_G(X, f^*L, f^*W)).$$
Since these exact functors don't send acyclic objects to acyclic ones in general, we need to take derived functors of them. In the following, we give a proposition that implies existences of derived functors and two lemmas about them, following  \cite{ls}. Since the proofs are same as \cite{ls}, we will omit proofs.

Denote the following compositions by same notation $f_*$ and $f^*$,
$$f_*: H^0({\rm Qcoh}_G(X, f^*L, f^*W))\rightarrow H^0({\rm Qcoh}_G(Y, L, W))\rightarrow{\rm DQcoh}_G(Y, L, W)$$
$$f^*:H^0({\rm Qcoh}_G(Y, L, W))\rightarrow H^0({\rm Qcoh}_G(X, f^*L, f^*W))\rightarrow{\rm DQcoh}_G(X, f^*L, f^*W).$$
By Lemma \ref{resolutions} and Theorem \ref{derived functor}, we have the following:

\vspace{2mm}
\begin{prop}[cf. \cite{ls} Theorem 2.35] 

\hfill

\begin{itemize} 

\item[$(1)$] The functor $f_*: H^0({\rm Qcoh}_G(X, f^*L, f^*W))\rightarrow{\rm DQcoh}_G(Y, L, W)$ admits a right derived functor ${\bf R}f_*:{\rm DQcoh}_G(X, f^*L, f^*W)\rightarrow{\rm DQcoh}_G(Y, L, W)$ with respect to ${\rm Acycl}_G(X, f^*L, f^*W)$.

\item[$(2)$] The functor  $f^*:H^0({\rm Qcoh}_G(Y, L, W))\rightarrow{\rm DQcoh}_G(X, f^*L, f^*W)$ has a left derived functor ${\bf L}f^*{\rm DQcoh}_G(Y, L, W)\rightarrow{\rm DQcoh}_G(X, f^*L, f^*W)$ with respect to ${\rm Acycl}_G(Y, L, W)$. This left derived functor maps to ${\rm Dcoh}_G(Y, L, W)$ to ${\rm Dcoh}_G(X, f^*L, f^*W)$. 
\end{itemize}
\end{prop}

\vspace{2mm}
The right derived functor ${\bf R}f_*$ doesn't map an object $E\in{\rm Dcoh}_G(X, f^*L, f^*W)$ to an object in ${\rm Dcoh}_G(Y, L, W)$ in general. But the following Lemma \ref{cohomology} implies  that  if $E$ has a $f$-proper support, then  ${\bf R}f_*(E)$ is isomorphic to an object in ${\rm Dcoh}_G(Y, L, W)$. In particular, if $f$ is proper morphism then ${\bf R}f_*$ maps an object in ${\rm Dcoh}_G(X, f^*L, f^*W)$ to an object which is isomorphic to an object in ${\rm Dcoh}_G(Y, L, W)$ and we also denote by ${\bf R}f_*$ the following composition  
$${\rm Dcoh}_G(X, f^*L, f^*W)\xrightarrow{{\bf R}f_*}{\rm D_{coh}Qcoh}_G(Y,L,W)\xrightarrow{\sim}{\rm Dcoh}_G(Y,L,W).$$

\vspace{2mm}
\begin{lem}[\cite{ls} Lemma 2.40]\label{cohomology}
Let $F\in {\rm Ch^{b}}(Z^0({\rm Qcoh}_G(Y,L,W))$. If each $H^i(F)\in{\rm DQcoh}_G(Y, L, W)$ is isomorphic to an object in ${\rm Dcoh}_G(Y, L, W)$, then so is ${\rm Tot}(F)$. 
\end{lem}

\vspace{2mm}

\vspace{2mm}
\begin{lem}[\cite{ls} Lemma 2.38]\label{f-acyclic}
Let $E=(E_1\rightarrow E_0\rightarrow E_1\otimes f^*L)\in H^0({\rm Qcoh}_G(X, f^*L, f^*W))$ and assume 
that ${\bf R}^{i}f_*(E_n)=0$ in ${\rm Qcoh}_G(Y)$ for any $i>0$ and each $n=0,1$. Then $E$ is right $f_*$-acyclic. In particular, if $f$ is affine morphism then we have a canonical isomorphism of functors $f_*\xrightarrow{\sim} {\bf R}f_*$.

Similarly, if $F=(F_1\rightarrow F_0\rightarrow F_1\otimes L)\in H^0({\rm Qcoh}_G(Y, L, W))$ and ${\bf L}^{j}f^*(F_m)=0$ in ${\rm Qcoh}_G(X)$ for any $j>0$ and each $m=0,1$, then $F$ is left $f^*$-acyclic.
In particular, if $f$ is flat morphism then ${\bf L}f^*\xrightarrow{\sim} f^*$. 

\end{lem}

\vspace{2mm} 
Since the direct image  $f_*:{\rm Qcoh}_G(X)\rightarrow{\rm Qcoh}_G(Y)$  is right cwp-adjoint to the inverse image $f^*:{\rm Qcoh}_G(Y)\rightarrow{\rm Qcoh}_G(X)$ with respect to potentials $(f^*L,f^*W)$ and $(L ,W)$,  $f_*: {\rm Qcoh}_G(X, f^*L, f^*W)\rightarrow {\rm Qcoh}_G(Y, L, W)$ is right adjoint to $f^*:{\rm Qcoh}_G(Y, L, W)\rightarrow{\rm Qcoh}_G(X, f^*L, f^*W)$, whose adjunction morphisms are of degree zero. Taking $H^0(-)$, we see that $f_*: H^0({\rm Qcoh}_G(X, f^*L, f^*W))\rightarrow H^0({\rm Qcoh}_G(Y, L, W))$ is right adjoint to $f^*:H^0({\rm Qcoh}_G(Y, L, W))\rightarrow H^0({\rm Qcoh}_G(X, f^*L, f^*W))$. Thus, by \cite[ Theorem 122]{murfet}, we obtain the following adjoint pair:
$${\bf L}f^*\dashv{\bf R}f_*$$

\vspace{2mm}
\subsubsection{Tensor product and local Hom} Let $L\in {\rm Pic }_G(X)$ and $V,W\in H^0(X,L)^G$.

Taking tensor product gives a bifunctor $(-)\otimes(-):{\rm Qcoh}_G(X)\times{\rm Qcoh}_G(X)\rightarrow{\rm Qcoh}_G(X)$. Note that this functor is compatible with potentials with respect to potentials $(L,V)$, $(L,W)$ and $(L,V+W)$ (see Definition \ref{cwp-cobifunctor}). So it induces a dg-bifuctor 
$$(-)\otimes(-):{\rm Qcoh}_G(X,L,V)\otimes{\rm Qcoh}_G(X,L,W)\rightarrow{\rm Qcoh}_G(X,L,V+W).$$
If we  fix an object  $P\in {\rm Qcoh}_G(X,L,W)$, we have an exact functor 
$$(-)\otimes P:H^0({\rm Qcoh}_G(X,L,V))\rightarrow {\rm DQcoh}_G(X,L,V+W).$$

\vspace{2mm}
\begin{prop}
The functor $(-)\otimes P:H^0({\rm Qcoh}_G(X,L,V))\rightarrow {\rm DQcoh}_G(X,L,V+W)$ has a left derived functor $(-)\otimes^{\bf L} P:{\rm DQcoh}_G(X,L,V)\rightarrow {\rm DQcoh}_G(X,L,V+W)$ with respect to ${\rm Acycl}_G(X,L,V)$. If $P\in{\rm coh}_G(X,L,W)$ then this left derived functor maps ${\rm Dcoh}_G(X,L,V)$ to ${\rm Dcoh}_G(X,L,V+W)$.

\begin{proof}
The proof is very similar to the proof of \cite[Theorem 2.35 (b)]{ls}, and the detail  is left to the reader.
\end{proof}
\end{prop}

\vspace{2mm}
\begin{dfn}\label{complextensor}
For any complex $C^{\text{\tiny{\textbullet}}}\in{\rm D^b}({\rm Qcoh}_G(X))$, we define an exact functor
$$(-)\otimes^{\bf L}C^{\text{\tiny{\textbullet}}} :{\rm DQcoh}_G(X,L,W)\rightarrow {\rm DQcoh}_G(X,L,W)$$
as $$E\otimes^{\bf L}C^{\text{\tiny{\textbullet}}}:=E\otimes^{\bf L}\Upsilon(C^{\text{\tiny{\textbullet}}}),$$
where $\Upsilon:{\rm D^b}({\rm Qcoh}_G(X))\rightarrow{\rm DQcoh}_G(X,L,0)$ is the  functor defined in Definition \ref{upsilon}.
We denote by $E\otimes C^{\text{\tiny{\textbullet}}}$ if $E\otimes^{\bf L}\Upsilon(C^{\text{\tiny{\textbullet}}})\cong E\otimes\Upsilon(C^{\text{\tiny{\textbullet}}})$.
\end{dfn}

\vspace{2mm}
Taking local Hom gives a bifunctor $\mathcal{H}om(-,-):{\rm coh}_G(X)^{\rm op}\times{\rm Qcoh}_G(X)\rightarrow{\rm Qcoh}_G(X)$. Note that this bifunctor is compatible with potentials with respect to potentials $(L,V)$, $(L,W)$ and $(L,W-V)$ (see Definition \ref{cwp-contrabifunctor}). So it induces a dg-bifunctor 
$$\mathcal{H}om(-,-):{\rm coh}_G(X,L,V)\otimes{\rm Qcoh}_G(X,L,W)\rightarrow{\rm Qcoh}_G(X,L,W-V).$$
If we fix an object $Q\in {\rm coh}_G(X,L,V)$, we have an exact functor
$$\mathcal{H}om(Q,-):H^0({\rm Qcoh}_G(X,L,W))\rightarrow {\rm DQcoh}_G(X,L,W-V)$$

\vspace{2mm}
\begin{prop}
The functor $\mathcal{H}om(Q,-):H^0({\rm Qcoh}_G(X,L,W))\rightarrow {\rm DQcoh}_G(X,L,W-V)$ has a right derived functor ${\bf R}\mathcal{H}om(Q,-):{\rm DQcoh}_G(X,L,W)\rightarrow {\rm DQcoh}_G(X,L,W-V)$ with respect to ${\rm Acycl}_G(X,L,W)$.
\begin{proof}
The proof is very similar to the proof of \cite[Theorem 2.35 (a)]{ls}, and the detail  is left to the reader.
\end{proof}
\end{prop}

\vspace{2mm}
By Lemma \ref{cohomology}, if $E\in{\rm Dcoh}_G(X,L,W)$, then ${\bf R}\mathcal{H}om(Q,E)\in  {\rm D_{coh}Qcoh}_G(X,L,W-V)$. We use same notation ${\bf R}\mathcal{H}om(Q,-)$ for the composition
$${\rm Dcoh}_G(X,L,W)\xrightarrow{{\bf R}\mathcal{H}om(Q,-)} {\rm D_{coh}Qcoh}_G(X,L,W-V)\xrightarrow{\sim}{\rm Dcoh}_G(X,L,W-V).$$

\vspace{2mm}
\begin{lem} 
Let $E=(E_1\rightarrow E_0\rightarrow E_1\otimes L)\in H^0({\rm Qcoh}_G(X,L,V))$ and $P=(P_1\rightarrow P_0\rightarrow P_1\otimes L)\in {\rm Qcoh}_G(X,L,W)$. If $\mathcal{T}or^i(E_n,P_m)=0$ for any $i>0$ and any $n,m\in\{0,1\}$, then $E$ is $(-)\otimes P$-acyclic object. In particular, if $P\in {\rm LFr}_G(X,L,W)$, then there is an isomorphism of exact functors $(-)\otimes^{\bf L} P\xrightarrow{\sim}(-)\otimes P$.

Let $F=(F_1\rightarrow F_0\rightarrow F_1\otimes L)\in H^0({\rm Qcoh}_G(X,L,W))$ and $Q=(Q_1\rightarrow Q_0\rightarrow Q_1\otimes L)\in {\rm coh}_G(X,L,V)$. If $\mathcal{E}xt^i(Q_n,F_m)=0$ for each $i>0$ and any $n,m\in\{0,1\}$, then $F$ is $\mathcal{H}om(Q,-)$-acyclic object. In particular, if $Q\in {\rm lfr}_G(X,L,V)$, there is an isomorphism of exact functors $\mathcal{H}om(Q,-)\xrightarrow{\sim}{\bf R}\mathcal{H}om(Q,-)$.

\begin{proof}
The proof is similar to \cite[Lemma 2.38]{ls}, and we leave the detail to the reader.
\end{proof}
\end{lem}

 \vspace{2mm}
\begin{rem}
In the above lemma, we can take $P$ and $Q$ as objects whose components are flat sheaves.
\end{rem}

\vspace{2mm}
\begin{prop}[\cite{bfk} Proposition 3.27]
Let $R\in {\rm coh}_G(X,L,V)$. Then $\mathcal{H}om(R,-):{\rm Qcoh}_G(X,L,W)\rightarrow {\rm Qcoh}_G(X,L,W-V)$ is right adjoint to $(-)\otimes R:{\rm Qcoh}_G(X,L,W-V)\rightarrow {\rm Qcoh}_G(X,L,W)$.
\end{prop}

\vspace{2mm}
$\mathcal{H}om(R,-):H^0({\rm Qcoh}_G(X,L,W))\rightarrow H^0({\rm Qcoh}_G(X,L,W-V))$ is right adjoint to $(-)\otimes R:H^0({\rm Qcoh}_G(X,L,W-V))\rightarrow H^0({\rm Qcoh}_G(X,L,W))$ by the above proposition. If $I\in {\rm Inj}_G(X,L,W)$, $J\in {\rm Inj}_G(X,L,W-V)$ and $F\in{\rm lfr}_G(X,L,V)$, then $\mathcal{H}om(F,I)\in{\rm Inj}_G(X,L,W-V)$ and $J\otimes R\in{\rm Inj}_G(X,L,W)$. Hence by Lemma \ref{inj} and Lemma \ref{lfr} we obtain an adjoint pair,
$$(-)\otimes^{\bf L} R \dashv {\bf R}\mathcal{H}om(R,-).$$

\vspace{2mm}
\begin{dfn}
Let $\mathcal{O}_X:=(0\rightarrow\mathcal{O}_X\rightarrow0)\in{\rm coh}_G(X,L,0)$. Then we define functors
$$(-)^{\vee}:=\mathcal{H}om(-,\mathcal{O}_X):{\rm coh}_G(X,L,W)^{\rm op}\rightarrow{\rm coh}_G(X,L,-W)$$
 $$(-)^{{\bf L}\vee}:= {\bf R}\mathcal{H}om(-,\mathcal{O}_X):{\rm Dcoh}_G(X,L,W)^{\rm op}\rightarrow{\rm Dcoh}_G(X,L,-W).$$
\end{dfn}

\vspace{2mm}
\begin{lem}[\cite{bfk} Lemma 3.30, 3.11]
The functor, 
$$(-)^{{\bf L}\vee}:{\rm Dcoh}_G(X,L,W)^{\rm op}\rightarrow{\rm Dcoh}_G(X,L,-W)$$
is an equivalence. 

For $F\in{\rm lfr}_G(X,L,W)$, we have an isomorphism of functors,
$$F^{\vee}\otimes(-)\cong\mathcal{H}om(F,-).$$ 
For $E\in {\rm Dcoh}_G(X,L,W)$, there is an isomorphism of functors,
$$E^{{\bf L}\vee}\otimes^{\bf L}(-)\cong {\bf R}\mathcal{H}om(E,-).$$
\end{lem}

\vspace{2mm}
\begin{lem}\label{tensor preserve ps}
Let $E\in{\rm Dcoh}_G(X,L,V)$ and  $F\in{\rm Dcoh}_G(X,L,W)$. Let $Y$ be a smooth quasi-projective variety and let $f:X\rightarrow Y$ be a morphism. If $E$ has a $f$-proper support, both of $E\otimes^{\bf L}F$ and ${\bf R}\mathcal{H}om(E,F)$ have  $f$-proper supports. In particular, if $E$ has a $f$-proper support, so is $E^{{\bf L}{\vee}}$.

\begin{proof}
By the assumption, there exists an object $E^{\text{\tiny{\textbullet}}}\in {\rm D^{b}}(Z^0({\rm coh}_G(X,L,V)))$ such that ${\rm Tot}(E^{\text{\tiny{\textbullet}}})\cong E$ and the morphism ${\rm Supp}(E^{\text{\tiny{\textbullet}}})\rightarrow Y$ is proper. Since $X$ is smooth and for any $M\in {\rm coh}_G(X)$ there exists a locally free equivariant sheaf $P$ and  a surjection $P\rightarrow M$, there exists an object $P^{\text{\tiny{\textbullet}}}\in{\rm D^{b}}(Z^0({\rm lfr}_G(X,L,V)))$ which is isomorphic to $E^{\text{\tiny{\textbullet}}}$ in ${\rm D^{b}}(Z^0({\rm coh}_G(X,L,V)))$. Then we have 
$$E\otimes^{\bf L}F\cong {\rm Tot}(P^{\text{\tiny{\textbullet}}})\otimes F\cong {\rm Tot}(P^{\text{\tiny{\textbullet}}}\otimes F)$$
and 
$${\bf R}\mathcal{H}om(E,F)\cong \mathcal{H}om( {\rm Tot}(P^{\text{\tiny{\textbullet}}}),F)\cong  {\rm Tot}(\mathcal{H}om(P^{\text{\tiny{\textbullet}}},F)).$$
Hence it is sufficient to prove that closed subsets ${\rm Supp}(P^{\text{\tiny{\textbullet}}}\otimes F)$ and ${\rm Supp}(\mathcal{H}om(P^{\text{\tiny{\textbullet}}},F))$ are contained in ${\rm Supp}(P^{\text{\tiny{\textbullet}}})$. But this  follows from equalities 
$${\rm Supp}(P^{\text{\tiny{\textbullet}}}\otimes F)=\bigcup_{i,j=0,1}{\rm Supp}(P^{\text{\tiny{\textbullet}}}_i\otimes F_j)$$
$${\rm Supp}(\mathcal{H}om(P^{\text{\tiny{\textbullet}}},F))=\bigcup_{k,l=0,1}{\rm Supp}(\mathcal{H}om(P^{\text{\tiny{\textbullet}}}_k,F_l))$$
 and the fact that  for $A^{\text{\tiny{\textbullet}}},B^{\text{\tiny{\textbullet}}}\in{\rm D^b}(X)$, we have  ${\rm Supp}(A^{\text{\tiny{\textbullet}}}\otimes^{\bf L}B^{\text{\tiny{\textbullet}}})\subset{\rm Supp}(A^{\text{\tiny{\textbullet}}})$ and ${\rm Supp}({\bf R}\mathcal{H}om(A^{\text{\tiny{\textbullet}}},B^{\text{\tiny{\textbullet}}}))\subset{\rm Supp}(A^{\text{\tiny{\textbullet}}})$.\end{proof}
\end{lem}

\vspace{2mm}
\subsubsection{Projection formula, flat base change and Grothendieck duality}

Let $X$ and $Y$ be smooth quasi-projective varieties and let $G$ be an affine algebraic group acting on $X$ and $Y$. Let $f:X\rightarrow Y$ be an equivariant morphism. Take $L\in {\rm Pic}_G(Y)$ and $W\in H^0(Y,L)^G$. 

\vspace{2mm}
The following proposition is a version of projection formula for  factorization categories.

\begin{prop}[\cite{bfk} Lemma 3.38] For $E\in {\rm DQcoh}_G(Y,L,W)$ and $F\in {\rm DQcoh}_G(X,f^*L,f^*W)$, we have a natural isomorphism of exact functors,
$${\bf R}f_* F\otimes^{\bf L}E\cong{\bf R}f_*(F\otimes^{\bf L}{\bf L}f^*E).$$
\end{prop}

\vspace{2mm}
Let $Z$ be another smooth quasi-projective variety with $G$-action and let $u:Z\rightarrow Y$ be an equivariant flat morphism. Consider the fiber product $W:=X\times_{Y} Z$,
 \[\xymatrix{
W \ar[r]^{f'} \ar[d]_{u'}& Z\ar[d]^{u}\\
 X\ar[r]^{f}&Y. 
  }\]
  
\vspace{2mm}
\begin{lem}[cf. \cite{bfk} Lemma 2.19]
We have a natural isomorphism of functors between coherent sheaves,
$$u^*\circ f_*\cong f'_*\circ u'^*: {\rm Qcoh}_G(X)\rightarrow{\rm Qcoh}_G(Z).$$
\end{lem}

\vspace{2mm}
 Note that the above natural isomorphism of functors is a cwp-functor morphism. By Lemma \ref{cwp-functor morphism} (2), we have an induced isomorphism of functors between factorizations,
 $$u^*\circ f_*\cong f'_*\circ u'^*: {\rm Qcoh}_G(X,f^*L,f^*W)\rightarrow{\rm Qcoh}_G(Z,u^*L,u^*W).$$  
Since this isomorphism of dg-functors is of degree zero, there is a natural isomorphism of exact functors, 
 $$u^*\circ f_*\cong f'_*\circ u'^*: H^0({\rm Qcoh}_G(X,f^*L,f^*W))\rightarrow H^0({\rm Qcoh}_G(Z,u^*L,u^*W)).$$ 
 Since $u$ and $u'$ are flat, we have ${\bf L}u^*\cong u^*$ and  ${\bf L}u'^*\cong u'^*$. For $E\in{\rm DQcoh}_G(X,f^*L,f^*W)$, let  $I\in{\rm Inj}_G(X,f^*L,f^*W)$ be an object which is quasi-isomorphic to $E$. Then we have
 $$u^*\circ{\bf R}f_*(E)\cong u^*(f_*(I))\cong f'_*(u'^*(I)).$$
 By the second property of right derived functor in Theorem \ref{derived functor} and Lemma \ref{f-acyclic}, we see that $u'^*(I)$ is right $f_*$-acyclic, which implies $f'_*(u'^*(I))\cong {\bf R}f'_*(u'^*(I))$.  Hence we have the following:

 \vspace{2mm}
 \begin{lem}
 We have a natural isomorphism of functors
 $$u^*\circ{\bf R}f_*\cong{\bf R}f'_*\circ u'^*: {\rm DQcoh}_G(X,f^*L,f^*W)\rightarrow{\rm DQcoh}_G(Z,u^*L,u^*W).$$
 \end{lem}
 
  \vspace{2mm}

\begin{dfn}
Let $\varphi:X_1\rightarrow X_2$ be a equivariant morphism of smooth $G$-varieties. We define the \textbf{relative dualizing bundle} $\omega_{\varphi}\in{\rm Pic}_G(X_1)$ as
$$\omega_{\varphi}:=\omega_{X_1}\otimes \varphi^*\omega_{X_2}^{\vee},$$
where $\omega_{X_i}\in{\rm Pic}_G(X_i)$ is the canonical bundle on $X_i$ with tautological equivariant structure.
\end{dfn}

In  \cite{efi-posi}, Positselski proved a version of Grothendieck duality for derived factorization categories.
In the following we give an immediate consequence of the Positselski's result.

 \begin{thm}[cf. \cite{efi-posi} Theorem 3.8]\label{grothendieck}
 If $f$ is proper, direct image ${\bf R}f_*:{\rm DQcoh}(X,f^*L,f^*W)\rightarrow{\rm DQcoh}(Y,L,W)$ has a right adjoint functor $f^!:{\rm DQcoh}(Y,L,W)\rightarrow{\rm DQcoh}(X,f^*L,f^*W)$. An explicit form of the functor $f^!$ is  the following:
  $$f^!(-)\cong{\bf L}f^*(-)\otimes \omega_f[{\rm dim}(X)-{\rm dim}(Y)],$$
where the tensor product on the right hand side is given by Definition \ref{complextensor}.
 \begin{proof}
 Let $D_Y^{\text{\tiny{\textbullet}}}$ be a dualizing complex on $Y$ and write $D_X^{\text{\tiny{\textbullet}}}:=f^{+}D_Y^{\text{\tiny{\textbullet}}}$, where $f^{+}$ is a right adjoint functor of the direct image ${\bf R}f_*:{\rm D^{b}}(X)\rightarrow{\rm D^{b}}(Y)$ of derived categories of coherent sheaves.  By  \cite[Theorem 3.8]{efi-posi}, for any object $E\in{\rm D^{co}Qcoh}(X,f^*L,f^*W)$ and an object $F\in{\rm D^{co}Qcoh}(Y,L,W)$ whose components $F_i$ are flat sheaves, we have an isomorphism
 $${\rm Hom}_{{\rm D^{co}Qcoh}(Y,L,W)}({\bf R}f_*E, F\otimes D_Y^{\text{\tiny{\textbullet}}})\cong{\rm Hom}_{{\rm D^{co}Qcoh}(X,f^*L,f^*W)}(E,f^*(F)\otimes  D_X^{\text{\tiny{\textbullet}}}).$$
Since $X$ and $Y$ are smooth, co-derived factorization categories are equal to absolute derived factorization categories  by Remark \ref{co=abs}, and   the structure sheaf $\mathcal{O}_Y$ is quasi-isomorphic to a dualizing complex. We have  $f^{+}\mathcal{O}_Y\cong \omega_X\otimes f^*\omega_Y^{-1}[{\rm dim}(X)-{\rm dim}(Y)]$. Since for any object of ${\rm DQcoh}(Y,L,W)$ is isomorphic to an object whose components are locally free, in particular, flat, we obtain the theorem.
 \end{proof}
 
  \end{thm}

\vspace{2mm}
\subsubsection{Extension by zero}
 In this section we construct a relative left adjoint functor $i_{!}$ of the inverse image $i^*$ of an open immersion $i$.

 Let $U$ be an open subvariety of $X$ and let $i:U\hookrightarrow X$ be the open immersion. In what follows we don't consider $G$-actions until the next section.

\begin{dfn}
For $F\in{\rm coh}(U)$, let  $\overline{F}$ be coherent sheaf on $X$ such that $\overline{F}|_U\cong F$. Let $\widetilde{\mathbb{Z}_{\geq 0}}$ be the category such that ${\rm Ob}(\widetilde{\mathbb{Z}_{\geq 0}})=\mathbb{Z}_{\geq 0}$ and whose sets of morphisms are defined as follows:
$${\rm Hom}_{\widetilde{\mathbb{Z}_{\geq 0}}}(n,m)=\begin{cases}
 \emptyset & {\rm if} \hspace{2mm}n<m\\
 \{\geq^n_m\} &  {\rm if}  \hspace{2mm}n\geq m
 \end{cases}$$
Then we define an object $i_!(F)\in{\rm Pro}({\rm coh}(X))$  as a functor $i_!(F):\widetilde{\mathbb{Z}_{\geq 0}}\rightarrow{\rm coh}(X)$ defined by $$i_!(F)(n):=\mathcal{I}^n\overline{F},$$ where $\mathcal{I}$ is the ideal sheaf defining the complement $X\setminus U$. Since the object $i_!(F)$ doesn't depend on the choice of an extension $\overline{F}$ by the following Lemma \ref{independent}, this gives an exact functor
$$i_!:{\rm coh}(U)\rightarrow{\rm Pro}({\rm coh}(X)).$$
The functor $i_!$ is called \textbf{the extension by zero} of $i$. We also denote by $i_!$ the composition
$${\rm coh}(U)\xrightarrow{i_!}{\rm Pro}({\rm coh}(X))\hookrightarrow{\rm Pro}({\rm Qcoh}(X)).$$
\end{dfn}
 
 \begin{lem}\label{independent}
 Let $F\in{\rm coh}(U)$ be an coherent sheaf on $U$, and let $N\in{\rm coh}(X)$ and $M\in{\rm Qcoh}(X)$ be subsheaves of $i_*(F)\in{\rm Qcoh}(X)$. If $i^*(N)$ is contained in $i^*(M)$, then there is a positive integer $n$ such that $\mathcal{I}^n N$ is contained in $M$.
 \begin{proof}
Since we can take finite affine covering, it is enough to  prove it for the case $X={\rm Spec}(A)$ and $U={\rm Spec}(A_f)$ for some ring $A$ and an element $f\in A$.  Then $\mathcal{I}$  corresponds to the ideal $I=\langle f \rangle$ of $A$ generated by $f$. We consider $F$, $N$ and $M$ as corresponding modules. Let $\{x_k\}_{1\leq k \leq r}\subset N$ be a generator of $N$. Since $i^*(N)=N\otimes_A A_f$ is contained in $i^*(M)=M\otimes_A A_f$, for each $k$, there is an element $y_k\in M$ and $n_k\geq0$ such that $x_k\otimes 1=y_k\otimes 1/{f^{n_k}}$ in $i_*(F)\otimes A_f$. This implies that $f^{n_k}x_k=y_k\in M$, since $i_*(F)\otimes A_f\cong F$. Set $n:={\rm max}\{ n_k | 1\leq k \leq r\}$. Then we have $I^nN\subset M$.
  \end{proof}
 \end{lem}
 
 \vspace{2mm}
 Deligne proved that the extension by zero $i_!$ is a relative left adjoint to the inverse image $i^*$.
 
\vspace{2mm}
\begin{prop}[cf.\,\cite{deli} Proposition 4]\label{deligne adjoint}
For any $F\in{\rm coh}(U)$ and $(G:\mathcal{I}\rightarrow {\rm Qcoh}(X))\in{\rm Pro}({\rm Qcoh}(X))$, we have an isomorphism
$${\rm Hom}_{{\rm Pro}({\rm Qcoh}(X))}(i_!(F),G)\cong{\rm Hom}_{{\rm Pro}({\rm Qcoh}(U))}(J(F),{\rm Pro}(i^*)(G)),$$
where $J:{\rm coh}(U)\rightarrow{\rm Pro}({\rm Qcoh}(U))$ is the natural inclusion. 

\begin{proof}
This is shown as follows;
\begin{align*}
{\rm Hom}_{{\rm Pro}({\rm Qcoh}(X))}(i_!(F),G)&=\varprojlim_{i\in\mathcal{I}}{\rm Hom}_{{\rm Pro}({\rm Qcoh}(X))}(i_!(F),G(i))\\
&\cong\varprojlim_{i\in\mathcal{I}}{\rm Hom}_{{\rm Qcoh}(U)}(F,i^*G(i))\\
&= {\rm Hom}_{{\rm Pro}({\rm Qcoh}(U))}(J(F),{\rm Pro}(i^*)(G)),
\end{align*}
where the isomorphism in the second line follows from \cite[Proposition 4]{deli}.
\end{proof}
\end{prop}

\vspace{3mm}
Let $L\in{\rm Pic}(X)$ and let $W\in H^0(X,L)$. Then $({\rm Pro}(L),{\rm Pro}(W))$ is a potential of ${\rm Pro}({\rm Qcoh}(X))$ and ${\rm Pro}({\rm coh}(X))$. We denote their factorization categories by 
$${\rm Qcoh}_{\rm Pro}(X,L,W):=\mathfrak{F}({\rm Pro}({\rm Qcoh}(X)),{\rm Pro}(L),{\rm Pro}(W))$$
$${\rm coh}_{\rm Pro}(X,L,W):=\mathfrak{F}({\rm Pro}({\rm coh}(X)),{\rm Pro}(L),{\rm Pro}(W)).$$
The extension by zero $i_!$ is compatible with potentials with respect to $(L|_U,W|_U)$ and $({\rm Pro}(L),{\rm Pro}(W))$. Hence the functor $i_!$ induces a dg-functor 
$$i_!:{\rm coh}(U,L|_U,W|_U)\rightarrow{\rm coh}_{\rm Pro}(X,L,W).$$
Since $i_!:{\rm coh}(U)\rightarrow{\rm Pro}({\rm coh}(X))$ is an exact functor of abelian categories, $i_!$ preserves acyclic objects. Hence $i_!:H^0({\rm coh}(U,L|_U,W|_U))\rightarrow H^0({\rm coh_{Pro}}(X,L,W))$ naturally induces an exact functor
$$i_!:{\rm Dcoh}(U,L|_U,W|_U)\rightarrow{\rm D}{\rm coh_{Pro}}(X,L,W).$$
On the other hand, there is a natural functor ${\rm D}{\rm coh_{Pro}}(X,L,W)\rightarrow{\rm Pro}({\rm D}{\rm coh}(X,L,W))$. Composing it with the embedding ${\rm Pro}({\rm D}{\rm coh}(X,L,W))\rightarrow{\rm Pro}({\rm D}{\rm Qcoh}(X,L,W))$ and  $i_!:{\rm Dcoh}(U,L|_U,W|_U)\rightarrow{\rm D}{\rm coh_{Pro}}(X,L,W)$, we construct a functor 
$$i_!:{\rm Dcoh}(U,L|_U,W|_U)\rightarrow{\rm Pro}({\rm D}{\rm Qcoh}(X,L,W)),$$
which is also denoted by the same notation $i_!$.

\begin{prop}\label{left adjoint of open}
$(1)$ The dg-functor $i_!:{\rm coh}(U,L|_U,W|_U)\rightarrow{\rm Qcoh_{Pro}}(X,L,W)$ is $J$-left adjoint to ${\rm Pro}(i^*):{\rm Qcoh_{Pro}}(X,L,W)\rightarrow{\rm Qcoh_{Pro}}(U,L|_U,W|_U)$, where $J$ is the natural embedding functor $J:{\rm coh}(U,L|_U,W|_U)\rightarrow{\rm Qcoh_{Pro}}(U,L|_U,W|_U)$.\\
$(2)$ For any $E\in{\rm Dcoh}(U,L|_U,W|_U)$ and $F\in{\rm DQcoh}(X,L,W)$, we have an isomorphism 
$${\rm Hom}_{{\rm Pro}({\rm D}{\rm Qcoh}(X,L,W))}(i_!(E),\iota_{\rm Pro}(F))\cong{\rm Hom}_{{\rm DQcoh}(U,L|_U,W|_U)}(E,i^*(F)).$$
\begin{proof}

(1) Consider the following diagram,
 \[\xymatrix{
{\rm coh}(U)\ar[rr]^{i_!} \ar[rd]_{J}&&{\rm Pro}({\rm Qcoh}(X))\ar[ld]^{{\rm Pro}(i^*)}\\
&{\rm Pro}({\rm Qcoh}(U))&
  }\]
where $J$ is the natural embedding. Then Proposition \ref{deligne adjoint} implies that $i_!$ is $J$-left  adjoint to ${\rm Pro}(i^*)$ (see Definition \ref{rel ad}). Hence, (1) holds since the front adjunction $J\rightarrow {\rm Pro}(i^*)\circ i_!$ is a cwp-functor morphism.\\
(2) Let $F\xrightarrow{\sim}J$ be an isomorphism in ${\rm D}{\rm Qcoh}(X,L,W)$ such that the  components of $I$ are injective quasi-coherent sheaves. Then $i^*(I)$ is an object whose components are injective quasi-coherent sheaves on $U$. By (1) and  Lemma \ref{inj}, the right hand side of the desired isomorphism is isomorphic to $$H^0(\{\varinjlim_{n\in\mathbb{Z}_{\geq0}}{\rm Hom}_{{\rm Qcoh}(X,L,W)}(\mathcal{I}^n\overline{E},I)\}^{\text{\tiny{\textbullet}}}).$$ Since taking direct limit is an exact functor, the above  abelian group  is isomorphic to $$\varinjlim_{n\in\mathbb{Z}_{\geq0}}H^0({\rm Hom}_{{\rm Qcoh}(X,L,W)}(\mathcal{I}^n\overline{E},I)^{\text{\tiny{\textbullet}}}),$$
which is isomorphic to the left hand side of the desired isomorphism by Lemma \ref{inj} again.
\end{proof}
\end{prop}

\vspace{3mm}
 For later use, we will extend the extension by zero $i_!:{\rm coh}(U)\rightarrow {\rm Pro}({\rm Qcoh}(X))$ to a functor defined on ${\rm Qcoh}(U)$. To do it, we need the following lemma.

 \begin{lem}[\cite{deli}, Proposition 2]\label{indcoh}
 Let $Y$ be a  Noetherian scheme, and  let $F\in {\rm Qcoh}(Y)$ be a quasi-coherent sheaf.  Denote by $\{ F_k\}_{k\in K}$ the family of all coherent subsheaves of $F$. Let $\theta(F)\in{\rm Ind}({\rm coh}(Y))$ be a functor given by 
 $$
\begin{array}{ccc}
 K & \stackrel{\theta(F)}{\longrightarrow} & {\rm coh}(Y) \\
\rotatebox{90}{$\in$} & & \rotatebox{90}{$\in$} \\
k & \longmapsto & F_k
\end{array}
$$
Then $\theta(-)$ gives an exact equivalence $$\theta:{\rm Qcoh}(Y)\xrightarrow{\sim} {\rm Ind}({\rm coh}(Y)).$$
  \end{lem}
 
  \vspace{2mm}
   \begin{dfn}\label{extended extension by zero}
   We define an exact functor $$i_{\#}:{\rm Qcoh}(U)\rightarrow {\rm Ind}({\rm Pro}({\rm Qcoh}(X)))$$ as the compositions $${\rm Qcoh}(U)\xrightarrow{\theta} {\rm Ind}({\rm coh}(U))\xrightarrow{{\rm Ind}(i_!)}{\rm Ind}({\rm Pro}({\rm coh}(X)))\hookrightarrow{\rm Ind}({\rm Pro}({\rm Qcoh}(X)))$$
   \end{dfn}
  
  \vspace{2mm}
  \begin{rem}
  By the construction of $i_{\#}$, we have a natural isomorphism of functors $$i_{\#}|_{{\rm coh}(U)}\cong \iota_{\rm Ind} i_!.$$
    \end{rem}
   
   \vspace{2mm}
  The following lemma will be necessary to prove Lemma \ref{completing lemma}.

   \begin{lem}\label{gammas} The notation is the same as above.
   \begin{itemize} 
   \item[$(1)$]  We have a natural  functor morphism
      $$\gamma_{\#}: i_{\#}\rightarrow \iota_{\rm InPr} i_*,$$
   where $\iota_{\rm InPr}:{\rm Qcoh}(X)\rightarrow {\rm Ind}({\rm Pro}({\rm Qcoh}(X)))$. Restricting $\gamma_{\#}$, we obtain a natural functor morphism
   $$\gamma_{!}: i_{!}\rightarrow \iota_{\rm Pro} i_*.$$
   such that  $\iota_{\rm Ind}\gamma_!=\gamma_{\#}|_{{\rm coh}(U)}$.
      \item[$(2)$] Consider the following cartesian square:
    \[\xymatrix{
    V:=U\times_X Y\ar[rr]^{j}\ar[d]_{q}&&Y\ar[d]^{p}\\
    U\ar[rr]^{i}&&X
  }\]
We have a morphism between functors from ${\rm coh}(V)$ to ${\rm Ind}({\rm Pro}({\rm Qcoh}(Y)))$
 $$\lambda : j_{\#}q^*q_*\rightarrow\iota_{\rm Ind}{\rm Pro}(p^*p_*)j_!$$
 such that the following diagram is commutative:
 \[\xymatrix{
 \iota_{\rm InPr}j_*q^*q_*\ar[r]^{\iota_{\rm InPr}\delta}&\iota_{\rm InPr}p^*p_*j_* \ar@{=}[r]&\iota_{\rm Ind}{\rm Pro}(p^*p_*)\iota_{\rm Pro}j_*\\
 j_{\#}q^*q_*\ar[rr]^{\lambda}\ar[u]^{\gamma  q^*q_*}&&\iota_{\rm Ind}{\rm Pro}(p^*p_*)j_!\ar[u]_{\iota_{\rm Ind}{\rm Pro}(p^*p_*)\gamma|_{{\rm coh}(V)}}
     }\]
     where $\delta:j_*q^*q_*\xrightarrow{\sim}p^*p_*j_*$ is a natural isomorphism of functors.
   \end{itemize}
   \begin{proof}
   (1) Let $F\in {\rm Qcoh}(U)$ be a quasi-coherent sheaf on $U$, and let $\{F_k\}_{k\in K}$ be the family of all coherent subsheaves of $F$. By definition, $i_{\#}(F):K\rightarrow{\rm Pro}({\rm Qcoh}(U))$ is a functor given by $i_{\#}(F)(k)=i_!(F_k)$, and the object $i_!(F_k)\in{\rm Pro}({\rm Qcoh}(U))$ is the  functor given by $$\mathbb{Z}_{\geq0}\ni n\mapsto \mathcal{I}^n \overline{F_k}\in {\rm coh}(U),$$
   where $\overline{F_k}$ is a coherent subsheaf of $i_*(F_k)$ such that $i^*(\overline{F_k})\cong F_k$. Hence, the natural inclusion $ \overline{F_k}\hookrightarrow i_*(F)$ gives a morphism of functors
   $$\gamma: i_{\#}\rightarrow \iota_{\rm InPr} i_*.$$
   
   (2) For $F\in{\rm coh}(V)$, we will define a morphism $\lambda(F): j_{\#}q^*q_*(F)\rightarrow\iota_{\rm Ind}{\rm Pro}(p^*p_*)j_!(F)$. Let $\{E_k\}_{k\in K}$ be the family of all coherent subsheaves of $q^*q_*(F)$. Then the object $j_{\#}q^*q_*(F)\in {\rm Ind}({\rm Pro}({\rm Qcoh}(Y)))$ is given by the following functor
    $$
\begin{array}{ccc}
 K & \stackrel{}{\longrightarrow} & {\rm Pro}({\rm Qcoh}(Y)) \\
\rotatebox{90}{$\in$} & & \rotatebox{90}{$\in$} \\
k & \longmapsto & j_!(E_k)
\end{array}
$$
 In order to define a morphism $\lambda(F): j_{\#}q^*q_*(F)\rightarrow\iota_{\rm Ind}{\rm Pro}(p^*p_*)j_!(F)$, it is enough to give a family of morphisms $\{\lambda(F)_k:j_!(E_k)\rightarrow{\rm Pro}(p^*p_*)j_!(F)\}_{k\in K}$ in ${\rm Pro}({\rm Qcoh}(Y))$ such that for any inclusion $v:E_k\hookrightarrow E_l$, the equation $\lambda(F)_k=\lambda(F)_l j_!(v)$ holds.  Let $\mathcal{J}$ be the ideal sheaf defining $Y\setminus V$, and let $\overline{E_k}$ and $\overline{F}$ be  coherent subsheaves of $j_*(E_k)$ and $j_*(F)$ with $j^*(\overline{E_k})\cong E_k$ and $j^*(\overline{F})\cong F$ respectively. Then the  object $j_!(E_k)$ and ${\rm Pro}(p^*p_*)j_!(F)$ are  the following functors 
 $$
\begin{array}{ccc}
 \mathbb{Z}_{\geq 0} & \stackrel{j_!(E_k)}{\longrightarrow} & {\rm Qcoh}(Y) \\
\rotatebox{90}{$\in$} & & \rotatebox{90}{$\in$} \\
n & \longmapsto & \mathcal{J}^n\overline{E_k}
\end{array}
$$
 $$
\begin{array}{ccc}
 \mathbb{Z}_{\geq 0} & \stackrel{{\rm Pro}(p^*p_*)j_!(F)}{\longrightarrow} & {\rm Qcoh}(Y) \\
\rotatebox{90}{$\in$} & & \rotatebox{90}{$\in$} \\
m & \longmapsto & \mathcal{J}^mp^*p_*(\overline{F})
\end{array}
$$
$\overline{E_k}$ is contained in $j_*q^*q_*(F)$ and $p^*p_*(\overline{F})$ can be considered as a subsheaf of $j_*q^*q_*(F)$ via the isomorphism $\delta(F):j_*q^*q_*(F)\xrightarrow{\sim}p^*p_*j_*(F)$. Since $j^*\overline{E_k}\cong E_k$ is contained in $j^*p^*p_*(\overline{F})\cong q^*q_*(F)$,  there is a positive integer $N$ such that $\mathcal{J}^N\overline{E_k}$ is a subsheaf of $p^*p_*(\overline{F})$ by Lemma \ref{independent}. Let  $\theta_{\lambda(F)_k}:\mathbb{Z}_{\geq 0}\ni n \mapsto n+N \in \mathbb{Z}_{\geq 0}$ be a map, and let 
 $\lambda(F)_k^n:j_!(E_k)(n+N)\rightarrow {\rm Pro}(p^*p_*)j_!(F)(n)$ be a morphism induced by the inclusion $\mathcal{J}^N\overline{E_k}\hookrightarrow p^*p_*(\overline{F})$. If we define a morphism $\lambda(F)_k:j_!(E_k)\rightarrow{\rm Pro}(p^*p_*)j_!(F)$ as a map of systems $(\{\lambda(F)_k^n\}_{n\in\mathbb{Z}_{\geq 0}}, \theta_{\lambda(F)_k})$ for each $k\in K$, then the family $\{\lambda(F)_k\}_{k\in K}$ defines a morphism $\lambda(F): j_{\#}q^*q_*(F)\rightarrow\iota_{\rm Ind}{\rm Pro}(p^*p_*)j_!(F)$, and this gives a functor morphism $$\lambda : j_{\#}q^*q_*\rightarrow\iota_{\rm Ind}{\rm Pro}(p^*p_*)j_!.$$ 
The commutativity of the diagram follows since $\gamma$ is induced by  natural inclusions, and $\lambda$ is also induced by natural inclusions via $\delta$. 
  \end{proof}
   \end{lem}

\vspace{2mm}
\subsubsection{Integral functor for factorization}
 
Let $X_1$ and $X_2$ be smooth quasi-project varieties with actions of affine algebraic group $G$. Take a character $\chi$ of $G$, and let $\mathcal{O}_i(\chi)$ be the corresponding equivariant line bundle on $X_i$.  Let $W_i\in{\rm H}^0(X_i,\mathcal{O}_i(\chi))^{G}$ be a $G$-invariant section. Then the corresponding regular function $W_i:X_i\rightarrow\mathbb{A}^1$ is $\chi$-semi invariant, i.e. $W(g\cdot x)=\chi(g)\cdot W(x)$ for any $g\in G$ and $x\in X_i$. Denote by $\pi_i:X_1\times X_2\rightarrow X_i$ the projection for each $i=1,2$.

Throughout this section 4.2.6, dropping the script $L$ from notation, we  write ${\rm Dcoh}_G(-,*)$ instead of ${\rm Dcoh}_G(-,L,*)$, because all equivariant line bundles in this section are the one corresponding to the character $\chi$.

\begin{dfn}
For  $P\in{\rm DQcoh}_G(X_1\times X_2,\pi_2^*W_2-\pi_1^*W_1)$, we define \textbf{the integral functor} $\Phi_P$ \textbf{with kernel} $P$ as 
$$\Phi_P:={\bf R}\pi_{2*}(\pi_1^*(-)\otimes^{\bf L}P):{\rm DQcoh}_G(X_1,W_1)\rightarrow{\rm DQcoh}_G(X_2,W_2).$$
\end{dfn}

\vspace{2mm}
\begin{rem}
If $Q\in{\rm Dcoh}_G(X_1\times X_2,\pi_2^*W_2-\pi_1^*W_1)$ has a $\pi_2$-proper support, then $\Phi_Q$ maps an object in ${\rm Dcoh}_G(X_1,W_1)$ to an object in ${\rm D_{coh}Qcoh}_G(X_2,W_2)$. We also denote by $\Phi_Q$ the following composition
$${\rm Dcoh}_G(X_1,W_1)\xrightarrow{\Phi_Q}{\rm D_{coh}Qcoh}_G(X_2,W_2)\xrightarrow{\sim}{\rm Dcoh}_G(X_2,W_2).$$
\end{rem}

\vspace{2mm}
For an object $P\in{\rm Dcoh}_G(X_1\times X_2,\pi_2^*W_2-\pi_1^*W_1)$, we define objects $P_R$ and $P_L$ in ${\rm Dcoh}_G(X_1\times X_2,\pi_1^*W_1-\pi_2^*W_2)$ as
$$P_R:=P^{{\bf L}\vee}\otimes\pi_1^*\omega_{X_1}[{\rm dim}(X_1)]$$
$$P_L:=P^{{\bf L}\vee}\otimes\pi_2^*\omega_{X_2}[{\rm dim}(X_2)].$$
If $G$ is trivial, we see that there are  relative adjoint pairs of integral functors.

\vspace{2mm}
\begin{prop}\label{integral adjoint}
Let $P\in{\rm Dcoh}(X_1\times X_2,\pi_2^*W_2-\pi_1^*W_1)$ be an  object which has a $\pi_1$-proper support. Then for any objects $E\in{\rm DQcoh}(X_1,W_1)$ and  $F\in{\rm Dcoh}(X_2,W_2)$, we have an isomorphism
$${\rm Hom}_{{\rm DQcoh}(X_2,W_2)}(F, \Phi_P(E))\cong{\rm Hom}_{{\rm DQcoh}(X_1,W_1)}(\Phi_{P_L}(F),E).$$
In particular, if $P$ has a  $\pi_2$-proper support, then $\Phi_{P_L}:{\rm Dcoh}(X_2,W_2)\rightarrow{\rm Dcoh}(X_1,W_1)$ (resp. $\Phi_{P_R}$) is a left (resp. right) adjoint functor of $\Phi_P:{\rm Dcoh}(X_1,W_1)\rightarrow{\rm Dcoh}(X_2,W_2)$. 

\begin{proof}
Since we already have the adjunction $\pi_2^*\dashv {\bf R}\pi_{2*}$, it is enough to obtain the following isomorphism 
$${\rm Hom}_{{\rm DQcoh}(X_1\times X_2,\pi_2^*W_2)}(D,\pi_{1}^*E\otimes^{\bf L}P)\cong{\rm Hom}_{{\rm DQcoh}(X_1,W_1)}({\bf R}\pi_{1*}(D\otimes^{\bf L}P_L),E)$$
for any objects $D\in {\rm Dcoh}(X_1\times X_2, \pi_2^*W_2)$ and $E\in{\rm DQcoh}(X_1,W_1)$. This is proved in a similar way to the proof of \cite[Lemma 4]{logvinenko}. Compactify $X_2$ and denote by $\overline{X_2}$ a smooth proper variety containing $X_2$ as an open subvariety. Let $\iota:X_1\times X_2\hookrightarrow X_1\times \overline{X_2}$ be the open immersion, and let $\overline{\pi_1}: X_1\times \overline{X_2}\rightarrow X_1$ be  the projection. Then $\pi_1=\overline{\pi_1}\circ\iota$ and $\overline{\pi_1}$ is a proper morphism. By Theorem \ref{grothendieck} and Proposition \ref{left adjoint of open}, we obtain the following isomorphism:
$${\rm Hom}_{{\rm DQcoh}(X_1\times X_2,\pi_2^*W_2)}(D,\pi_{1}^*E\otimes^{\bf L}P)\cong{\rm Hom}_{{\rm Pro}({\rm DQcoh}(X_1,W_1))}({\bf R}\overline{\pi_{1}}_*(\iota_{!}(D\otimes^{\bf L}P^{{\bf L}\vee})\otimes\omega_{\overline{\pi_1}}[{\rm dim}(X_2)]),E)$$
Since the object $P^{{\bf L}\vee}$  has a $\pi_1$-proper support, there exists an object $P^{\text{\tiny{\textbullet}}}\in {\rm D^b}(Z^0({\rm coh}(X_1\times X_2,\pi_1^*W_1-\pi_2^*W_2)))$ such that $P^{{\bf L}\vee}\cong {\rm Tot}(P^{\text{\tiny{\textbullet}}})$ and ${\rm Supp}(P^{\text{\tiny{\textbullet}}})$ is $\pi_1$-proper, in particular, $\iota$-proper. 
By a similar reasoning to one of  \cite[Lemma 4]{logvinenko}, we see that there is an isomorphism 
$$\iota_{!}((-)\otimes^{\bf L}P^{\text{\tiny{\textbullet}}})\xrightarrow{\sim}\iota_{*}((-)\otimes^{\bf L}P^{\text{\tiny{\textbullet}}})$$
of functors from ${\rm D^b}(Z^0({\rm coh}(X_1\times X_2,\pi_2^*W_2)))$ to ${\rm Pro}({\rm D^b}(Z^0({\rm Qcoh}(X_1\times \overline{X_2},\pi_1^*W_1))))$. 
By taking totalizations of the above isomorphism, we obtain an isomorphism of functors
 $$\iota_{!}((-)\otimes^{\bf L}P^{{\bf L}\vee})\xrightarrow{\sim}\iota_{*}((-)\otimes^{\bf L}P^{{\bf L}\vee}).$$
Hence, we have an isomorphism ${\bf R}\overline{\pi_{1}}_*(\iota_{!}(D\otimes^{\bf L}P^{{\bf L}\vee})\otimes\omega_{\overline{\pi_1}}[{\rm dim}(X_2)])\cong{\bf R}\pi_{1*}(D\otimes^{\bf L}P_L)$. 
If $P$ has a $\pi_2$-proper support, the integral functor $\Phi_P$ maps ${\rm Dcoh}(X_1,W_1)$ to ${\rm Dcoh}(X_2,W_2)$, and $\Phi_{P_L}$ maps ${\rm Dcoh}(X_2,W_2)$ to ${\rm Dcoh}(X_1,W_1)$ since $P_L$ has a $\pi_1$-proper support by Lemma \ref{tensor preserve ps}. Hence we have $\Phi_{P_L}\dashv\Phi_P$. Since $(P_L)_R\cong P$, we obtain the other adjunction $\Phi_P \dashv \Phi_{P_R}$.
\end{proof}

\end{prop}

\vspace{2mm}
We will show that the composition of integral functors is also an integral functor. Let  $X_3$ be another smooth quasi-projective $G$-variety, $\mathcal{O}_3(\chi)$ be the equivariant line bundle corresponding to the character $\chi$, and $W_3\in H^0(X_3,\mathcal{O}_3(\chi))^G$ be an invariant section. We define morphisms of varieties by the following diagram;
\[\xymatrix{
&&X_1\times X_2\ar[lld]_{\pi_1}\ar[rrd]^{\pi_2}&&\\
X_1&&&&X_2\\
&&X_1\times X_2\times X_3\ar[uu]^{\pi_{12}}\ar[lld]_{\pi_{13}}\ar[rrd]^{\pi_{23}}&&\\
X_1\times X_3\ar[uu]^{q_1}\ar[rrd]_{q_3}&&&&X_2\times X_3\ar[uu]_{p_2}\ar[lld]^{p_3}\\
&&X_3&&
  }\]
where all morphisms are projections.
 For two objects
$$P\in{\rm Dcoh}_G(X_1\times X_2,\pi_2^*W_2-\pi_1^*W_1)$$
$$Q\in{\rm Dcoh}_G(X_2\times X_3,p_3^*W_3-p_2^*W_2),$$
we set another object
$$P\star Q:=\pi_{13*}(\pi_{12}^*P\otimes^{\bf L}\pi_{23}^*Q)\in{\rm Dcoh}_G(X_1\times X_3, q_3^*W_3-q_1^*W_1).$$
For two complexes $P^{\text{\tiny{\textbullet}}}\in{\rm D^b}({\rm coh}X_1\times X_2)$ and $Q^{\text{\tiny{\textbullet}}}\in {\rm D^b}({\rm coh}X_2\times X_3)$, we also define another object $$P^{\text{\tiny{\textbullet}}}\star Q^{\text{\tiny{\textbullet}}}\in{\rm D^b}({\rm coh}X_1\times X_3)$$
in the same manner. Then we have the following:

\vspace{2mm}
\begin{prop}\label{composition of integral2}The notation is the same as above.
The composition of integral functors 
$${\rm Dcoh}_G(X_1,W_1)\xrightarrow{\Phi_P}{\rm Dcoh}_G(X_2,W_2)\xrightarrow{\Phi_Q}{\rm Dcoh}_G(X_3,W_3)$$
is isomorphic to the following integral functor
$${\rm Dcoh}_G(X_1,W_1)\xrightarrow{\Phi_{P\star Q}}{\rm Dcoh}_G(X_3,W_3).$$
The similar result holds for integral functors of derived categories of coherent sheaves.
\begin{proof}
For the proof of the  result for derived categories of coherent sheaves, see \cite[Proposition 5.10]{huybrechts}, for example. We can prove the result for derived factorization categories in the same way.
\end{proof}
\end{prop}

\vspace{3mm}

\subsection{Comonads induced by restriction and induction functors}

We construct restriction and  induction functors and study comonads induced by these functors.

We continue to assume that the quasi-projective variety $X$ is smooth. 
Let $G\times^l X$ and $G\times ^d X$ be the varieties $G\times X$ with different $G$-actions which are defined as follows;

$$
\begin{array}{ccc}
G\times G\times^l X & \stackrel{}{\longrightarrow} & G\times^l X \\
\rotatebox{90}{$\in$} & & \rotatebox{90}{$\in$} \\
(g,g',x) & \longmapsto & (gg',x)
\end{array}
$$
and
$$
\begin{array}{ccc}
G\times G\times^d X & \stackrel{}{\longrightarrow} & G\times^d X \\
\rotatebox{90}{$\in$} & & \rotatebox{90}{$\in$} \\
(g,g',x) & \longmapsto & (gg',gx).
\end{array}
$$
\vspace{2mm}
Then the following morphisms 
$$
\begin{array}{ccc}
\varphi: G\times^l X & \stackrel{}{\longrightarrow} & G\times^d X \\
\rotatebox{90}{$\in$} & & \rotatebox{90}{$\in$} \\
(g,x) & \longmapsto & (g,gx)
\end{array}
$$ and 
$$
\begin{array}{ccc}
\pi: G\times^d X & \stackrel{}{\longrightarrow} & X \\
\rotatebox{90}{$\in$} & & \rotatebox{90}{$\in$} \\
(g,x) & \longmapsto & x
\end{array}
$$
  are $G$-equivariant. The action $\sigma:G\times X\rightarrow X$ on $X$ is  the composition $\pi\circ\varphi$.

Let $\iota:X\rightarrow G\times X$ be a morphism defined by  $$X\ni x\longmapsto (e,x)\in G\times X.$$ We define an exact functor $\iota^*:{\rm Qcoh}_G(G\times^lX)\rightarrow{\rm Qcoh}X$ as 
$$\begin{array}{ccc}
{\rm Qcoh}_G(G\times^lX) & \stackrel{}{\longrightarrow} &{\rm Qcoh}X \\
\rotatebox{90}{$\in$} & & \rotatebox{90}{$\in$} \\
(\mathcal{F},\theta) & \longmapsto & \iota^*\mathcal{F}.
\end{array}
$$

\vspace{2mm}
\begin{lem}\label{equivalences}
\begin{itemize}
\item[$(1)$]
The functor $\iota^*:{\rm Qcoh}_G(G\times^lX)\rightarrow{\rm Qcoh}X$ is an equivalence.
\item[$(2)$]
The functors $\varphi^*:{\rm Qcoh}_G(G\times^dX)\rightarrow{\rm Qcoh}_G(G\times^lX)$ and $\varphi_*:{\rm Qcoh}_G(G\times^lX)\rightarrow{\rm Qcoh}_G(G\times^dX)$ are equivalences.
\item[$(3)$]The functors $\pi^*:{\rm Qcoh}_G(X)\rightarrow{\rm Qcoh}_G(G\times^dX)$ and $\pi_*:{\rm Qcoh}_G(G\times^dX)\rightarrow{\rm Qcoh}_G(X)$ are exact functors.
\end{itemize}
\begin{proof}
(1)This is a special case of \cite[Lemma 1.3.]{thomason}\\
(2)The morphism $\varphi$ is an isomorphism.\\
(3)Since $\pi$ is smooth, in particular flat, and affine, $\pi^*$ and $\pi_*$ are  exact functors.
\end{proof}
\end{lem}

\vspace{2mm}
\begin{dfn}
We define \textbf{the restriction functor} ${\rm Res}_G:{\rm Qcoh}_G(X)\rightarrow{\rm Qcoh}X$ and  \textbf{the induction functor} ${\rm Ind}_G:{\rm Qcoh}X\rightarrow{\rm Qcoh}_G(X)$ as 
$${\rm Res}_G:=\iota^*\circ\sigma^*\hspace{5mm}{\rm and}\hspace{5mm}{\rm Ind}_G:=\sigma_*\circ(\iota^*)^{-1}.$$
\end{dfn}

\vspace{2mm}
\begin{rem}
Note that the restriction functor ${\rm Res}_G:{\rm Qcoh}_G(X)\rightarrow{\rm Qcoh}X$ is isomorphic to the forgetful functor, i.e. ${\rm Res}_G(\mathcal{F},\theta)\cong \mathcal{F}$.

\end{rem}

\vspace{2mm}
Let $L$ be an invertible $G$-equivariant  sheaf, and let $W$ be an invariant section of $L$. Then the pair $(L,W)$ defines potentials of ${\rm Qcoh}_G(X)$ and ${\rm Qcoh}X$. Since the functors ${\rm Res}_G$ and ${\rm Ind}_G$ are cwp-functors, these functors induce  functors of factorization categories
$${\rm Res}_G:{\rm Qcoh}_G(X,L,W)\rightarrow{\rm Qcoh}(X,L,W)$$
$${\rm Ind}_G:{\rm Qcoh}(X,L,W)\rightarrow{\rm Qcoh}_G(X,L,W)$$
Since $\iota^*$ is an equivalence, the adjoint pair $\sigma^*\dashv\sigma_*$ induces the  adjoint pair
$${\rm Res}_G\dashv{\rm Ind}_{G}.$$
Since the functors ${\rm Res}_G$ and ${\rm Ind}_G$ are exact functors, we obtain the exact functor of derived factorization categories
$$\Pi^*_G:={\rm Res}_G:{\rm DQcoh}_G(X,L,W)\rightarrow{\rm DQcoh}(X,L,W)$$
$$\Pi_{G*}:={\rm Ind}_G:{\rm DQcoh}(X,L,W)\rightarrow{\rm DQcoh}_G(X,L,W),$$
and these defines an adjoint pair
$$\Pi_G:=(\hspace{1mm}\Pi^*_G\dashv\Pi_{G*}\hspace{1mm}).$$

\begin{rem} The functor $\Pi_G^*$ sends  objects in ${\rm Dcoh}_G(X,L,W)$ to objects in ${\rm Dcoh}(X,L,W)$. But the functor $\Pi_{G*}$ does not preserve coherentness of components of factorizations.
\end{rem}

\begin{dfn}
We define a comonad $\mathbb{T}_G$ on ${\rm DQcoh}(X,L,W)$ induced by $G$-action as the one induced by the adjoint pair $\Pi_G:=(\Pi^*_G\dashv\Pi_{G*})$;
$$\mathbb{T}_G:=\mathbb{T}(\Pi_G),$$
where the notation is the same as in Example \ref{adj}. Denote  by $\Gamma_G$ is the comparison functor of the adjoint pair $\Pi_G:=(\Pi^*_G\dashv\Pi_{G*})$,
$$\Gamma_G:{\rm DQcoh}_G(X,L,W)\rightarrow{\rm DQcoh}(X,L,W)_{\mathbb{T}_G}.$$ 
\end{dfn}

\vspace{2mm} We recall the definition of (linearly) reductiveness of algebraic groups.

\begin{dfn}
Let $H$ be an affine algebraic group over a field $K$.
\begin{itemize}
\item[$(1)$]
 $H$ is called \textbf{reductive} if the radical of $H$ is a torus.
 \item[$(2)$] $H$ is called \textbf{linearly reductive} if every rational representations of $H$ over $K$ is completely reducible.
 \end{itemize}
\end{dfn}

\begin{prop}[\cite{git} Appendix A]
Let $H$ be an affine algebraic group over   a field $K$ of characteristic zero.
Then $H$ is reductive if and only if linearly reductive.
\end{prop}

\begin{lem}\label{reductive} If $G$ is linearly reductive, then the adjunction morphism ${\rm id}\rightarrow \Pi_{G*}\Pi^*_G$ is a split mono. In particular, the comparison functor $\Gamma_G:{\rm DQcoh}_G(X,L,W)\rightarrow{\rm DQcoh}(X,L,W)_{\mathbb{T}_G}$ is an equivalence.

\begin{proof}
Since the adjunction morphism ${\rm id}\rightarrow \Pi_{G*}\Pi^*_G$ coincide with the adjunction morphism  ${\rm id}\rightarrow \pi_*\pi^*$, and the morphism $E\rightarrow \pi_*\pi^*E$ is equal to the morphism $E\otimes(\mathcal{O}_X\rightarrow \pi_*\pi^*\mathcal{O}_X$) via the projection formula, it is enough to show that $\mathcal{O}_X\rightarrow \pi_*\pi^*\mathcal{O}_X$ is split mono. Since $G$ is linearly reductive, the homomorphism $k\rightarrow \mathcal{O}_G(G)$ of $G$-modules is split mono. This means that the adjunction $\mathcal{O}_{{\rm Spec}(k)}\rightarrow p_*p^*\mathcal{O}_{{\rm Spec}(k)}$ is split mono, where $p:G\rightarrow{\rm Spec}(k)$ is the morphism defining the base space.
Hence by the cartesian square,
\[\xymatrix{
G\times X\ar[r]^{\pi}\ar[d]&X\ar[d]\\
G\ar[r]^{p}&{\rm Spec}(k),
  }\]
we see that $\mathcal{O}_X\rightarrow \pi_*\pi^*\mathcal{O}_X$ is also a split mono. The latter statement  follows from Proposition \ref{comparison theorem} and Lemma \ref{cocomplete}.
\end{proof}

\end{lem}

\begin{lem}\label{linearizable2}
Let $X'$ be another smooth quasi-projective variety with $G$-action and let $f:X'\rightarrow X$ be a $G$-equivariant morphism. Let $\mathbb{T}_{G}'=\mathbb{T}(\Pi'_G)$ be the comonad on  ${\rm DQcoh}(X',f^*L,f^*W)$ induced by its $G$-action. Let $P\in{\rm Dcoh}_G(X,L,W)$ be an object.

Then there exist functor isomorphisms $\Omega^*:{\bf R}f_*{\Pi'}_G^*\xrightarrow{\sim}\Pi_G^*{\bf R}f_*$ and $\Omega_*:{\bf R}f_*\Pi'_{G*}\xrightarrow{\sim}\Pi_{G*}{\bf R}f_*$ such that the following diagrams are commutative;
\[\xymatrix{
\Pi^*_G{\bf R}f_*\Pi'_{G*}\ar[rr]^{\Pi^*_G\Omega_*}&&\Pi_G^*\Pi_{G*}{\bf R}f_*\ar[d]^{\varepsilon{\bf R}f_*}& &{\bf R}f_*\ar[rr]^{\eta{\bf R}f_*}\ar[d]_{{\bf R}f_*\eta'}&&\Pi_{G*}\Pi_G^*{\bf R}f_*\\
{\bf R}f_*{\Pi'}_G^*\Pi'_{G*}\ar[rr]^{{\bf R}f_*\varepsilon'}\ar[u]^{\Omega^*\Pi'_{G*}}&&{\bf R}f_*& &{\bf R}f_*\Pi'_{G*}{\Pi'}^*_G\ar[rr]^{\Omega_*{\Pi'}_G^*}&&\Pi_{G*}{\bf R}f_*{\Pi'}_G^*\ar[u]_{\Pi_{G*}\Omega^*},
 }\]
 where $\varepsilon$, $\varepsilon'$, $\eta$ and $\eta'$ are  adjunction morphisms.
In particular, the direct image ${\bf R}f_*:{\rm DQcoh}(X',f^*L,f^*W)\rightarrow{\rm DQcoh}(X,L,W)$ is a linearizable functor with respect to $\mathbb{T}'_G$ and $\mathbb{T}_G$ with a linearization $\Omega:=\Pi_G^*\Omega_*\circ\Omega^*\Pi'_{G*}$, and the following diagram is commutative:
\[\xymatrix{
{\rm DQcoh}(X',f^*L,f^*W)_{\mathbb{T}'_G}\ar[rr]^{{\bf R}f_{*\Omega}}&&{\rm DQcoh}(X,L,W)_{\mathbb{T}_G}\\
{\rm DQcoh}_G(X',f^*L,f^*W)\ar[u]^{\Gamma'_G}\ar[rr]^{{\bf R}f_*}&&{\rm DQcoh}_G(X,L,W)\ar[u]_{\Gamma_G}.
  }\]

 The similar results hold for the inverse image ${\bf L}f^*:{\rm DQcoh}(X,L,W)\rightarrow{\rm DQcoh}(X',f^*L,f^*W)$ and the tensor product $(-)\otimes^{\bf L}\Pi^*_GP:{\rm DQcoh}(X,L,V)\rightarrow{\rm DQcoh}(X,L,V+W)$.
\begin{proof}
We only give a proof for the case of the direct image. Let $\pi:G\times X\rightarrow X$ and $\pi':G\times X'\rightarrow X'$ be natural projections, and set $\overline{f}:={\rm id}_G\times f:G\times X'\rightarrow G\times X$. By Lemma \ref{equivalences} (1) and (2), we have the following equivalences;
$$\Phi:{\rm DQcoh}_G(G\times^dX,\pi^*L,\pi^*W)\xrightarrow{\sim}{\rm DQcoh}(X,L,W)$$
$$\Phi':{\rm DQcoh}_G(G\times^dX',{\pi'}^*f^*L,{\pi'}^*f^*W)\xrightarrow{\sim}{\rm DQcoh}(X',f^*L,f^*W)$$
such that $\Pi_G^*\cong\Phi\pi^*$, $\Pi_{G*}\cong\pi_{*}\Phi^{-1}$, ${\Pi'_G}^*\cong\Phi'{\pi'}^*$ and $\Pi'_{G*}\cong\pi'_*{\Phi'}^{-1}$. By the following cartesian square 
\[\xymatrix{
G\times X'\ar[r]^{\overline{f}}\ar[d]_{\pi'}& G\times X\ar[d]^{\pi}\\
X'\ar[r]^{f}&X,
  }\]
we have isomorphisms of functors between categories of quasi-coherent sheaves;
$$\overline{\omega}^*:\overline{f}_*{\pi'}^*\xrightarrow{\sim}\pi^*f_*\hspace{3mm}{\rm and}\hspace{3mm}\overline{\omega}_*:f_*\pi'_*\xrightarrow{\sim}\pi_*\overline{f}_*.$$
By  easy computation, we see that the following diagrams are commutative;
\[\xymatrix{
\pi^*f_*\pi'_*\ar[rr]^{\pi^*\overline{\omega}_*}&&\pi^*\pi_*\overline{f}_*\ar[d]^{\varepsilon_{\pi}\overline{f}_*}& &f_*\ar[rr]^{\eta_{\pi}f_*}\ar[d]_{f_*\eta_{\pi'}}&&\pi_*\pi^*f_*\\
\overline{f}_*{\pi'}^*\pi'_*\ar[rr]^{\overline{f}_*\varepsilon_{\pi'}}\ar[u]^{\overline{\omega}^*\pi'_*}&&\overline{f}_*& &f_*\pi'_*{\pi'}^*\ar[rr]^{\overline{\omega}_*{\pi'}^*}&&\pi_*\overline{f}_*{\pi'}^*\ar[u]_{\pi_*\overline{\omega}^*}.
 }\]
 Since the functor morphisms in the above diagrams are cwp-functor morphisms, taking $H^0(\mathfrak{F}(-))$, we obtain  similar isomorphisms of functors between homotopy categories of factorization categories,  and similar commutative diagrams of morphisms of  exact functors between homotopy categories of factorization categories. These isomorphisms of functors and commutative diagrams induce  isomorphisms of functors between derived factorization categories 
 $$\overline{\Omega}^*:{\bf R}\overline{f}_*{\pi'}^*\xrightarrow{\sim}\pi^*{\bf R}f_*\hspace{3mm}{\rm and}\hspace{3mm}\overline{\Omega}_*:{\bf R}f_*\pi'_*\xrightarrow{\sim}\pi_*{\bf R}\overline{f}_*$$
and the following commutative diagrams
\[\xymatrix{
\pi^*{\bf R}f_*\pi'_*\ar[rr]^{\pi^*\overline{\Omega}_*}&&\pi^*\pi_*{\bf R}\overline{f}_*\ar[d]^{\varepsilon_{\pi}{\bf R}\overline{f}_*}& &{\bf R}f_*\ar[rr]^{\eta_{\pi}{\bf R}f_*}\ar[d]_{{\bf R}f_*\eta_{\pi'}}&&\pi_*\pi^*{\bf R}f_*\\
{\bf R}\overline{f}_*{\pi'}^*\pi'_*\ar[rr]^{{\bf R}\overline{f}_*\varepsilon_{\pi'}}\ar[u]^{\overline{\Omega}^*\pi'_*}&&{\bf R}\overline{f}_*& &{\bf R}f_*\pi'_*{\pi'}^*\ar[rr]^{\overline{\Omega}_*{\pi'}^*}&&\pi_*{\bf R}\overline{f}_*{\pi'}^*\ar[u]_{\pi_*\overline{\Omega}^*}.
 }\]
Since ${\bf R}f_*\Phi'\cong\Phi{\bf R}\overline{f}_*$, applying the equivalences $\Phi$ and ${\Phi'}^{-1}$ to the above functor isomorphisms and commutative diagrams,  we obtain the desired functor isomorphisms and commutative diagrams.

The results for the inverse image and the tensor product are proved similarly.
\end{proof}

\end{lem}







\vspace{5mm}
\section{Main results}

At first, we prepare notation used throughout this section. Let $X_1$ and $X_2$ be smooth quasi-project varieties  with actions of reductive affine algebraic group $G$ over an algebraically closed field $k$ of characteristic zero. For a character $\chi: G\rightarrow\mathbb{G}_m$ of $G$, take  $\chi$-semi invariant regular functions $W_i\in{\rm H}^0(X_i,\mathcal{O}_{X_i}(\chi))^{G}$ on $X_i$. Let $\pi_i:X_1\times X_2\rightarrow X_i$ and $q_i:X_1\times_{\mathbb{A}^1} X_2\rightarrow X_i$ be the projections and let $j:X_1\times_{\mathbb{A}^1} X_2\hookrightarrow X_1\times X_2$ be the embedding. We have the following commutative diagram:
\[\xymatrix{
&&X_1\times X_2\ar[lldd]_{\pi_1}\ar[rrdd]^{\pi_2}&&\\
&&X_1\times_{\mathbb{A}^1} X_2\ar[u]^{j}\ar[lld]^{q_1}\ar[rrd]_{q_2}&&\\
X_1\ar[rrd]_{W_1}&&&&X_2\ar[lld]^{W_2}\\
&&\mathbb{A}^1&&
}\]
Abbreviating $\mathcal{O}_{X_i}(\chi)$, we  write 
$${\rm DQcoh}_G(X_i,W_i):={\rm DQcoh}_G(X_i,\mathcal{O}_{X_i}(\chi),W_i)$$
and
$${\rm Dcoh}_G(X_i,W_i):={\rm Dcoh}_G(X_i,\mathcal{O}_{X_i}(\chi),W_i).$$

\vspace{2mm}
\subsection{Equivariantization}
The action of $G$ on $X_i$ induces an adjoint pair 
$$\Pi_i:=(\hspace{1mm}\Pi^*_i\dashv\Pi_{i*}\hspace{1mm}),$$
where the functor $\Pi^*_i$ and $\Pi_{i*}$ are given by restriction and induction functors respectively;
$$\Pi^*_i:={\rm Res}_G:{\rm DQcoh}_G(X_i,W_i)\rightarrow{\rm DQcoh}(X_i,W_i)$$
$$\Pi_{i*}:={\rm Ind}_G:{\rm DQcoh}(X_i,W_i)\rightarrow{\rm DQcoh}_G(X_i,W_i).$$
Denote by $\mathbb{T}_i$ be the comonad on ${\rm DQcoh}(X_i,W_i)$ induced by the adjoint pair $\Pi_i=(\Pi^*_i\dashv\Pi_{i*})$ and let $\Gamma_i$ be the comparison functor of the adjoint pair $\Pi_i$,
$$\Gamma_i:{\rm DQcoh}_G(X_i,W_i)\rightarrow{\rm DQcoh}(X_i,W_i)_{\mathbb{T}_i}.$$

\vspace{2mm}
\begin{thm}\label{equivariantization}
Let $P_G\in{\rm Dcoh}_G(X_1\times X_2, \pi_2^*W_2-\pi_1^*W_1)$ be an object and set $P:={\rm Res}_G(P_G)\in{\rm Dcoh}(X_1\times X_2,  \pi_2^*W_2-\pi_1^*W_1)$. Assume that $P$ has a $\pi_i$-proper support $(i=1,2)$. If the integral functor $\Phi_{P}:{\rm Dcoh}(X_1,W_1)\rightarrow{\rm Dcoh}(X_2,W_2)$ is fully faithful (resp. equivalence), then the integral functor $\Phi_{{P_G}}:{\rm Dcoh}_G(X_1,W_1)\rightarrow{\rm Dcoh}_G(X_2,W_2)$ is fully faithful (resp. equivalence).

\begin{proof}
Set $(P_L)_G:=(P_G)^{\bf L{\vee}}\otimes \pi_2^*\omega_{X_2}[{\rm dim}(X_2)]\in{\rm Dcoh}_G(X_1\times X_2,\pi_1^*W_1-\pi_2^*W_2)$ and $P_L:={\rm Res}_G((P_L)_G)\in{\rm Dcoh}(X_1\times X_2,\pi_1^*W_1-\pi_2^*W_2)$. Then we have the following diagram.

 \[\xymatrix{
{\rm Dcoh}(X_1,W_1)\ar@/^/[rrrr]^{\Phi_{P}} \ar@{^{(}-{>}}[rd]&&&&{\rm Dcoh}(X_2,W_2)\ar@/^/[llll]^{\Phi_{P_L}} \ar@{^{(}-{>}}[ld]\\
 &{\rm DQcoh}(X_1,W_1)\ar@/^/[rr]^{\Phi_P}\ar@/^/[d]^{\Pi_{1*}}  && {\rm DQcoh}(X_2,W_2) \ar@/^/[d]^{\Pi_{2*}} \ar@/^/[ll]^{\Phi_{P_L}}& \\
&{\rm DQcoh}_G(X_1,W_1) \ar@/^/[u]^{\Pi_1^*} \ar@/^/[rr]^{\Phi_{P_G}}&& {\rm DQcoh}_G(X_2,W_2) \ar@/^/[u]^{\Pi_2^*} \ar@/^/[ll]^{\Phi_{(P_L)_G}}&\\
{\rm Dcoh}_G(X_1,W_1)\ar@{^{(}-{>}}[ru] \ar@/^/[rrrr]^{\Phi_{P_G}}  \ar[uuu]^{\Pi_1^*}&&&&{\rm Dcoh}_G(X_2,W_2) \ar@/^/[llll]^{\Phi_{(P_L)_G}} \ar@{^{(}-{>}}[lu] \ar[uuu]_{\Pi_2^*}
 }\]
 By Lemma \ref{linearizable2}, there exist functor isomorphisms $\Omega^*:\Phi_P\Pi^*_1\xrightarrow{\sim}\Pi^*_2\Phi_{P_G}$, $\Omega_*:\Phi_{P_G}\Pi_{1*}\xrightarrow{\sim}\Pi_{2*}\Phi_P$, $\Omega_{L}^*:\Phi_{P_L}\Pi_2^*\xrightarrow{\sim}\Pi_1^*\Phi_{(P_L)_G}$ and $\Omega_{L*}:\Phi_{(P_L)_G}\Pi_{2*}\xrightarrow{\sim}\Pi_{1*}\Phi_{P_L}$ such that the diagrams corresponding to (i) and (ii) in Lemma \ref{compatible lemma}, namely the following diagrams,  are commutative.
 \[\xymatrix{
\Pi_2^*\Phi_{P_G}\Pi_{1*}\ar[rr]^{\Pi^*_2\Omega_*}&&\Pi_2^*\Pi_{2*}\Phi_P\ar[d]^{\varepsilon_2\Phi_P}& &\Phi_{P_G}\ar[rr]^{\eta_2\Phi_{P_G}}\ar[d]_{\Phi_{P_G}\eta_1}&&\Pi_{2*}\Pi_2^*\Phi_{P_G}\\
\Phi_P{\Pi}_1^*\Pi_{1*}\ar[rr]^{\Phi_P\varepsilon_1}\ar[u]^{\Omega^*\Pi_{1*}}&&\Phi_P& &\Phi_{P_G}\Pi_{1*}{\Pi}^*_1\ar[rr]^{\Omega_*{\Pi}_1^*}&&\Pi_{2*}\Phi_P{\Pi}_1^*\ar[u]_{\Pi_{2*}\Omega^*},
 }\]
  \[\xymatrix{
\Pi_1^*\Phi_{(P_L)_G}\Pi_{2*}\ar[rr]^{\Pi^*_1\Omega_{L*}}&&\Pi_1^*\Pi_{1*}\Phi_{P_L}\ar[d]^{\varepsilon_1\Phi_{P_L}}& &\Phi_{(P_L)_G}\ar[rr]^{\eta_1\Phi_{(P_L)_G}}\ar[d]_{\Phi_{(P_L)_G}\eta_2}&&\Pi_{1*}\Pi_1^*\Phi_{(P_L)_G}\\
\Phi_{P_L}{\Pi}_2^*\Pi_{2*}\ar[rr]^{\Phi_{P_L}\varepsilon_2}\ar[u]^{\Omega_L^*\Pi_{2*}}&&\Phi_{P_L}& &\Phi_{(P_L)_G}\Pi_{2*}{\Pi}^*_2\ar[rr]^{\Omega_{L*}{\Pi}_2^*}&&\Pi_{1*}\Phi_{P_L}{\Pi}_2^*\ar[u]_{\Pi_{1*}\Omega_L^*},
 }\]
 where $\varepsilon_i$ and $\eta_i$ are adjunction morphisms of the adjoint pair $(\Pi_i^*\dashv\Pi_{i*})$.
 Combining Lemma \ref{main lemma} with Proposition \ref{integral adjoint} and Lemma \ref{reductive}, we see that if $\Phi_{P}:{\rm Dcoh}(X_1,W_1)\rightarrow{\rm Dcoh}(X_2,W_2)$ is fully faithful, then $\Phi_{{P_G}}:{\rm Dcoh}_G(X_1,W_1)\rightarrow{\rm Dcoh}_G(X_2,W_2)$ is also fully faithful. 
 
 Assume $\Phi_P$ is an equivalence. Then, $\Phi_{P_L}$ is fully faithful functor. Applying the above argument to $\Phi_{P_L}$, we see that $\Phi_{(P_L)_G}$ is also fully faithful. Set $\Omega:=\Pi_2^*\Omega_*\circ \Omega^*\Pi_{1*}$ and $\Omega_L:=\Pi_1^*\Omega_{L*}\circ\Omega_L^*\Pi_{2*}$. By Lemma \ref{main lemma}, we see  that $\Phi_{P_G}$ is an equivalence by the following Lemma \ref{completing lemma}.
 \end{proof}
\end{thm}

 \begin{lem}\label{completing lemma}
 With notation same as above,   the following diagram of functors from ${\rm Dcoh}_G(X_1,W_1)$ to ${\rm DQcoh}(X_1,W_1)$ is commutative;
 \[\xymatrix{
 (\ast): &\Phi_{P_L}\Phi_P\Pi_1^* \ar[rr]^{\Phi_{P_L}\Phi_P\Pi_1^*\eta_1}\ar[d]_{\omega \Pi_1^*}&&\Phi_{P_L}\Phi_P\Pi_1^*\Pi_{1*}\Pi_1^* \ar[r]^{\Phi_{P_L}\Omega \Pi_1^*}&\Phi_{P_L}\Pi_2^*\Pi_{2*}\Phi_P\Pi_1^* \ar[r]^{\Omega_L\Phi_P\Pi_1^*} &\Pi_1^*\Pi_{1*}\Phi_{P_L}\Phi_P\Pi_1^* \ar[d]^{\Pi_1^*\Pi_{1*}\omega \Pi_1^*}\\
 & \Pi_1^* \ar[rrrr]^{\Pi_1^*\eta_1} &&&& \Pi_1^*\Pi_{1*}\Pi_1^*, 
 }\]
 where $\omega:\Phi_{P_L}\Phi_P\rightarrow{\rm id}_{{\rm Dcoh}(X_1,W_1)}$ is the adjunction morphism of $(\Phi_{P_L}\dashv \Phi_P).$\qed
\end{lem}

We will prove the above lemma in the next section.

\subsection{Proof of Lemma \ref{completing lemma}}
In what follows, we will prove the above Lemma \ref{completing lemma}. Since it seems difficult to verify the commutativity of the diagram $(\ast)$ directly, we will replace it with another diagram $(\ast)'$, and  decompose the diagram $(\ast)'$ into several diagrams whose commutativity are easier to verify.

Take   a smooth proper variety $\overline{X_2}$ containing $X_2$ as an open subvariety as in the proof of Lemma \ref{integral adjoint}.  Let  $i:X_1\times X_2\hookrightarrow X_1\times\overline{X_2}$ be the open immersion, and let $\overline{\pi_1}:X_1\times\overline{X_2}\rightarrow X_1 $ be the natural projection.  Denote natural projections by $p_i:G\times X_i\rightarrow X_i$, $p_{12}:G\times X_1\times X_2\rightarrow X_1\times X_2$ and $\overline{p_{12}}:G\times X_1\times \overline{X_2}\rightarrow X_1\times \overline{X_2}$, and set 
$$\pi_i':=1_G\times \pi_i:G\times X_1\times X_2\rightarrow G\times X_i$$
$$i':=1_G\times i:G\times X_1\times X_2\rightarrow G\times X_1\times \overline{X_2}$$
$$\overline{\pi_1}':=1_G\times\overline{\pi_1}: G\times X_1\times \overline{X_2}\rightarrow G\times X_1.$$
Then objects $Q_G:=p_{12}^*P_G\in{\rm Dcoh}_G(G\times X_1\times X_2,{\pi_2'}^*p_2^*W_2-{\pi_1'}^*p_1^*W_1)$ and $(Q_L)_G:=p_{12}^*(P_L)_G\in{\rm Dcoh}_G(G\times X_1\times X_2,{\pi_1'}^*p_1^*W_1-{\pi_2'}^*p_2^*W_2)$ define functors
$$
\begin{array}{ccc}
\Psi_{Q_G}:{\rm DQcoh}_G(G\times X_1,p_1^*W_1) & \stackrel{}{\longrightarrow} & {\rm DQcoh}_G(G\times X_2,p_2^*W_2)\\
&&\\
F & \longmapsto & {\pi_2'}_*(\pi_1'^*(F)\otimes^{\bf{L}}Q_G)
\end{array}
$$
and
$$
\begin{array}{ccc}
\Psi_{(Q_L)_G}:{\rm DQcoh}_G(G\times X_2,p_2^*W_2) & \stackrel{}{\longrightarrow} & {\rm DQcoh}_G(G\times X_1,p_1^*W_1)\\
&&\\
E & \longmapsto & {\pi_1'}_*(\pi_2'^*(E)\otimes^{\bf{L}}(Q_L)_G).
\end{array}
$$
Note that $Q_G$  has a $\pi_i'$-proper support ($i=1,2$). Hence the functors $\Psi_{Q_G}$ and $\Psi_{(Q_L)_G}$ preserve coherent factorizations.

Similarly, the objects $Q:={\rm Res}_G(Q_G)\in{\rm Dcoh}(G\times X_1\times X_2,{\pi_2'}^*p_2^*W_2-{\pi_1'}^*p_1^*W_1)$ and $Q_L:={\rm Res}_G((Q_L)_G)\in{\rm Dcoh}(G\times X_1\times X_2,{\pi_1'}^*p_1^*W_1-{\pi_2'}^*p_2^*W_2)$ defines functors
$$
\begin{array}{ccc}
\Psi_{Q}:{\rm DQcoh}(G\times X_1,p_1^*W_1) & \stackrel{}{\longrightarrow} & {\rm DQcoh}(G\times X_2,p_2^*W_2)\\
&&\\
F & \longmapsto & {\pi_2'}_*(\pi_1'^*(F)\otimes^{\bf{L}}Q)
\end{array}
$$
and
$$
\begin{array}{ccc}
\Psi_{Q_L}:{\rm DQcoh}(G\times X_2,p_2^*W_2) & \stackrel{}{\longrightarrow} & {\rm DQcoh}(G\times X_1,p_1^*W_1)\\
&&\\
E & \longmapsto & {\pi_1'}_*(\pi_2'^*(E)\otimes^{\bf{L}}Q_L).
\end{array}
$$

By Lemma \ref{equivalences} the composition,  
$$\iota^*\circ \varphi^*:{\rm DQcoh}_G(G\times X_1,p_1^*W_1)\xrightarrow{\sim}{\rm DQcoh}(X_1,W_1),$$
is an equivalence, and the following diagrams are commutative,
\[\xymatrix{
{\rm DQcoh}_G(G\times X_1,p_1^*W_1)\ar[rr]^{\Psi_{Q_G}}\ar[d]_{\iota^*\circ\varphi^*}&&{\rm DQcoh}_G(G\times X_2,p_2^*W_2)\ar[d]^{\iota^*\circ\varphi^*}\\
{\rm DQcoh}(X_1,W_1)\ar[rr]^{\Phi_P}&&{\rm DQcoh}(X_2,W_2)
 }\]
and
\[\xymatrix{
{\rm DQcoh}_G(G\times X_2,p_2^*W_2)\ar[rr]^{\Psi_{(Q_L)_G}}\ar[d]_{\iota^*\circ\varphi^*}&&{\rm DQcoh}_G(G\times X_1,p_1^*W_1)\ar[d]^{\iota^*\circ\varphi^*}\\
{\rm DQcoh}(X_2,W_2)\ar[rr]^{\Phi_P}&&{\rm DQcoh}(X_1,W_1).
 }\]
 Let $\Omega':\Psi_{Q_G}p_1^*p_{1*}\xrightarrow{\sim}p_2^*p_{2*}\Psi_{Q_G}$ and $\Omega'_L:\Psi_{(Q_L)_G}p_2^*p_{2*}\xrightarrow{\sim}\Psi_{(Q_L)_G}p_1^*p_{1*}$ be functor isomorphisms induced by the functor isomorphisms $\Omega:\Phi_P\Pi_1^*\Pi_{1*}\xrightarrow{\sim}\Pi_2^*\Pi_{2*}\Phi_P$ and $\Omega':\Phi_{P_L}\Pi_2^*\Pi_{2*}\xrightarrow{\sim}\Pi_1^*\Pi_{1*}\Phi_{P_L}$ via the equivalence $\iota^*\circ\varphi^*$ respectively.
Via the equivalence $\iota^*\circ\varphi^*$, the diagram $(\ast)$ is commutative if and only if the following diagram is commutative;
\[\xymatrix{
 \Psi_{(Q_L)_G}\Psi_{Q_G}p_1^* \ar[rrr]^{\Psi_{(Q_L)_G}\Psi_{Q_G}p_1^*\eta_{p_1}}\ar[d]_{\omega' _Gp_1^*}&&&\Psi_{(Q_L)_G}\Psi_{Q_G} p_1^*p_{1*}p_1^* \ar[r]^{\Psi_{(Q_L)_G}\Omega' p_1^*}&\Psi_{(Q_L)_G}p_2^*p_{2*}\Psi_{Q_G} p_1^* \ar[r]^{{\Omega'}_L\Psi_{Q_G} p_1^*} &p_1^*p_{1*}\Psi_{(Q_L)_G}\Psi_{Q_G} p_1^* \ar[d]^{p_1^*p_{1*}\omega'_G p_1^*}\\
  p_1^* \ar[rrrrr]^{p_1^*\eta_{p_1}} &&&&& p_1^*p_{1*}p_1^*, 
 }\]
 where $\omega'_G:\Psi_{(Q_L)_G}\Psi_{Q_G}\rightarrow{\rm id}_{{\rm Dcoh}_G(G\times X_1,p_1^*W_1)}$ is the adjunction morphism of $(\Psi_{(Q_L)_G}\dashv \Psi_{Q_G}).$ Furthermore, since the restriction functor
 $${\rm Res}_G:{\rm DQcoh}_G(G\times X_2,{\pi'_2}^*p_2^*W_2)\rightarrow{\rm DQcoh}(G\times X_2,{\pi'_2}^*p_2^*W_2)$$
  is faithful functor,  in order to prove that the above diagram is commutative, it is enough to show that the following diagram is commutative,
  \[\xymatrix{
 (\ast)':&\Psi_{Q_L}\Psi_{Q}p_1^* \ar[rr]^{\Psi_{Q_L}\Psi_{Q}p_1^*\eta_{p_1}}\ar[d]_{\omega' p_1^*}&&\Psi_{Q_L}\Psi_{Q} p_1^*p_{1*}p_1^* \ar[r]^{\Psi_{Q_L}\Omega' p_1^*}&\Psi_{Q_L}p_2^*p_{2*}\Psi_{Q} p_1^* \ar[r]^{{\Omega'}_L\Psi_{Q} p_1^*} &p_1^*p_{1*}\Psi_{Q_L}\Psi_{Q} p_1^* \ar[d]^{p_1^*p_{1*}\omega' p_1^*}\\
 & p_1^* \ar[rrrr]^{p_1^*\eta_{p_1}} &&&& p_1^*p_{1*}p_1^*, 
 }\]
 where $\omega':\Psi_{Q_L}\Psi_{Q}\rightarrow{\rm id}_{{\rm Dcoh}(G\times X_1,p_1^*W_1)}$ is the adjunction morphism of $(\Psi_{Q_L}\dashv \Psi_Q).$

To decompose the diagram $(\ast)'$, we give the following:

\begin{lem}\label{decomposition}
Given the following diagram of functors
 \[\xymatrix{
\mathcal{A}_1\ar@/^/[dd]^{P_{1*}} \ar@/^/[rr]^{F_1} &&  \mathcal{A}_2 \ar@/^/[ll]^{G_1}\ar@/^/[dd]^{P_{2*}}\ar@/^/[rr]^{F_2} && \mathcal{A}_3\ar@/^/[ll]^{G_2}\ar@/^/[dd]^{P_{3*}} \\\\
\mathcal{A}_1' \ar@/^/[uu]^{P_1^*} \ar@/^/[rr]^{F'_1}&& \mathcal{A}_2' \ar@/^/[ll]^{G'_1}\ar@/^/[uu]^{P_2^*}\ar@/^/[rr]^{F'_2} &&\mathcal{A'}_3\ar@/^/[uu]^{P_3^*}\ar@/^/[ll]^{G'_2}
 }\]
 and isomorphisms of functors $\Omega_{F_i}:F_iP^*_iP_{i*}\xrightarrow{\sim}P^*_{i+1}P_{i+1*}F_{i}$ and $\Omega_{G_i}:G_iP^*_{i+1}P_{i+1*}\xrightarrow{\sim}P^*_{i}P_{i*}G_{i}$,
 assume the adjunction $(G_i\dashv F_i)$ and $(P_i^*\dashv P_{i*})$  for each $i=1,2$. Set $F:=F_2\circ F_1$ and  $G:=G_1\circ G_2$, and denote by $\omega:GF\rightarrow{\rm id}$  the functor morphism given by the composition $GF=G_1G_2F_2F_1\xrightarrow{G_1\omega_2F_1}G_1F_1\xrightarrow{\omega_1}{\rm id}$, where $\omega_i:G_iF_i\rightarrow{\rm id}$ is the adjunction morphism. Let $\Omega_F:FP_1^*P_{1*}\rightarrow P_3^*P_{3*}F$ and $\Omega_G:GP_3^*P_{3*}\rightarrow P_1^*P_{1*}G$ be the functor isomorphisms induced by $\Omega_{F_i}$ and $\Omega_{G_i}$, i.e. $\Omega_F:=\Omega_{F_2}F_1\circ F_2\Omega_{F_1}$ and $\Omega_{G}:=\Omega_{G_1}G_2\circ G_1\Omega_{G_2}$. For each $i=1,2$, consider the following diagrams of functor morphisms
  \[\xymatrix{
 (\diamondsuit)_i: &G_iF_iP_i^* \ar[rr]^{G_iF_iP_i^*\eta_i}\ar[d]_{\omega_i P_i^*}&&G_iF_iP_i^*P_{i*}P_i^*\ar[r]^{G_i\Omega_{F_i} P_i^*}&G_iP_{i+1}^*P_{i+1*}F_iP_i^*  \ar[r]^{\Omega_{G_i}F_iP_i^*} &P_i^*P_{i*}G_iF_iP_i^* \ar[d]^{P_i^*P_{i*}\omega_i P_i^*}\\
 & P_i^* \ar[rrrr]^{P_i^*\eta_i} &&&& P_i^*P_{i*}P_i^*, 
 }\]
  where $\eta_i:{\rm id}\rightarrow P_{i*}P_i^*$ is the adjunction.
  
   If the above diagrams $(\diamondsuit)_1$ and $(\diamondsuit)_2$ are commutative, and there exist  isomorphisms of functors $\mu:F_1'P_{1*}P_1^*\xrightarrow{\sim} P_{2*}P^*_2F_1'$ and $\nu :F_1P_1^*\xrightarrow{\sim}P_2^*F_1'$ with the following diagrams 
  \[\xymatrix{
(\dagger):&F_1P_1^*P_{1*}P^*_1\ar[rr]^{\Omega_{F_1}P_1^*}\ar[d]_{\nu P_{1*}P_1^*}&&P_2^*P_{2*}F_1P_1^*\ar[d]^{P_2^*P_{2*}\nu}&  (\dagger\dagger):&F_1'P_{1*}P_1^*
\ar[rr]^{\mu}&&P_{2*}P_2^*F_1'\\
&P_2^*F_1'P_{1*}P_1^*\ar[rr]^{P_2^*\mu}&&P_2^*P_{2*}P_2^* F_1'&  & &F_1'\ar[ul]^{F_1'\eta_1}\ar[ur]_{\eta_2F_1'}&
}\]
 commutative, then the following diagram $(\diamondsuit)$ is also commutative.
 \[\xymatrix{
 (\diamondsuit): &GFP_1^* \ar[rr]^{GFP_1^*\eta_1}\ar[d]_{\omega P_1^*}&&GFP_1^*P_{1*}P_1^*\ar[r]^{G\Omega_{F} P_1^*}&GP_{3}^*P_{3*}FP_1^*  \ar[r]^{\Omega_{G}FP_1^*} &P_1^*P_{1*}GFP_1^* \ar[d]^{P_1^*P_{1*}\omega P_1^*}\\
 & P_1^* \ar[rrrr]^{P_1^*\eta_1} &&&& P_1^*P_{1*}P_1^*.
 }\]
 \begin{proof}
 At first, we  show that the following diagram is commutative;
 \[\xymatrix{
 (\clubsuit):&G_2FP_1^* \ar[rr]^{G_2FP_1^*\eta_1}\ar[d]_{\omega_2F_1 P_1^*}&&G_2FT_1P_1^*\ar[r]^{G_2F_2\Omega_{F_1} P_1^*}&G_2F_2T_2F_1P_1^*  \ar[r]^{G_2\Omega_{F_2}F_1P_1^*} &G_2T_3FP_1^*\ar[r]^{\Omega_{G_2}FP_1^*} & T_2G_2FP_1^* \ar[d]^{T_2\omega_2 F_1P_1^*}\\
  &F_1P_1^* \ar[rr]^{F_1P_1^*\eta_1} &&F_1T_1P_1^*\ar[rrr]^{\Omega_{F_1}P_1^*}&&& T_2F_1P_1^*,
 }\] 
 where $T_i:=P_i^*P_{i*}$ for $i=1,2,3$.
 By the commutativity of the diagram $(\dagger)$, the following diagram is commutative;
  \[\xymatrix{
F_1P_1^*\ar[rr]^{F_1P_1^*\eta_1}\ar[d]_{\nu}&&F_1P_1^*P_{1*}P_1^*\ar[rr]^{\Omega_{F_1}P_1^*}\ar[d]^{\nu P_{1*}P_1^*}&&P_2^*P_{2*}F_1P_1^*\ar[d]^{P_2^*P_{2*}\nu}\\
P_2^*F_1'\ar[rr]^{P_2^*F_1'\eta_1}&&P_2^*F_1'P_{1*}P_1^*\ar[rr]^{P_2^*\mu}&&P_2^*P_{2*}P_2^*F_1',
 }\] 
 and we have $P_2^*\mu\circ P_2^*F_1'\eta_1=P_2^*\eta_2F_1'$ by the commutativity of the diagram $(\dagger\dagger)$.
Hence we see that, via the isomorphism of functors $\nu:F_1P_1^*\xrightarrow{\sim}P_2^*F_1'$,  the commutativity of the  diagram $(\clubsuit)$ is equivalent to the commutativity of the following diagram
  \[\xymatrix{
  G_2F_2P_2^*F_1'\ar[rr]^{G_2F_2P_2^*\eta_2F_1'}\ar[d]_{\omega_2 P_2^*F_1'}&&G_2F_2P_2^*P_{2*}P_2^*F_1'\ar[rr]^{G_2\Omega_{F_2}P_2^*F_1'}&&G_2P_3^*P_{3*}F_2P_2^*F_1'\ar[rr]^{\Omega_{G_2}F_2P_2^*F_1'}&&P_2^*P_{2*}G_2F_2P_2^*F_1'\ar[d]^{P_2^*P_{2*}\omega_2P_2^*F_1'}\\
  P_2^*F_1'\ar[rrrrrr]^{P_2^*\eta_2F_1'}&&&&&&P_2^*P_{2*}P_2^*F_1'.
 }\] 
 This diagram is commutative by the commutativity of the diagram $(\diamondsuit)_2$.
 
 Now we see that the diagram $(\diamondsuit)$ is commutative as follows;
 \begin{align}
P_1^*\eta_1\circ\omega P_1^*&=T_1\omega_1P_1^*\circ\Omega_{G_1}F_1P_1^*\circ G_1\Omega_{F_1}P_1^*\circ G_1F_1P_1^*\eta_1\circ G_1\omega_2F_1P_1^*\notag\\
&=T_1\omega_1P_1^*\circ\Omega_{G_1}F_1P_1^*\circ G_1T_2\omega_2 F_1P_1^*\circ G_1\Omega_{G_2}FP_1^*\circ G\Omega_{F_2}F_1P_1^*\circ GF_2\Omega_{F_1}P_1^*\circ GFP_1^*\eta_1\notag\\
&=T_1\omega_1P_1^*\circ T_1G_1\omega_2F_1P_1^*\circ \Omega_{G_1}G_2FP_1^*\circ G_1\Omega_{G_2}FP_1^*\circ G\Omega_{F_2}F_1P_1^*\circ GF_2\Omega_{F_1}P_1^*\circ GFP_1^*\eta_1\notag\\
&=T_1\omega P_1^*\circ\Omega_{G}FP_1^*\circ G\Omega_FP_1^*\circ GFP_1^*\eta_1,\notag
\end{align}
where the first equation (resp. the second equation)  follows from the commutativity of the diagram $(\diamondsuit)_1$ (resp. $(\clubsuit)$), and the third equation  follows from  the functoriality of the functor isomorphism $\Omega_{G_1}$. 
 \end{proof}
\end{lem}

\vspace{2mm}

The adjoint pair $$ \Psi_{P_L}\dashv\Psi_Q$$
 \[\xymatrix{
{\rm Dcoh}(G\times X_1,p_1^*W_1)\ar@/^/[rr]^{\Psi_Q}&&{\rm Dcoh}(G\times X_2,p_2^*W_2)\ar@/^/[ll]^{\Psi_{Q_L}}
  }\] 
is induced by the following three adjoint pairs 

\vspace{2mm}
$(1): \overline{\pi_1}'_!\dashv{\overline{\pi_1}'}^*$
 \[\xymatrix{
 {\rm Dcoh}(G\times X_1,p_1^*W_1)\ar@/^/[rr]^{{\overline{\pi_1}'}^*}&&{\rm Dcoh}(G\times X_1\times\overline{X_2},{\overline{\pi_1}'}^*p_1^*W_1)\ar@/^/[ll]^{\overline{\pi_1}'_!},
  }\] 
  where $\overline{\pi_1}'_!:={\bf R}\overline{\pi_1}'_*((-)\otimes\overline{p_{12}}^*\omega_{\overline{\pi_1}}[{\rm dim}(\overline{X_2})])$.
 
 \vspace{2mm}
  $(2):i'_*((-)\otimes^{\bf L} Q^{\bf L\vee})\dashv {i'}^*(-)\otimes^{\bf L} Q$
   \[\xymatrix{
&{\rm Dcoh}(G\times X_1\times\overline{X_2},{\overline{\pi_1}'}^*p_1^*W_1)\ar@/^/[rr]^{{i'}^*(-)\otimes^{\bf L} Q}&&{\rm Dcoh}(G\times X_1 \times X_2,{\pi_2'}^*p_2^*W_2)\ar@/^/[ll]^{i'_*((-)\otimes^{\bf L} Q^{\bf L\vee})}
  }\] 
  and 

  \vspace{2mm}
    $(3):{\pi_2'}^*\dashv{\bf R}{\pi_2'}_*$
\[\xymatrix{
{\rm Dcoh}(G\times X_1 \times X_2,{\pi_2'}^*p_2^*W_2)\ar@/^/[rr]^{{\bf R}{\pi_2'}_*}&&{\rm Dcoh}(G\times X_2,p_2^*W_2)\ar@/^/[ll]^{{\pi_2'}^*}.
  }\]

 Hence the adjunction morphism $\omega':\Psi_{Q_L}\Psi_{Q}\rightarrow{\rm id}_{{\rm Dcoh}(G\times X_1,p_1^*W_1)}$ in the diagram $(\ast)'$ is the composition 

\begin{align}
\Psi_{Q_L}\Psi_{Q}&={\bf R}\overline{\pi}'_{1*}(i'_*({\pi_2'}^*{\bf R}\pi'_{2*}({i'}^*{\overline{\pi_1}'}^*(-)\otimes^{\bf L} Q)\otimes^{\bf L} Q^{{\bf L}\vee})\otimes\overline{p_{12}}^*\omega_{\overline{\pi_1}}[{\rm dim}(\overline{X_2})])\notag\\
&\xrightarrow{\zeta_3} {\bf R}\overline{\pi}'_{1*}(i'_*({i'}^*{\overline{\pi_1}'}^*(-)\otimes^{\bf L} Q\otimes^{\bf L} Q^{{\bf L}\vee})\otimes\overline{p_{12}}^*\omega_{\overline{\pi_1}}[{\rm dim}(\overline{X_2})])\notag\\
&\xrightarrow{\zeta_2}{\bf R}\overline{\pi}'_{1*}({\overline{\pi_1}'}^*(-)\otimes\overline{p_{12}}^*\omega_{\overline{\pi_1}}[{\rm dim}(\overline{X_2})])\notag\\
&\xrightarrow{\zeta_1}{\rm id}_{{\rm Dcoh}(G\times X_1,p_1^*W_1)},\notag
\end{align}
where for each $i=1,2,3$, $\zeta_i$ is the functor morphism induced by the adjunction morphism of the above adjunction pair $(i)$.
Hence, by Lemma \ref{decomposition} and Lemma \ref{linearizable2}, to prove that the diagram $(\ast)'$ is commutative, it is enough to prove that the following diagrams $(\ast)'_i$ are commutative;

\vspace{3mm}
$(\ast)'_1:$
  \[\xymatrix{
\overline{\pi_1}'_{!}{\overline{\pi_1}'}^*p_1^*\ar[rr]^{\overline{\pi_1}'_{!}{\overline{\pi_1}'}^*p_1^*\eta_{p_1}}\ar[d]_{\omega'_1 p_1^*}&&\overline{\pi_1}'_{!}{\overline{\pi_1}'}^*p_1^*p_{1*}p_1^* \ar[r]^{\overline{\pi_1}'_{!}\Omega_1 p_1^*}&\overline{\pi_1}'_{!}\overline{p_{12}}^*\overline{p_{12}}_{*}{\overline{\pi_1}'}^*p_1^*\ar[r]^{\Omega_{L1}{\overline{\pi_1}'}^* p_1^*} &p_1^*p_{1*}\overline{\pi_1}'_{!}{\overline{\pi_1}'}^* p_1^* \ar[d]^{p_1^*p_{1*}\omega'_1 p_1^*}\\
  p_1^* \ar[rrrr]^{p_1^*\eta_{p_1}} &&&& p_1^*p_{1*}p_1^*, 
 }\]
where $\Omega_1:{\overline{\pi_1}'}^*p_1^*p_{1*}\xrightarrow{\sim}\overline{p_{12}}^*\overline{p_{12}}_{*}{\overline{\pi_1}'}^*$ and $\Omega_{L1}:\overline{\pi_1}'_{!}\overline{p_{12}}^*\overline{p_{12}}_{*}\xrightarrow{\sim}p_1^*p_{1*}\overline{\pi_1}'_{!}$ are the functor isomorphisms given by Lemma \ref{linearizable2}, and $\omega'_1:\overline{\pi_1}'_{!}{\overline{\pi_1}'}^*\rightarrow{\rm id}$ is the adjunction morphism of the adjoint pair (1).

\vspace{5mm}
$(\ast)'_2:$
  \[\xymatrix{
i'_{Q*}{i'_Q}^*\overline{p_{12}}^*\ar[rrr]^{i'_{Q*}{i'_Q}^*\overline{p_{12}}^*\eta_{\overline{p_{12}}}}\ar[d]_{\omega'_2 \overline{p_{12}}^*}&&&i'_{Q*}{i'_Q}^*\overline{p_{12}}^*\overline{p_{12}}_{*}\overline{p_{12}}^*\ar[r]^{i'_{Q*}\Omega_2\overline{p_{12}}^*}&i'_{Q*}{p_{12}}^*{p_{12}}_{*}{i'_Q}^*\overline{p_{12}}^*\ar[r]^{\Omega_{L2}{i'_Q}^*\overline{p_{12}}^*} &\overline{p_{12}}^*\overline{p_{12}}_{*} i'_{Q*}{i'_Q}^*\overline{p_{12}}^*\ar[d]^{ \overline{p_{12}}^*\overline{p_{12}}_{*}\omega'_2 \overline{p_{12}}^*}\\
\overline{p_{12}}^* \ar[rrrrr]^{\overline{p_{12}}^*\eta_{\overline{p_{12}}}} &&&&& \overline{p_{12}}^*\overline{p_{12}}_{*}\overline{p_{12}}^*, 
 }\]
where $i'_{Q*}(-):=i'_*((-)\otimes^{\bf L} Q^{\bf L\vee})$ and ${i'_Q}^*:={i'}^*(-)\otimes^{\bf L} Q$, and $\Omega_2:{i'_Q}^*\overline{p_{12}}^*\overline{p_{12}}_{*}\xrightarrow{\sim}{p_{12}}^*{p_{12}}_{*}{i'_Q}^*$ and $\Omega_{L2}:i'_{Q*}{p_{12}}^*{p_{12}}_{*}\xrightarrow{\sim}{p_{12}}^*\overline{p_{12}}_{*} i'_{Q*}$ are the functor isomorphisms given by Lemma \ref{linearizable2}, and $\omega'_2:i'_{Q*}{i'_Q}^*\rightarrow{\rm id}$ is the adjunction morphism of the adjoint pair (2).

\vspace{5mm}
$(\ast)'_3:$
  \[\xymatrix{
{\pi'_2}^*{\bf R}{\pi'_2}_*{p_{12}}^*\ar[rr]^{{\pi'_2}^*{\bf R}{\pi'_2}_*{p_{12}}^*\eta_{{p_{12}}}}\ar[d]_{\omega'_3 {p_{12}}^*}&&{\pi'_2}^*{\bf R}{\pi'_2}_*{p_{12}}^*{p_{12}}_{*}{p_{12}}^*\ar[r]^{{\pi'_2}^*\Omega_3{p_{12}}^*}&{\pi'_2}^*{p_{2}}^*{p_{2}}_{*}{\bf R}{\pi'_2}_*{p_{12}}^*\ar[r]^{\Omega_{L3}{\bf R}{\pi'_2}_*{p_{12}}^*} &{p_{12}}^*{p_{12}}_{*} {\pi'_2}^*{\bf R}{\pi'_2}_*{p_{12}}^*\ar[d]^{ {p_{12}}^*{p_{12}}_{*}\omega'_3{p_{12}}^*}\\
{p_{12}}^* \ar[rrrr]^{{p_{12}}^*\eta_{{p_{12}}}} &&&&{p_{12}}^*{p_{12}}_{*}{p_{12}}^*, 
 }\]
where $\Omega_3:{\bf R}{\pi'_2}_*{p_{12}}^*{p_{12}}_{*}\xrightarrow{\sim}{p_{2}}^*{p_{2}}_{*}{\bf R}{\pi'_2}_*$ and $\Omega_{L3}:{\pi'_2}^*{p_{2}}^*{p_{2}}_{*}\xrightarrow{\sim}{p_{12}}^*{p_{12}}_{*} {\pi'_2}^*$ are the functor isomorphisms given by Lemma \ref{linearizable2}, and $\omega'_3:{\pi'_2}^*{\bf R}{\pi'_2}_*\rightarrow{\rm id}$ is the adjunction morphism of the adjoint pair (3).

\vspace{5mm}
In the following, for each $i=1,2,3$, we will prove that the diagram $(\ast)'_i$ is commutative.

\vspace{3mm}
\begin{flushleft}
\textbf{$\bullet$ Proof of the commutativity of $(\ast)'_1$}
\end{flushleft}

\vspace{2mm}
Since the adjunction morphism $\omega'_1:\overline{\pi_1}'_{!}{\overline{\pi_1}'}^*\rightarrow{\rm id}_{{\rm Dcoh}(G\times X_1,p_1^*W_1)}$ is a restriction of the adjunction morphism $\omega'_1:\overline{\pi_1}'_{!}{\overline{\pi_1}'}^*\rightarrow{\rm id}_{{\rm DQcoh}(G\times X_1,p_1^*W_1)}$ of the adjoint pair
\[\xymatrix{
 {\rm DQcoh}(G\times X_1,p_1^*W_1)\ar@/^/[rr]^{{\overline{\pi_1}'}^*}&&{\rm DQcoh}(G\times X_1\times\overline{X_2},{\overline{\pi_1}'}^*p_1^*W_1)\ar@/^/[ll]^{\overline{\pi_1}'_!},
  }\] 
 we have the functor morphism $$\omega'_1p_1^*p_{1*}p_1^*:\overline{\pi_1}'_{!}{\overline{\pi_1}'}^*p_1^*p_{1*}p_1^*\rightarrow p_1^*p_{1*}p_1^*.$$ By the functoriality of $\omega'_1$, to prove the commutativity of $(\ast)'_1$ it is enough to prove that  the following diagram is commutative;
 \[\xymatrix{
(\ast)'_{1a}:&\overline{\pi_1}'_{!}{\overline{\pi_1}'}^*p_1^*p_{1*} \ar[rr]^{\overline{\pi_1}'_{!}\Omega_1 } \ar[drr]_{\omega'_1p_1^*p_{1*}} &&\overline{\pi_1}'_{!}\overline{p_{12}}^*\overline{p_{12}}_{*}{\overline{\pi_1}'}^*\ar[rr]^{\Omega_{L1}{\overline{\pi_1}'}^* } &&p_1^*p_{1*}\overline{\pi_1}'_{!}{\overline{\pi_1}'}^*  \ar[dll]^{p_1^*p_{1*}\omega'_1 }\\
&&& p_1^*p_{1*}&&
 }\]
The adjunction morphism $$\omega'_1:\overline{\pi_1}'_{!}{\overline{\pi_1}'}^*\rightarrow{\rm id}$$ is given by the composition of the following functor morphisms;

$$\varphi:\overline{\pi_1}'_{!}{\overline{\pi_1}'}^*(-)={\bf R}\overline{\pi_1}'_*({\overline{\pi_1}'}^*(-)\otimes\overline{p_{12}}^*\omega_{\overline{\pi_1}}[d_2])\longrightarrow(-)\otimes^{{\bf L}} p_1^*{\bf R}\overline{\pi_1}_*\omega_{\overline{\pi_1}}[d_2]$$
and
$$\psi:(-)\otimes ^{{\bf L}} p_1^*{\bf R}\overline{\pi_1}_*\omega_{\overline{\pi_1}}[d_2]\longrightarrow(-),$$
 where $d_2:={\rm dim}(\overline{X_2})$, the functor morphism $\varphi$ is given by the projection formula and an isomorphism ${\bf R}\overline{\pi_1}'_*\overline{p_{12}}^*\cong p_1^*{\bf R}\overline{\pi_1}_*$, whence $\varphi$ is a functor {\it isomorphism}, and $\psi$ is given as follows.
  Let 
 $$\sigma:{\bf R}\overline{\pi_1}_*\omega_{\overline{\pi_1}}[d_2]\longrightarrow \mathcal{O}_{ X_1}$$
 be the following composition of morphisms  in ${\rm D^b}(X_1)$;
 $${\bf R}\overline{\pi_1}_*\omega_{\overline{\pi_1}}[d_2]\xrightarrow{\sim}{\bf R}\overline{\pi_1}_*{\overline{\pi_1}}^!(\mathcal{O}_{X_1})\longrightarrow \mathcal{O}_{ X_1},$$
  where the morphism ${\bf R}\overline{\pi_1}_*{\overline{\pi_1}}^!(\mathcal{O}_{ X_1})\rightarrow \mathcal{O}_{ X_1}$ is  induced by the adjunction morphism of the adjoint pair,
\[\xymatrix{
 {\bf R}\overline{\pi_1}_*\dashv{\overline{\pi_1}}^!&{\rm D^b}({\rm coh} X_1)\ar@/^/[rr]^{{\overline{\pi_1}}^!}&&{\rm D^b}({\rm coh} X_1\times\overline{X_2})\ar@/^/[ll]^{{\bf R}\overline{\pi_1}_*}.
  }\] 
 Then the  functor morphism $\psi$ is given as $$\psi:=(-)\otimes p_1^*\Upsilon(\sigma),$$
where $\Upsilon:{\rm D^b}( X_1)\rightarrow{\rm Dcoh}( X_1,0)$ is the functor defined in Definition \ref{upsilon}. 
Hence it is enough to prove that for any object $F\in{\rm Dcoh}(G\times X_1,p_1^*W)$ the following two diagrams are commutative, 
\[\xymatrix{
(\ast)'_{1b}: &{\bf R}\overline{\pi_1}'_*({\overline{\pi_1}'}^*p_1^*p_{1*}(F)\otimes\overline{p_{12}}^*\omega_{\overline{\pi_1}}[d_2])\ar[d]_{\overline{\pi_1}'_{!}\Omega_1}\ar[rr]^{\varphi p_1^*p_{1*}}&&p_1^*p_{1*}(F)\otimes^{{\bf L}} p_1^*{\bf R}\overline{\pi_1}_*\omega_{\overline{\pi_1}}[d_2]\ar[d]\\
&{\bf R}\overline{\pi_1}'_*(\overline{p_{12}}^*\overline{p_{12}}_*{\overline{\pi_1}'}^*(F)\otimes\overline{p_{12}}^*\omega_{\overline{\pi_1}}[d_2])\ar[d]_{\Omega_{L1}{\overline{\pi_1}'}^*}&&p_1^*(p_{1*}(F)\otimes^{{\bf L}} {\bf R}\overline{\pi_1}_*\omega_{\overline{\pi_1}}[d_2])\ar[d]\\
&p_1^*p_{1*}{\bf R}\overline{\pi_1}'_*({\overline{\pi_1}'}^*(F)\otimes\overline{p_{12}}^*\omega_{\overline{\pi_1}}[d_2])\ar[rr]^{p_1^*p_{1*}\varphi}&&p_1^*p_{1*}((F)\otimes^{{\bf L}} p_1^*{\bf R}\overline{\pi_1}_*\omega_{\overline{\pi_1}}[d_2])
  }\]
  and
\[\xymatrix{
(\ast)'_{1c}: &p_1^*p_{1*}(F)\otimes^{{\bf L}} p_1^*{\bf R}\overline{\pi_1}_*\omega_{\overline{\pi_1}}[d_2]\ar[rrr]^{p_1^*p_{1*}(F)\otimes p_1^*\Upsilon(\sigma)}\ar[d]&&&p_1^*p_{1*}(F)\otimes p_1^*\mathcal{O}_{X_1}\ar[d]\ar[rrd]&&\\
&p_1^*(p_{1*}(F)\otimes^{{\bf L}} {\bf R}\overline{\pi_1}_*\omega_{\overline{\pi_1}}[d_2])\ar[d]\ar[rrr]^{p_1^*(p_{1*}(F)\otimes \Upsilon(\sigma))}&&&p_1^*(p_{1*}(F)\otimes \mathcal{O}_{X_1})\ar[d]&&p_1^*p_{1*}(F)\\
&p_1^*p_{1*}((F)\otimes^{{\bf L}} p_1^*{\bf R}\overline{\pi_1}_*\omega_{\overline{\pi_1}}[d_2])\ar[rrr]^{p_1^*p_{1*}((F)\otimes p_1^*\Upsilon(\sigma))}&&&p_1^*p_{1*}((F)\otimes p_1^*\mathcal{O}_{X_1}),\ar[rru]&&
  }\]  
where  arrows with no symbols are natural isomorphisms.

At first, we show the diagram $(\ast)'_{1b}$ is commutative. Since functor morphisms in the diagram $(\ast)'_{1b}$ are natural in $F$ and $\omega_{\overline{\pi_1}}[d_2]$, we can replace the objects $F$ and $\omega_{\overline{\pi_1}}[d_2]$ with  objects $E\in{\rm Dlfr}(G\times X_1,p_1^*W)$ and  $I\in {\rm DQcoh}(X_1\times \overline{X_2},\overline{\pi_1}^*W)$ whose components $I_1$ and $I_0$ are injective sheaves  respectively.   Then derived functors in $(\ast)'_{1b}$ are isomorphic to underived functors, since  the derived functor in the lowest row on the right side in $(\ast)'_{1b}$ is isomorphic to underived functor, and the direct images $p_{1*}$ and $\overline{p_{12}}_*$ maps locally free sheaves to locally free sheaves, and the projection formulae  for $p_1$ and $\overline{p_{12}}$ hold in  categories of quasi-coherent sheaves  without assuming locally freeness of sheaves. So it is enough to prove that the commutativity of the similar diagram in the abelian category ${\rm Qcoh}(G\times X_1)$. But this is checked by easy computations.
 
Next, we show the diagram $(\ast)'_{1c}$ is commutative. The commutativity of two square diagrams on the left side  follows  automatically by the functoriality. So we have only to verify that the triangular diagram on the right side is commutative. But this is verified by easy computations, and the detail is left to the reader.

\vspace{3mm}
\begin{flushleft}
\textbf{$\bullet$ Proof of the commutativity of $(\ast)'_2$}
\end{flushleft}
To decompose the diagram $(\ast)'_2$, we will embed the diagram $(\ast)'_2$ to a larger category.
Before embedding it, we provide some functors and some functor morphisms.

\vspace{2mm}
Since the functor $i'_{\#}:{\rm Qcoh}(G\times X_1\times X_2)\rightarrow {\rm Ind}({\rm Pro}({\rm Qcoh}(G\times X_1\times \overline{X_2})))$, constructed in Definition \ref{extended extension by zero},  is exact and compatible with potentials, it induces a functor 
$$i'_{\#}:{\rm DQcoh}(G\times X_1\times X_2,{\pi'_1}^*p_1^*W_1)\rightarrow {\rm Ind}({\rm Pro}({\rm DQcoh}(G\times X_1\times \overline{X_2},\overline{\pi'_1}^*p_1^*W_1))).$$
Let $i'_!:{\rm Dcoh}(G\times X_1\times X_2,{\pi'_1}^*p_1^*W_1)\rightarrow {\rm Pro}({\rm DQcoh}(G\times X_1\times \overline{X_2},\overline{\pi'_1}^*p_1^*W_1))$ be the extension by zero, and set 
$$i'_{Q!}(-):=i'_!((-)\otimes^{\bf L} Q^{\bf L{\vee}})\hspace{5mm}{\rm and}\hspace{5mm}i'_{Q\#}(-):=i'_{\#}((-)\otimes^{\bf L} Q^{\bf L{\vee}}).$$
Functor morphisms constructed in Lemma \ref{gammas} (1) induces functor isomorphism 
$$\gamma_{Q!}:i'_{Q!}\xrightarrow{\sim}\iota_{\rm Pro}i'_{Q*}$$
and functor morphism
$$\gamma_{Q\#}:i'_{Q\#}\rightarrow\iota_{\rm InPr}i'_{Q*}.$$
with $\iota_{\rm Ind}\gamma_{Q!}=\gamma_{Q\#}|_{{\rm Dcoh}(G\times X_1\times X_2,{\pi'_1}^*p_1^*W_1)}$.  Let  $\omega'_2: i'_{Q*}i^{'*}_Q\rightarrow {\rm id}$ be the adjunction morphism. Then the morphism $\iota_{\rm Pro}\omega'_2: \iota_{\rm Pro}i'_{Q*}i'^{*}_Q\rightarrow \iota_{\rm Pro}$ is decomposed into the following compositions
$$\iota_{\rm Pro}i'_{Q*}i'^{*}_Q\xrightarrow{\gamma_{Q!}^{-1}}i'_{Q!}i'^{*}_Q\xrightarrow{i'_!i^{'*}\omega_Q}i'_!i'^{*}\xrightarrow{\omega_{i_!}}\iota_{\rm Pro},$$
where $\omega_Q:(-)\otimes^{\bf L}Q\otimes^{\bf L}Q^{\bf L{\vee}}\rightarrow (-)$ and $\omega_{i'_!}:i'_!i^{'*}\rightarrow\iota_{\rm Pro}$ are the adjunction morphisms.
Furthermore, the functor morphism constructed in  Lemma \ref{gammas} (2) induces   a functor morphism 
$$\lambda: i'_{\#}p^*_{12}p_{12*}\rightarrow \iota_{\rm Ind}{\rm Pro}(\overline{p_{12}}^*\overline{p_{12}}_*)i'_!.$$

\vspace{4mm}
Now we are ready to decompose the diagram $(\ast)'_2$.  Let $\Omega_{i'*}:i'_*p_{12}^*p_{12*}\xrightarrow{\sim}\overline{p_{12}}^*\overline{p_{12}}_{*}i'_*$ and $\Omega_{i'}^*:i'^*\overline{p_{12}}^*\overline{p_{12}}_{*}\xrightarrow{\sim}p_{12}^*p_{12*}i'^*$ be natural functor isomorphisms. Set $i'^*_{Q\otimes Q^{\vee}}(-):=i'^*(-)\otimes^{\bf L}Q\otimes^{\bf L}Q^{\bf L{\vee}}$, and let  $\Omega^{i'*}_{Q\otimes Q^{\vee}}:i'^*_{Q\otimes Q^{\vee}}\overline{p_{12}}^*\overline{p_{12}}_{*}\xrightarrow{\sim}{p_{12}}^*{p_{12}}_{*}i'^*_{Q\otimes Q^{\vee}}$ be the functor isomorphism given by natural functor isomorphims $\Omega_2:{i'^*_Q}\overline{p_{12}}^*\overline{p_{12}}_{*}\xrightarrow{\sim}p_{12}^*{p_{12}}_{*}{i'_Q}^*$ and $\Omega_{Q^{\vee}}:p_{12}^*p_{12*}(-)\otimes^{\bf L}Q^{\bf L{\vee}}\xrightarrow{\sim}p_{12}^*p_{12*}((-)\otimes^{\bf L}Q^{\bf L{\vee}})$.
Embedding the diagram $(\ast)'_2$ into the category ${\rm Ind}({\rm Pro}({\rm DQcoh}(G\times X_1\times \overline{X_2},\overline{\pi'_1}^*p_1^*W_1)))$ by the inclusion $$\iota_{\rm InPr}:{\rm DQcoh}(G\times X_1\times \overline{X_2},\overline{\pi'_1}^*p_1^*W_1)\hookrightarrow {\rm Ind}({\rm Pro}({\rm DQcoh}(G\times X_1\times \overline{X_2},\overline{\pi'_1}^*p_1^*W_1))),$$ the diagram $(\ast)'_2$ is decomposed into the following diagram

  \[\xymatrix{
  &&&\ar@{}[rrd]|{(a)}&i'_{Q*}{p_{12}}^*{p_{12}}_{*}{i'_Q}^*\overline{p_{12}}^*\ar[rd]^{\Omega_{L2}}&\\
\ar@{}[rrrd]|{(b)}i'_{Q*}{i'^*_Q}\overline{p_{12}}^*\ar[rrr]^{\eta_{\overline{p_{12}}}}\ar[d]_{\gamma_{Q!}^{-1}}&&&\ar@{}[rd]|{(c)}i'_{Q*}{i'^*_Q}\overline{p_{12}}^*\overline{p_{12}}_{*}\overline{p_{12}}^*\ar[r]^{\Omega^{i'*}_{Q\otimes Q^{\vee}}}\ar[ru]^{\Omega_2}&\ar@{}[rd]|{(d)}i'_{*}{p_{12}}^*{p_{12}}_{*}{i'^*_{Q\otimes Q^{\vee}}}\overline{p_{12}}^*\ar[r]^{\Omega_{i'*}} &\overline{p_{12}}^*\overline{p_{12}}_{*} i'_{Q*}{i'^*_Q}\overline{p_{12}}^*\ar[d]^{\gamma_{Q!}^{-1}}\\
\ar@{}[rrrd]|{(e)}i'_{Q!}{i'^*_Q}\overline{p_{12}}^*\ar[rrr]^{\eta_{\overline{p_{12}}}}\ar[d]_{\omega_Q}&&&i'_{Q\#}{i'^*_Q}\overline{p_{12}}^*\overline{p_{12}}_{*}\overline{p_{12}}^*\ar[r]^{\Omega^{i'*}_{Q\otimes Q^{\vee}}}\ar[u]^{\gamma_{Q\#}}\ar[d]_{\omega_Q}&\ar@{}[rd]|{(f)}i'_{\#}{p_{12}}^*{p_{12}}_{*}{i'^*_{Q\otimes Q^{\vee}}}\overline{p_{12}}^*\ar[r]^{\lambda}\ar[u]_{\gamma_{Q\#}}\ar[d]^{\omega_Q} &\overline{p_{12}}^*\overline{p_{12}}_{*} i'_{Q!}{i'^*_Q}\overline{p_{12}}^*\ar[d]^{\omega_Q}\\
i'_!i'_*\overline{p_{12}}^*\ar[rrr]^{\eta_{\overline{p_{12}}}}\ar[d]_{\omega_{i'_!}}&&&i'_{\#}i'^*\overline{p_{12}}^*\overline{p_{12}}_{*}\overline{p_{12}}^*\ar[r]^{\Omega_{i'}^*}& i'_{\#}{p_{12}}^*{p_{12}}_{*}i'^*\overline{p_{12}}^*\ar[r]^{\lambda}&\overline{p_{12}}^*\overline{p_{12}}_{*}i'_!i'^*\overline{p_{12}}^*\ar[d]^{\omega_{i'_!}}\\
\overline{p_{12}}^* \ar[rrrrr]^{\eta_{\overline{p_{12}}}} &&&&& \overline{p_{12}}^*\overline{p_{12}}_{*}\overline{p_{12}}^*, 
 }\]
where functor morphisms attached to  arrows are the ones which induce the functor morphisms, and we omit  embedding functors $\iota_{\rm InPr}$ and $\iota_{\rm Pro}$ from the above diagram. The diagram $(a)$ is commutative, since $\Omega_{L2}$ is given by $\Omega_{Q^{\vee}}$ and $\Omega_{i'*}$. The commutativity of the diagrams $(b)$, $(c)$, $(e)$ and $(f)$ follows from the functoriality of functor morphisms, and the diagram $(d)$ is commutative by Lemma \ref{gammas} (2). Hence, it is enough to verify the commutativity of the following diagrams
 \[\xymatrix{
 (\ast)'_{2a}:&i'^*\overline{p_{12}}^*\overline{p_{12}}_{*}(-)\otimes^{\bf L}Q\otimes^{\bf L}Q^{\bf L{\vee}}\ar[rr]^{\Omega^{i'*}_{Q\otimes Q^{\vee}}}\ar[d]_{\omega_Q i'^*\overline{p_{12}}^*\overline{p_{12}}_{*}}&&p_{12}^*p_{12*}(i'^*(-)\otimes^{\bf L}Q\otimes^{\bf L}Q^{\bf L{\vee}})\ar[d]^{p_{12}^*p_{12*}\omega_Q i'^*}\\
& i'\overline{p_{12}}^*\overline{p_{12}}_{*}(-)\ar[rr]^{\Omega_{i'}^*}&&p_{12}^*p_{12*}i'^*
  }\]
and 
 \[\xymatrix{
 (\ast)'_{2b}:&\iota_{\rm Ind}i'_!i'_*\ar[rrr]^{i'_{\#}i'^*\eta_{\overline{p_{12}}}}\ar[d]_{\iota_{\rm Ind}\omega_{i'_!}}&&&i'_{\#}i'^*\overline{p_{12}}^*\overline{p_{12}}_{*}\ar[r]^{i'_{\#}\Omega_{i'}^*}& i'_{\#}{p_{12}}^*{p_{12}}_{*}i'^*\ar[r]^-{\lambda i'^*}&\iota_{\rm Ind}{\rm Pro}(\overline{p_{12}}^*\overline{p_{12}}_{*})i'_!i'^*\ar[d]^{\iota_{\rm Ind}{\rm Pro}(\overline{p_{12}}^*\overline{p_{12}}_{*})\omega_{i'_!}}\\
&\iota_{\rm InPr} \ar[rrrrr]^{\iota_{\rm InPr}\eta_{\overline{p_{12}}}} &&&&& \iota_{\rm InPr}\overline{p_{12}}^*\overline{p_{12}}_{*}
  }\]
  
  We show that the diagram $(\ast)'_{2a}$ is commutative. Let $\Omega_Q:p_{12}^*p_{12*}(-)\otimes^{\bf L}Q\xrightarrow{\sim}p_{12}^*p_{12*}((-)\otimes^{\bf L}Q)$ be the natural functor isomorphism. Then, the  functor morphism $\Omega_2:{i'_Q}^*\overline{p_{12}}^*\overline{p_{12}}_{*}\xrightarrow{\sim}{p_{12}}^*{p_{12}}_{*}{i'_Q}^*$ is the  following compositions of functor morphisms
  $$i'^*\overline{p_{12}}^*\overline{p_{12}}_{*}(-)\otimes^{\bf L}Q\xrightarrow{((-)\otimes^{\bf L}Q)\Omega_{i'}^*}p_{12}^*p_{12*}i'^*(-)\otimes^{\bf L}Q\xrightarrow{\Omega_Qi'^*}p_{12}^*p_{12*}(i'^*(-)\otimes^{\bf L}Q).$$
 Moreover, the following diagram
  \[\xymatrix{
i'^*\overline{p_{12}}^*\overline{p_{12}}_{*}(-)\otimes^{\bf L}Q\otimes^{\bf L}Q^{\bf L{\vee}}\ar[rr]^{((-)\otimes^{\bf L}Q\otimes^{\bf L}Q^{\bf L{\vee}})\Omega_{i'}^*}\ar[d]_{\omega_Q i'^*\overline{p_{12}}^*\overline{p_{12}}_{*}}&&p_{12}^*p_{12*}i'^*(-)\otimes^{\bf L}Q\otimes^{\bf L}Q^{\bf L{\vee}}\ar[d]^{\omega_Q p_{12}^*p_{12*}i'^*}\\
 i'\overline{p_{12}}^*\overline{p_{12}}_{*}(-)\ar[rr]^{\Omega_{i'}^*}&&p_{12}^*p_{12*}i'^*
  }\]
is commutative by the functoriality of the functor morphism $\omega_Q$. Hence, to show that the diagram $(\ast)'_{2a}$ is commutative, we have only to show the commutativity of the following diagram

 \[\xymatrix{
 p_{12}^*p_{12*}i'^*(-)\otimes^{\bf L}Q\otimes^{\bf L}Q^{\bf L{\vee}}\ar[rrd]_{\omega_Q p_{12}^*p_{12*}i'^*}\ar[rr]^{((-)\otimes^{\bf L}Q^{\bf L{\vee}})\Omega_Qi'^*} && p_{12}^*p_{12*}(i'^*_Q(-))\otimes^{\bf L}Q^{\bf L{\vee}}\ar[rr]^{\Omega_{Q^{\vee}}(i'^*(-)\otimes^{\bf L}Q)} &&p_{12}^*p_{12*}(i'^*_Q(-)\otimes^{\bf L}Q^{\bf L{\vee}})\ar[lld]^{p_{12}^*p_{12*}\omega_Q i'^*}\\
  &&p_{12}^*p_{12*}i'^*&&
  }\]
Replacing the object $P\in{\rm Dcoh}(X_1\times X_2, \pi_2^*W_2-\pi_1^*W_1)$ with an object in ${\rm Dlfr}(X_1\times X_2, \pi_2^*W_2-\pi_1^*W_1)$, we may assume that the object $Q=p_{12}^*P$ is an object whose components are locally free sheaves. Then, the functors in the above diagram are underived functors. Hence, the commutativity of the diagram is verified by easy diagram chasing of morphisms between quasi-coherent sheaves, which is left to the reader.

Since all of the functors in  $(\ast)'_{2b}$ are underived functors, the diagram $(\ast)'_{2b}$ is also verified by diagram chasing of map of systems, which is also left to the reader. 


\vspace{3mm}
\begin{flushleft}
\textbf{$\bullet$ Proof of the commutativity of $(\ast)'_3$}
\end{flushleft}

By the functoriality of $\omega'_3$, the following diagram is commutative:
 \[\xymatrix{
{\pi'_2}^*{\bf R}{\pi'_2}_*{p_{12}}^*\ar[rrr]^{{\pi'_2}^*{\bf R}{\pi'_2}_*{p_{12}}^*\eta_{{p_{12}}}}\ar[d]_{\omega'_3 {p_{12}}^*}&&&{\pi'_2}^*{\bf R}{\pi'_2}_*{p_{12}}^*{p_{12}}_{*}{p_{12}}^*\ar[d]^{\omega'_3p_{12}^*p_{12*}p_{12}^*}\\
{p_{12}}^* \ar[rrr]^{{p_{12}}^*\eta_{{p_{12}}}} &&&{p_{12}}^*{p_{12}}_{*}{p_{12}}^*
 }\]
 Hence, to prove that the diagram $(\ast)'_3$ is commutative, it is enough to prove the following diagram is commutative:
 
  \[\xymatrix{
{\pi'_2}^*{\bf R}{\pi'_2}_*{p_{12}}^*{p_{12}}_{*}\ar[rrd]_{\omega'_3 p_{12}^*p_{12*}}\ar[rr]^{{\pi'_2}^*\Omega_3} && {\pi'_2}^*p_2^*p_{2*}{\bf R}\pi'_{2*}\ar[rr]^{\Omega_{L3}{\bf R}\pi'_{2*}} &&p_{12}^*p_{12*}{\pi'_2}^*{\bf R}\pi'_{2*}\ar[lld]^{p_{12}^*p_{12*}\omega'_3}\\
  &&p_{12}^*p_{12*}&&
  }\]
 Since we may replace any object in ${\rm Dcoh}(G\times X_1\times X_2,{\pi'_2}^*p_2^*W_2)$ with an object whose components are injective sheaves, the commutativity of the above diagram can be checked by easy diagram chasing of morphisms between quasi-coherent sheaves, which is left to the reader.

\vspace{5mm}
\subsection{Main Theorem}

At first, to state the main theorem, we give the definition of $G$-linearizable objects.

\begin{dfn}
Let $X$ be a variety with $G$-action. An object $F$ of ${\rm D^b}({\rm coh}X)$ is called $G$\textbf{-linearizable}, if $F$ is in the essential image of the forgetful functor ${\rm D^b}({\rm coh}_GX)\rightarrow{\rm D^b}({\rm coh}X)$.

\end{dfn}

\vspace{2mm}
Now we are ready to state and prove the main theorem.
\begin{thm}\label{main theorem 2}

Let $P\in{\rm D^b}({\rm coh}X_1\times_{\mathbb{A}^1} X_2)$ be a $G$-linearizable object whose support is proper over $X_1$ and $X_2$.  If the integral functor  $\Phi_{j_*(P)}:{\rm D^b}({\rm coh}X_1)\rightarrow{\rm D^b}({\rm coh}X_2)$ is an equivalence (resp. fully faithful), then there is an integral functor  $$\Phi_{\widetilde{P}_G}:{\rm Dcoh}_G(X_1,W_1)\rightarrow{\rm Dcoh}_G(X_2,W_2)$$ which is also an equivalence (resp. fully faithful) for some ${\widetilde{P}_G}\in {\rm Dcoh}_G(X_1\times X_2,\pi_2^*W_2-\pi_1^*W_1)$.

\begin{proof}Since $P$ is $G$-linearizable, we may assume that there is  an object $P_G\in{\rm D^b}({\rm coh}_GX_1\times_{\mathbb{A}^1} X_2)$ such that $\Pi(P_G)=P$, where $\Pi:{\rm D^b}({\rm coh}_GX_1\times_{\mathbb{A}^1} X_2)\rightarrow{\rm D^b}({\rm coh}X_1\times_{\mathbb{A}^1} X_2)$ is the forgetful functor. Set $$\widetilde{P}_G:=j_*(\Upsilon(P_G))\in {\rm Dcoh}_G(X_1\times X_2,\pi_2^*W_2-\pi_1^*W_1),$$
where $\Upsilon:{\rm D^b}({\rm coh}_GX_1\times_{\mathbb{A}^1} X_2)\rightarrow{\rm Dcoh}_G(X_1\times_{\mathbb{A}^1} X_2,0)$ is the exact functor defined in Definition \ref{upsilon}, and $j_*:{\rm Dcoh}_G(X_1\times_{\mathbb{A}^1} X_2,0)\rightarrow{\rm Dcoh}_G(X_1\times X_2,\pi_2^*W_2-\pi_1^*W_1)$ is the direct image of embedding $j:X_1\times_{\mathbb{A}^1} X_2\rightarrow X_1\times X_2$.
Let  $\widetilde{P}:={\rm Res}_G(\widetilde{P}_G)\in{\rm Dcoh}(X_1\times X_2,\pi_2^*W_2-\pi_1^*W_1)$. Then we have $$\widetilde{P}=j_*(\Upsilon(P))=j_*({\rm Tot}(\tau(P)))\cong{\rm Tot}(j_*(\tau(P))),$$ where $\tau:{\rm D^b}({\rm coh}X_1\times_{\mathbb{A}^1} X_2)\rightarrow{\rm D^b}(Z^0({\rm coh}(X_1\times_{\mathbb{A}^1} X_2,0)))$ is the functor given by the same manner as in just before Definition \ref{upsilon}, and $j_*$ in the last one is  the direct image  $$j_*:{\rm D^b}(Z^0({\rm coh}(X_1\times_{\mathbb{A}^1} X_2,0)))\rightarrow{\rm D^b}(Z^0({\rm coh}(X_1\times X_2,0)))$$ induced by an exact functor $j_*:Z^0({\rm coh}(X_1\times_{\mathbb{A}^1} X_2,0))\rightarrow Z^0({\rm coh}(X_1\times X_2,0))$ between abelian categories. Since  ${\rm Supp}(j_*(\tau(P)))={\rm Supp}(P)$,  $\widetilde{P}$ has a $\pi_i$-proper support ($i=1,2$). By Theorem \ref{equivariantization}, it is enough to show that if the integral functor  $\Phi_{j_*(P)}:{\rm D^b}({\rm coh}X_1)\rightarrow{\rm D^b}({\rm coh}X_2)$ is an equivalence (resp. fully faithful), then the integral functor $$\Phi_{\widetilde{P}}:{\rm Dcoh}(X_1,W_1)\rightarrow{\rm Dcoh}(X_2,W_2)$$ is an equivalence (resp. fully faithful).

Assume that the integral functor $\Phi_{j_*(P)}:{\rm D^b}({\rm coh}X_1)\rightarrow{\rm D^b}({\rm coh}X_2)$ is fully faithful.
The integral functor $\Phi_{j_*(P)}$ induces the extended functor $\Phi_{j_*(P)}':{\rm D}({\rm Qcoh}X_1)\rightarrow{\rm D}({\rm Qcoh}X_2)$. Then the functor $\Phi_{j_*(P)}'$ is also fully faithful since it preserves any direct limit and any object in the unbounded derived category ${\rm D}({\rm Qcoh}X_1)$ is isomorphic to the direct limit of a direct system of  objects in ${\rm D^b}({\rm coh}X_1)$ by \cite[Proposition 2.3.2]{tt}. Hence by the argument in the proof of  \cite[Theorem 5.15]{bdfik}, we obtain an isomorphism of functors $$\Phi'_{\widetilde{P}_R}\circ\Phi'_{\widetilde{P}}\cong {\rm id}_{{\rm DQcoh}(X_1,W_1)},$$
where $\Phi_{\tilde{P}}':{\rm DQcoh}(X_1,W_1)\rightarrow{\rm DQcoh}(X_2,W_2)$ and $\Phi'_{\widetilde{P}_R}:{\rm DQcoh}(X_2,W_2)\rightarrow{\rm DQcoh}(X_1,W_1)$ are the extended functors from $\Phi_{\widetilde{P}}$ and its right adjoint $\Phi_{\widetilde{P}_R}$ respectively.  This isomorphism of functors induces the restricted isomorphism of functors 
$$\Phi_{\widetilde{P}_R}\circ\Phi_{\widetilde{P}}\cong {\rm id}_{{\rm Dcoh}(X_1,W_1)}.$$ 
Since $\Phi_{\widetilde{P}}\dashv \Phi_{\widetilde{P}_R}$ by Proposition \ref{integral adjoint}, this isomorphism implies that the functor $\Phi_{\widetilde{P}}:{\rm Dcoh}(X_1,W_1)\rightarrow{\rm Dcoh}(X_2,W_2)$ is fully faithful by \cite[Lemma 4.6]{hirano}.

If the integral functor $\Phi_{j_*(P)}:{\rm D^b}({\rm coh}X_1)\rightarrow{\rm D^b}({\rm coh}X_2)$ is an equivalence, its left adjoint functor $\Phi_{j_*(P)_L}$ is fully faithful. Hence, by the above argument, we see that a left adjoint functor $\Phi_{\widetilde{P}_L}:{\rm Dcoh}(X_2,W_2)\rightarrow{\rm Dcoh}(X_1,W_1)$ of the fully faithful functor $\Phi_{\widetilde{P}}$ is also fully faithful. Hence $\Phi_{\widetilde{P}}$ is an equivalence. 
\end{proof}
\end{thm}


 \vspace{2mm}
\subsection{Applications}
 In this last section, we give two applications  of the main theorem.
 
 \subsubsection{Flops of three folds} Let $X$ and $X^{+}$ be smooth quasi-projective threefolds, and let the diagram $$X\xrightarrow{f} Y\xleftarrow{f^+} X^{+}$$ be a flop. Set $Z:=X\times_{Y}X^{+}$ and let $\iota:Z\rightarrow X\times X^{+}$ be the embedding. 

In \cite{bridgeland}, Bridgeland shows the following theorem:

\vspace{2mm}
 \begin{thm}[\cite{bridgeland}]
 The integral functor $$\Phi_{\iota_*(\mathcal{O}_Z)}:{\rm D^b}({\rm coh}X)\rightarrow{\rm D^b}({\rm coh}X^{+})$$
is an equivalence.

 \end{thm}

 Let $G$ be a reductive affine algebraic group acting on $X$, $X^{+}$ and $Y$ with the morphisms $f$ and $f^+$ equivariant. 
 Take a semi invariant regular function $W_Y:Y\rightarrow\mathbb{A}^1$, and set $W:=f^*W_Y$ and $W^{{+}}:=f^{+*}W_Y$. Consider the following cartesian square;
 \[\xymatrix{
&&X\times_{\mathbb{A}^1} X^{+}\ar[lld]_{}\ar[rrd]^{}&&\\
X\ar[rrd]_{W}&&&&X^{+}\ar[lld]^{W^{+}}\\
&&\mathbb{A}^1&&
}\]
The embedding $\iota:Z\rightarrow X\times X^{+}$ factors through $X\times_{\mathbb{A}^1} X^{+}$, i.e.  $\iota$ is the  composition of embeddings $i:Z\rightarrow X\times_{\mathbb{A}^1} X^{+}$ and $j:X\times_{\mathbb{A}^1} X^{+}\rightarrow X\times X^{+}$. Set $$P:=i_*(\mathcal{O}_Z)\in{\rm D^b}({\rm coh}X\times_{\mathbb{A}^1} X^{+}).$$
Since flopping contractions $f$ and $f^{+}$ are proper morphisms, the support of $P$ is proper over $X$ and $X^{+}$. Furthermore, the object $\mathcal{O}_Z\in{\rm D^b}({\rm coh}Z)$ has a tautological $G$-equivariant structure. Hence, $P$ is a $G$-linearizable object. Consequently, we obtain the following  corollary of Theorem \ref{main theorem 2}:

\vspace{2mm}
\begin{cor}\label{flop equivalence}
We have an equivalence of derived factorization categories;
$$\Phi_{\widetilde{P}_G}:{\rm Dcoh}_G(X,W)\xrightarrow{\sim}{\rm Dcoh}_G(X^{+},W^{+}).$$
\end{cor}

\vspace{2mm}

We define $K$-equivalence of gauged LG models. The above gauged LG models $(X,W)^G$ and $(X^{+},W^{+})^G$ are $K$-equivalent.
\vspace{2mm}
\begin{dfn} Let $X_1$ and $X_2$ be smooth varieties with  group $G$-actions, and let $W_1:X_1\rightarrow\mathbb{A}^1$ and $W_2:X_2\rightarrow\mathbb{A}^1$ be $\chi$-semi invariant regular functions for some character $\chi:G\rightarrow\mathbb{G}_m$. The gauged LG models $(X_1,\mathcal{O}(\chi),W_1)^G$ and $(X_2,\mathcal{O}(\chi),W_2)^G$ are called {\it K}\textbf{-equivalent}, if there exists a common $G$-equivariant resolution of $X_1$ and $X_2$
\[\xymatrix{
&&Z\ar[lld]_{p}\ar[rrd]^{q}&&\\
X_1&&&&X_2
}\]
such that $p^*W_1=q^*W_2$ and $p^*\omega_{X_1}\cong q^*\omega_{X_2}$.
\end{dfn} 

By Corollary \ref{flop equivalence} or \cite[Conjecture 2,15]{segal}, it is natural to expect the following conjecture:

\vspace{2mm}
\begin{conj}\label{k implies d}
If two gauged LG models $(X,\mathcal{O}(\chi),W_X)^G$ and $(Y,\mathcal{O}(\chi),W_Y)^G$ are $K$-equivalent, then their derived factorization categories are equivalent;
$${\rm Dcoh}_G(X,W_X)\cong{\rm Dcoh}_G(Y,W_Y).$$
\end{conj}

 The above conjecture for $K$-equivalent gauged LG models of trivial $\sigma$-type is proposed  by Kawamata \cite{kawamata}. The converse of the above conjecture is not true in general. A counterexample to the converse of the Kawamata's conjecture  is given by Uehara \cite{uehara}.

 \vspace{3mm}
 \subsubsection{Equivariantizations of derived equivalences}
 Let $G$ be a reductive affine algebraic group, and let  $X_1$ and $X_2$ be smooth quasi-projective varieties with $G$-actions. 
 
 \vspace{2mm}
 \begin{cor}\label{last}
 Let $P\in {\rm D^b}({\rm coh}X_1\times X_2)$ be an object. Assume that  $P$ is $G$-linearizable object and  the support of $P$ is proper over $X_1$ and $X_2$. Choose an object $P_G\in {\rm D^b}({\rm coh}_GX_1\times X_2)$ such that $\Pi(P_G)\cong P$, where $\Pi:{\rm D^b}({\rm coh}_GX_1\times X_2)\rightarrow{\rm D^b}({\rm coh}X_1\times X_2)$ is the forgetful functor.
 If the integral functor $\Phi_P:{\rm D^b}({\rm coh}X_1)\rightarrow{\rm D^b}({\rm coh}X_1)$ is an equivalence (resp. fully faithful), then the integral functor  $$\Phi_{P_G}:{\rm D^b}({\rm coh}_GX_1)\rightarrow{\rm D^b}({\rm coh}_GX_2)$$
 is also an equivalence (resp. fully faithful).
 
 \begin{proof}
 Extend the $G$-action to $G\times \mathbb{G}_m$-action  by  $\mathbb{G}_m$ acting trivially. Then $P$ is $G\times \mathbb{G}_m$-linearizable.  By Theorem \ref{main theorem 2}, there is an object $\widetilde{P}_{G\times \mathbb{G}_m}\in{\rm Dcoh}_{G\times \mathbb{G}_m}(X_1\times X_2,0)$ which induces an equivalence (resp. fully faithful)
 $$\Phi_{\widetilde{P}_{G\times \mathbb{G}_m}}:{\rm Dcoh}_{{G\times \mathbb{G}_m}}(X_1,0)\rightarrow{\rm Dcoh}_{{G\times \mathbb{G}_m}}(X_2,0).$$
 By Proposition \ref{sigma LG} and equivalences ${\rm coh}_GX_i\cong{\rm coh}[X_i/G]$ for each $i=1,2$, we have equivalences 
 $$\Omega_i:{\rm Dcoh}_{{G\times \mathbb{G}_m}}(X_i,0)\cong{\rm D^b}({\rm coh}_GX_i).$$
 Since the following diagram 
  \[\xymatrix{
{\rm Dcoh}_{{G\times \mathbb{G}_m}}(X_1,0)\ar[rr]^{\Phi_{\widetilde{P}_{G\times \mathbb{G}_m}}}\ar[d]_{\Omega_1}&&{\rm Dcoh}_{{G\times \mathbb{G}_m}}(X_2,0)\ar[d]^{\Omega_2}\\
{\rm D^b}({\rm coh}_GX_1)\ar[rr]^{\Phi_{P_G}}&&{\rm D^b}({\rm coh}_GX_2)
}\]
is commutative, the integral functor $\Phi_{P_G}$ is also an equivalence (resp. fully faithful).
 \end{proof}
 
 \end{cor}
 
 Corollary \ref{last} is shown if the group $G$ is finite by Ploog \cite[Lemma 5]{ploog}. We can also prove Corollary \ref{last} for finite group actions by the result of \cite{elagin2}. 
 
\vspace{4mm}
\address{Department of Mathematics and Information Sciences, Tokyo Metropolitan University, 1-1 Minamiohsawa, Hachioji-shi, Tokyo, 192-0397, Japan}\\
{\it E-mail address}: \email{yuki-hirano@ed.tmu.ac.jp}

\end{document}